\documentclass[11pt]{article}
\usepackage[a4paper,vmargin={3.5cm,3.5cm},hmargin={2.5cm,2.5cm}]{geometry}
\usepackage[T1]{fontenc}
\usepackage[utf8]{inputenc}

\usepackage{hyperref,dsfont}
\hypersetup{colorlinks=true}

\usepackage{mathtools}
\usepackage{graphicx}
\usepackage{framed}
\usepackage[normalem]{ulem}
\usepackage{amsmath,amsthm}
\usepackage{physics}
\usepackage{tikz}
\usepackage{mathdots}
\usepackage{yhmath}
\usepackage{cancel}
\usepackage{color}
\usepackage{siunitx}
\usepackage{array}
\usepackage{multirow}
\usepackage{enumitem}   
\DeclareMathOperator*{\argmax}{argmax}

\DeclareMathOperator*{\argmaxx}{\overset{[2]}{argmax}  }
\usepackage{amsthm}
\usepackage{amssymb}
\usepackage{gensymb}
\usepackage{tabularx}
\usepackage{extarrows}
\usepackage{booktabs}
\usetikzlibrary{fadings}
\usetikzlibrary{patterns}
\usetikzlibrary{shadows.blur}
\usetikzlibrary{shapes}
\usepackage{bbm}
\usepackage{mathtools}
\usepackage{forest}
\forestset{multiple directions/.style={for tree={#1}, phantom, for relative level=1{no edge, delay={!c.content/.pgfmath=content("!u")}, before computing xy={l=0,s=0}}},
    multiple directions/.default={},
    grow subtree/.style={for tree={grow=#1}}, 
    grow' subtree/.style={for tree={grow'=#1}}}
\usepackage{tikz}
\usepackage{capt-of}
\newtheorem{theorem}{Theorem}
\newtheorem{coro}{Corollary}
\newtheorem{prop}{Proposition}[section]
\newtheorem{lemma}[prop]{Lemma}

\newtheorem{definition}[prop]{Definition}
\newtheorem*{theorem*}{Theorem}
\newtheoremstyle{example}{}{}{}{}{\bfseries}{\smallskip}{\newline}{}

\usepackage[maxbibnames=99]{biblatex}
\bibliography{References.bib}
\DeclareMathOperator{\perfE}{perf_E}
\DeclareMathOperator{\perf}{perf}
\DeclareMathOperator{\perfV}{perf_V}
\begin{document}
\title{Optimal unimodular matchings}
\author{Nathanaël Enriquez, Mike Liu, Laurent Ménard and Vianney Perchet}

\theoremstyle{remark}
\newtheorem*{remark}{Remark}

\maketitle
\begin{abstract}
    We consider sequences of finite weighted random graphs that converge locally to unimodular i.i.d. weighted random trees. When the weights are atomless, we prove that the matchings of maximal weight converge locally to a matching on the limiting tree. For this purpose, we introduce and study unimodular matchings on weighted unimodular random trees as well as a notion of optimality for these objects. In this context, we prove that, in law, there is a  unique optimal unimodular matching for a given unimodular tree. We then prove that this law is the local limit of the sequence of matchings of maximal weight. Along the way, we also show that this law is characterised by an equation derived from a message-passing algorithm.
\end{abstract}

\paragraph{Keywords:} Optimal matchings, sparse random graphs, unimodularity, local convergence.

\paragraph{Mathematics Subjects Classification:} 05C70, 05C82, 60C05, 60K35.

\tableofcontents

\section{Introduction}

Optimal matchings in random graphs is a problem with a long history dating back to the paper by Karp and Sipser \cite{4568355}.
In this seminal paper, they obtain the asymptotic size of a maximal matching in a sparse Erd\H{o}s-R\'enyi random graph with average degree $c \in (0,+\infty)$. Their method relies on the analysis of a greedy algorithm that constructs a matching in a graph. It turns out that the algorithm is simple when $c \leq e$ and less straightforward when $c > e$, making the analysis trickier in this regime.

More recently, Gamarnik, Nowick and Swirscsz \cite{gamarnik2003maximum} studied the edge-weighted version of this maximal matching problem on the Erd\H{o}s-R\'enyi random graph and on regular random graphs. Their approach does not rely on the Karp--Sipser algorithm, but on a local optimality ansatz inspired by the earlier introduction of the so-called objective method by Aldous~\cite{Aldous2004,aldous2018processes,Aldous1992AsymptoticsIT,aldous2000zeta2}. Gamarnik, Nowick and Swirscsz were able to compute the asymptotic weight of the optimal matching in Erd\H{o}s-R\'enyi random graph when $c \leq e$ and in the special case where the weights have exponential distribution.

The fundamental tool for studying these kinds of problems is the belief propagation algorithm. Informally, this algorithm assigns numbers to all edges of the graph, with an explicit dependency between neighbours, hence the name ``propagation''. Moreover, they can be interpreted as beliefs: the lower the number, the more likely an edge belongs to the optimal matching. In homogeneous random graphs, the above dependency gives a fixed point equation in law for the beliefs, called the ``belief propagation equation''. 
The restriction $c \leq e$ in \cite{gamarnik2003maximum} comes from the fact that the employed method takes as input the uniqueness of solutions of iterates of the belief propagation equation.

One of the principles of Aldous's objective method is that some asymptotic results are manifestations of limiting objects. In many settings involving random graphs, the objects are limits in law of the graphs for the local topology. In the special case where graphs are rooted uniformly, the graphs and their limits are called unimodular. This line of research has found tremendous success since the seminal paper by Benjamini and Schramm \cite{BenjaminiSchramm1}.

In this article, we study optimal matchings on weighted random graphs that converge to a random tree with iid weights in the Benjamin--Schramm sense. Note that we only require the weights to be independent in the limit. This includes classical models of random graphs such as Erd\H{o}s-R\'enyi random graphs, regular random graphs or configuration models that all converge to unimodular versions of Bienaymé Galton Watson random trees. See for instance the books by van der Hofstad \cite{vanderHofstad1,vanderHofstadBook2}, Benjamini Lyons and Schramm \cite{BENJAMINI_LYONS_SCHRAMM_2015}, or Aldous and Lyons \cite{aldous2018processes}, and references therein.

Strikingly, in this setting and under mild assumptions, the uniqueness of the solution to the iterates of the belief propagation equation mentioned above is not required for the convergence in law of the optimal matching. Informally, the limit $\mathbb M$ is a random matching on the unimodular limiting tree $\mathbb T$ rooted at a distinguished oriented edge $\overset{\rightarrow}{o} = (o_-,o_+) \in \overset{\rightarrow}{E}(\mathbb{T})$ and with i.i.d. edge-weights $w=(w(e))_{e\in E(\mathbb T)}$ with law $\omega$, whose law is invariant by re-rooting. We call such random matchings unimodular (see Section~\ref{sec:defs} for precise definitions). In addition, the limit $\mathbb M$ is optimal in the sense that the quantity
\begin{align} \label{eq:perf}
\perf(\mathbb{T},\overset{\rightarrow}{o}, \mathbb{M}) &:= \mathbb{E}\left[ w(\overset{\rightarrow}{o})\mathbbm{1}_{\overset{\rightarrow}{o} \in \mathbb{M}} \right]
\end{align}
is maximal among all unimodular matchings of $(\mathbb T, \overset{\rightarrow}{o}, w)$. Note that for finite graphs, when $\overset{\rightarrow}{o}$ is chosen uniformly among directed edges, this performance is maximal when the total weight of the matching is maximal.

Before stating our main result, we need to define more precisely Unimodular Bienaymé Galton Watson random trees (UBGW). Let $\pi$ be a distribution on non-negative integers with finite variance and generating function $\phi$. We denote by $\hat{\pi}$ the size-biased version of $\pi$, whose generating function is $\hat \phi$, defined by $\hat \phi(x) = \frac{\phi'(x)}{\phi'(1)}$. We say that a rooted tree $(\mathbb{T},\overset{\rightarrow}{o})$ is a UBGW with offspring distribution $\pi$ if the law of the tree can be obtained as follows. Take $(\mathbb{T}_-,o_-)$ and $(\mathbb{T}_+,o_+)$ two vertex-rooted independent Bienaymé-Galton-Watson trees with offspring distribution $\hat{\pi}$ (that is, each vertex of these trees have an independent number of children distributed according to $\hat \pi$) . The tree $(\mathbb{T},\overset{\rightarrow}{o})$ is then obtained by connecting $(\mathbb{T}_-,o_-)$ and $(\mathbb{T}_+,o_+)$ with the additional edge  $\overset{\rightarrow}{o} = (o_-,o_+)$.
We can now state our first main result:
\begin{theorem}[Existence and uniqueness of the optimal matching for iid weighted UBGW]\label{unicityUBGW}

Let $(\mathbb{T},\overset{\rightarrow}{o},w)$ be an edge-rooted Unimodular Bienaymé-Galton-Watson (UBGW) tree  with reproduction law $\pi$ and i.i.d edge-weights with law $\omega$ such that:
    \begin{itemize}
        \item both $\pi$ and $\omega$ have finite expectation,
        \item the law $\omega$ is atomless.
    \end{itemize}
    There exists a unique (in law) optimal matching $\mathbb{M}_{\mathrm{opt}}(\mathbb{T})$ on $(\mathbb{T},\overset{\rightarrow}{o},w)$. Furthermore, the distribution of $(\mathbb{T},\overset{\rightarrow}{o},w,\mathbb{M}_{\mathrm{opt}}(\mathbb{T}))$ is explicitly described in Proposition~\ref{prop:constructionMZ_h}~(ii).
\end{theorem}

The description of $\mathbb{M}_{\mathrm{opt}}(\mathbb{T})$ given in Proposition~\ref{prop:constructionMZ_h}~(ii) relies on the belief propagation algorithm mentioned above.
In the simple case of a finite deterministic tree,  the algorithm assigns to any  \textit{oriented} edge $(u,v)$ a belief variable $Z(u,v)$ that captures the penalty of excluding $v$ in the optimal matching of $T_{(u,v)}$, where $T_{(u,v)}$ is the connected component of $T \setminus \{ \{u,v \}\}$ containing $v$. When $T_{(u,v)}$ is finite, this corresponds to the difference between the maximal weight of a matching of 
$T_{(u,v)}$ and the maximal weight of a matching of $T_{(u,v)} \setminus \{v\}$. From this, we can infer (see Section~\ref{sec:matchUBGW} for details) the so-called decision rule of the belief propagation procedure:
\begin{align} \label{eq:decisionrule}
      \{ u,v \} \in \mathbb{M}_{\mathrm{opt}}(\mathbb{T}) \Leftrightarrow  Z(u,v)+Z(v,u) < w(u,v) .
\end{align}
We will also see that the tree structure induces the following recursive equation for the belief variables:
\begin{equation}\label{eq:Z_hrecursion}
Z(u,v)=\max \left( 0 \, , \, \max_{\substack{u' \sim v \\ u' \neq u}} \left(w(v,u')-Z(v,u') \right) \right) . 
\end{equation}

On an infinite random tree, such variables can still be defined, and this will be done in detail in Section~\ref{sec:message}.
In a UBGW tree with reproduction law $\pi$ and i.i.d weights of law $\omega$, the stationary distribution satisfies the following equation in law:
\begin{equation}\label{eq:Zlaw}
Z \overset{\text{(law)}}{=} \max \left( 0,\max_{1 \leq i \leq N} \left( w_i - Z_i \right) \right),
\end{equation}
where $N \sim \hat{\pi}$, the $w_i$'s and the $Z_i$'s are all independent, $w_i \sim \omega$ and $Z_i \sim Z$.

As we shall see, it is not very hard to show that this equation in law always has a solution, but the question of its uniqueness is trickier.  It is not hard to see that any solution to Equation~\eqref{eq:Zlaw} gives with the decision rule of Equation~\eqref{eq:decisionrule} a unimodular matching of the underlying tree. Surprisingly, it turns out that the matching associated to any solution of  Equation~\eqref{eq:Zlaw} is optimal.

\bigskip

Our second main result states the local convergence of maximal matchings on random graphs:



\begin{theorem}[Local convergence of maximal matchings]\label{maintheorem}
    Let $G_n=(V_n, E_n,w_n)$ be a sequence of weighted random graphs and let $\mathbb{M}_{\mathrm{opt}}(G_n)$ be any maximum matching on $G_n$. Choose a root $\overset{\rightarrow}{o}_n$ uniformly among the directed edges and suppose that the sequence of rooted weighted random graphs $(G_n,\overset{\rightarrow}{o}_n,w_n)$ converges locally to a vertex-rooted Unimodular Bienaymé-Galton-Watson (UBGW) tree $(\mathbb{T},\overset{\rightarrow}{o},w)$ with reproduction law $\pi$ and i.i.d edge-weights with law $\omega$ such that:
    \begin{itemize}
        \item both $\pi$ and $\omega$ have finite expectation,
        \item the law $\omega$ is atomless,
        \item the ratio $\frac{2|E_n|}{|V_n|}$ converges in probability to the expectation of $\pi$.
    \end{itemize}
    Then $(G_n,\overset{\rightarrow}{o}_n,w_n,\mathbb{M}_{\mathrm{opt}}(G_n))$ converges in law to $(\mathbb{T},\overset{\rightarrow}{o},w,\mathbb{M}_{\mathrm{opt}}(\mathbb{T}))$ for the local topology.
\end{theorem}




\bigskip
Theorem~\ref{maintheorem} describes the full local geometry of the optimal matching. In particular, it can be used to obtain the asymptotic weight and density of the optimal matching on the sequence of graphs:

\begin{coro}\label{coro:perf}
Let $G_n=(V_n, E_n,w_n,\overset{\rightarrow}{o}_n)$ be a sequence of uniformly rooted weighted graphs as in Theorem~\ref{maintheorem}. Let $(\mathbb{T},w,\overset{\rightarrow}{o})$ be the local limit of the previous sequence. Take Z and Z' to be independent copies of random variables satisfying~\eqref{eq:Zlaw} and $W$ with law $\omega$ independent of $Z$ and $Z'$. One has
\begin{align}
        &\lim_{n \rightarrow \infty} \mathbb{E}\left[ \frac{\sum_{e \in E(\mathbb{M}_{\mathrm{opt}}(G_n))}w(e) }{|E_n|} \right] = \mathbb{E}\left[w(\overset{\rightarrow}{o})\mathbbm{1}_{ \overset{\rightarrow}{o} \in \mathbb{M}_{\mathrm{opt}}(\mathbb{T})}\right] = \mathbb{E}\left[W\mathbbm{1}_{Z+Z'<W}  \right],\\
         &\lim_{n \rightarrow \infty} \mathbb{E}\left[ \frac{\left| E \left( \mathbb{M}_{\mathrm{opt}}(G_n) \right) \right|}{|E_n|} \right] = \mathbb{P}\left(\overset{\rightarrow}{o} \in \mathbb{M}_{\mathrm{opt}(\mathbb{T})}\right) = \mathbb{P}\left(Z+Z'<W\right) = \frac{1-\phi(\hat{\phi}^{-1}(\mathbb P (Z=0)))}{\phi'(1)}  . \label{eq:unifedge}
    \end{align}
\end{coro}

Notice that \eqref{eq:unifedge} is the probability that a uniform edge is in the optimal matching. It could also be stated as the probability that a uniform vertex is in the matching:
\[
         \lim_{n \rightarrow \infty} \mathbb{E}\left[ \frac{\left| V \left( \mathbb{M}_{\mathrm{opt}}(G_n) \right) \right|}{|V_n|} \right]  = 1-\phi \left( \hat{\phi}^{-1}(\mathbb P (Z=0)) \right).
\]
We can also establish a conditional version of this statement:
\begin{coro}
Let $G_n=(V_n, E_n,w_n)$ be a sequence of weighted graphs as in Theorem~\ref{maintheorem}. Let $o_n$ be a uniform vertex of $G_n$. One has
\begin{equation}
    \lim_{n \rightarrow \infty} \mathbb{P}\left( o_n \in V \left( \mathbb{M}_{\mathrm{opt}}(G_n) \right) \middle| \deg(o_n)=k\right) = 1- \left(\hat{\phi}^{-1}(\mathbb P (Z=0))\right)^{k}.
\end{equation}
\end{coro}

\bigskip

With the aid of Theorem~\ref{unicityUBGW}, we are able to give a probabilistic proof of the uniqueness of solutions to Equation~\eqref{eq:Zlaw}:
\begin{theorem}[Uniqueness of the belief propagation stationary measure]\label{th:uniqueness}
Let $\omega$ be a non atomic law that is not entirely singular with respect to Lebesgue measure in the neighbourhood of $\sup(\mathrm{supp}(\omega))$. 
Let $\hat{\pi}$ be a law on natural integers with finite expectation. 
Then there is a unique law $\zeta$ satisfying:

\begin{equation*}
Z \overset{\text{(law)}}{=} \max \left( 0,\max_{1 \leq i \leq N} \left( w_i - Z_i \right) \right),
\end{equation*}
where $N \sim \hat{\pi}$, the $w_i$'s and the $Z_i$'s are all independent, $w_i \sim \omega$ and $Z_i,Z \sim \zeta$.

\end{theorem}
As mentioned, we are not able to prove it directly, however, we are able to prove that for any matching distribution $\mathbb{M}$, there is at most one solution to Equation~\eqref{eq:Zlaw} that can yield $\mathbb{M}$ through the rule~\eqref{eq:decisionrule}. Combined with the uniqueness of the optimal matching given by Theorem~\ref{unicityUBGW}, we are then able to prove that there is only one solution to Equation~\eqref{eq:Zlaw}. Conversely, we emphasise that assuming uniqueness of the solution of Equation~\eqref{eq:Zlaw} does not help in proving Theorem~\ref{unicityUBGW}.

\bigskip

Let us conclude this introduction by mentioning a few papers that deal with matchings and Benjamini-Schramm convergence. Usually, the setting in these works consist of unweighted random graphs that converge for the topology and they prove the convergence of some matching statistics. Among these are the size of the maximal matchings \cite{4568355,bordenave2012matchings} or the number of maximal matchings and the matching polynomial \cite{matchingmeasure1,matchingmeasure2}.

\paragraph{Organisation of the paper}

In Section~\ref{sec:defs}, we recall basic definitions on unimodularity and local convergence of random graphs. We then introduce the notions of unimodular matchings and optimality for unimodular matchings. 

Section~\ref{sec:matchUBGW} carries out the study of optimal matchings on UBGW trees and contains the proof of Theorem~\ref{unicityUBGW}. Our construction of unimodular matchings $\mathbb M$ is based on the message passing algorithm and inspired by Aldous' work \cite{aldous2000zeta2} on the Poisson Weighted Infinite Tree (PWIT). However, we work in the more general context of unimodular random trees, and one of the main difficulties is that, contrary to Aldous' work, we deal with partial matchings rather than perfect matchings. Furthermore, there is no explicit solution to Equation~\eqref{eq:Zlaw}. We overcome these two difficulties by introducing self-loops in Subsection~\ref{sec:selfloops} and by proving just enough properties on solutions to Equation~\eqref{eq:Zlaw}. We then prove in Subsection~\ref{sec:optimality} that the unimodular matching $\mathbb M$ we constructed is in fact optimal. Finally, in Subsection~\ref{sec:uniqueness}, we prove that the law of $\mathbb M$ is the only optimal matching law and does not depend on the particular solution of Equation~\eqref{eq:Zlaw} used.
This bypasses previous restrictions on the belief propagation approach as in Gamarnik, Nowick and Swirscsz ~\cite{gamarnik2003maximum}, Aldous and Steele~\cite{Aldous2004} or Aldous and Bandyopadhyay~\cite{Aldous_2005} that mostly required the variables $Z(u,v)$ to be measurable with respect to the tree.

The main idea behind Theorem~\ref{maintheorem} is to construct almost optimal matchings of finite graphs from an optimal matching $\mathbb{M}_{\mathrm{opt}}(\mathbb{T})$ on their limiting tree. This is the purpose of section~\ref{sec:infinitetofinite}. A similar program has already been carried out by Aldous in the special case of the random assignment problem on bipartite graphs \cite{Aldous1992AsymptoticsIT}.

In section~\ref{sec:thuniqueness}, we deduce Theorem~\ref{th:uniqueness} by studying path conditionings on $(\mathbb{T},o,w,\mathbb{M}_{\mathrm
opt}(\mathbb{T}))$, in particular we show that if $\mathbb{M}_\mathrm{opt}(\mathbb{T})$ is given by a family of variables $(Z(u,v))_{(u,v) \in \overset{
\rightarrow}{E}}$ through the rule~\eqref{eq:decisionrule}, then the conditional law of the matching with respect to given weights on paths is sufficient to recover the law of the variables $(Z(u,v))_{(u,v) \in \overset{
\rightarrow}{E}}$.

In section~\ref{sec:applications}, we show how to apply Theorem~\ref{maintheorem} to obtain the announced corollaries. Finally, in Section~\ref{sec:extensions}, we discuss a few generalisations of the method to adjacent problems such as when the weights are on vertices instead of edges.

\paragraph{Acknowledgements.}
L.M.'s research is supported by the ANR grant ProGraM (ANR-19-CE40-0025). V.P. acknowledges the support of the ANR through the grant DOOM (ANR-23-CE23-0002) and through the PEPR IA FOUNDRY project (ANR-23-PEIA-0003).

\section{Unimodular matched random graphs} \label{sec:defs}
The goal of this section is to precisely define our objects of interest. The key concepts are local topology and unimodularity of rooted graphs. Depending on the context, it is often easier to work with graphs rooted either at a vertex or at an oriented edge. We will present both cases, as well as the classical correspondence between the two points of view. For a more complete exposition on the subject, we refer to \cite{aldous2018processes}.

Finally, we will define unimodular matchings and a concept of optimality with respect to a natural performance for vertex-rooted or edge-rooted graphs. We will see that the correspondence mentioned above conserves our notion of optimality.

\subsection{Rooted graphs and local convergence}

We will work on rooted weighted graphs defined as follows.
\begin{definition}
     Let $G=(V,E)$ be a (locally finite) graph.
     We call \emph{vertex-rooted weighted graph} (resp. \emph{edge-rooted weighted graph}), the triplet $(G,o,w)$, where $o \in V$ (resp. $o \in \overset{\rightarrow}{E}$) and $w$ is a function from $E$ to $\mathbb{R}$.
\end{definition}
It will be useful to add graph decorations, namely functions that map  $\overset{\rightarrow }{E}$ to $\overline{\mathbb{R}}$.
\begin{definition}
    Let $(G,o,w)$ be a vertex-rooted (resp.\ edge-rooted) weighted graph. Fix $I$ an integer and $(f_i)_{i \in \{1,...,I\}}$ some functions from $\overset{\rightarrow}{E}$ to $\overline{\mathbb{R}}$. We then say that
    $(G,o,w,(f_i)_{i \in \{1,...,I\}})$ is a \emph{decorated} vertex-rooted (resp.\ edge-rooted) weighted graph.  If $f_i$ is symmetrical then we will identify it with a map from $E$ to $\overline{\mathbb{R}}$.
\end{definition}

Since we are not interested in the labels of the vertices, we will work up to graph isomorphism.
We will say that $(G,o,w,(f_i)_{i \in \{1,...,I\}}) \simeq (G',o',w',(f_i')_{i \in \{1,...,I\}})$ if there exists some one-to-one function $g$ called graph isomorphism from $V$ to $V'$ such that if $g(o)=o'$ and for all $(u,v) \in V$, $w'(g(u),g(v))=w(u,v)$ and for all $ i \in \{1,...,I\}$, $f'_i(g(u),g(v))=f_i(u,v)$ .

\begin{definition}
Let $\mathcal{G}^{\star}$ be the space of locally finite decorated vertex-rooted weighted graphs up to isomorphism.
We will write $\mathcal{L}(\mathcal{G}^{\star})$ for the space of laws on this space. Similarly, we denote by  $\hat{\mathcal{G}}^{\star}$ and $\mathcal{L}(\hat{\mathcal{G}}^{\star})$ the corresponding edge-rooted space and laws.
\end{definition}
When we do not need to keep track of the weights, root, or some decorations of the graphs, we will denote elements of 
$\mathcal{G}^{\star}$ indifferently by $G$, $(G,o)$, $(G,o,w)$, $(G,o,w,(f_i)_{i \in \{1,...,I\}})$, ..., keeping only the quantities we are currently interested in.
\bigskip

The topology of local convergence on $\mathcal{G}^{\star}$ and $\hat{\mathcal{G}^{\star}}$ was first introduced by Benjamini-Schramm~\cite{BenjaminiSchramm1} and by Aldous-Steele \cite{Aldous2004}, whose  precise definitions follow. We denote by $d_{gr}$ the graph distance between vertices of a graph. For any vertex-rooted decorated weighted graph $(G=(V,E),o,w,(f_i)_{i \in I})$ and $H>0$, the $H$-neighborhood $B_H((G,o,w,(f_i)_{i \in \{1,...,I\}}))$ of $o$ in the graph is the vertex-rooted decorated weighted graph $(G_H,o,w_H,(f_{i,H})_{i \in \{1,...,I\}})$ with $G_H :=(V_H,E_H)$ such that
\begin{align*}
V_H &=\{v \in V, d_{gr}(o,v)\leq H\},\\
E_H &=\{ (u,v) \in E, (u,v)\in V_H^2 \},
\end{align*}
$w_H$ are the weights of edges in $E_H$ and, for all $i \in \{1,...,I\}$, the decoration $f_{i,H}$ is the decoration $f_{i}$ restricted to $\overset{\rightarrow}{E}_H$.


Let $(G,o,w,(f_{i})_{i \in \{1,...,I\}})$ and $(G',o',w',(f_i')_{i \in \{1,...,I\}}$ be two vertex-rooted decorated weighted graphs. Let $H \geq 0$ be the largest integer such that there exists a graph isomorphism $g$ from  $B_H(G,o)$ to $ B_H (G',o')$ such that $\| w'_H\circ g - w_H \|_{\infty} \leq \frac{1}{H}$ and $\| f'_{i,H} \circ g - f_{i,H} \|_{\infty} \leq \frac{1}{H}$ for all $i \in \{1, \ldots , I \}$. We set
\begin{equation} \label{eq:dloc}
   d_{\mathrm{loc}} \left( (G,o,w,(f_{i})_{i \in \{1,...,I\}}) , (G',o',w',(f_i')_{i \in \{1,...,I\}} \right) := \frac{1}{1+H}.
\end{equation}
The function $ d_{\mathrm{loc}}$ is a distance on $\mathcal G^\star$ and the space $(\mathcal{G}^{\star},d_\mathrm
{loc})$ is a Polish space. The topology induced by  $ d_{\mathrm{loc}}$ is called the local topology on $\mathcal{G}^{\star}$. Weak convergence on $\mathcal{L}(\mathcal{G}^{\star})$ for the local topology is called local weak convergence.
Similar definitions can be given for the edge-rooted versions.

\subsection{Vertex-rooted and edge-rooted unimodularity}

We now introduce the notion of unimodularity.
To this end, we define the space of doubly rooted graphs up to isomorphism $\mathcal{G}^{\star\star}$ similarly as before, but with two distinguished roots. It will be simpler to give separately specific definitions for vertex-rooted and edge-rooted graphs.
 We refer to Aldous-Lyons \cite{aldous2018processes} for a comprehensive exposition on the topic.

\begin{definition}[Vertex-rooted unimodularity]
We say that a probability measure $\mu$ on decorated vertex-rooted graphs is unimodular if the following statement holds for every measurable $f : \mathcal{G}^{\star\star} \mapsto \mathbb{R}_+$:
     
    \[
    \int_{\mathcal G^{\star}} \left( \sum_{v \in V(G)} f(G,o,v)  \right) \, 
    \mathrm{d} \mu (G,o)
     = \int_{\mathcal G^{\star}} \left( \sum_{v \in V(G)} f(G,v,o)  \right) \, 
    \mathrm{d} \mu (G,o)   .\]
    The subspace of unimodular laws on vertex rooted-graph will be noted $\mathcal{L}_U(\mathcal{G}^{\star})$
\end{definition}
The definition can be written alternatively as 
\[ \mathbb{E}_{(\mathbb{G},\mathbf{o}) \sim \mu  } \left[\sum_{v \in V(G)} f(\mathbb{G},\mathbf{o},v)\right]=\mathbb{E}_{(\mathbb{G},\mathbf{o}) \sim \mu }\left[\sum_{v \in V(G)} f(\mathbb{G},v,\mathbf{o})\right].    \]
For edge-rooted graphs, we will use the following definition: 
\begin{definition}[Edge-rooted unimodularity]
    Let $\hat{\mu}$ be a probability measure on decorated edge-rooted graphs.
    Let $(\mathbb{G},(\mathbf{o}_-,\mathbf{o}_+))$ be a random edge-weighted decorated graph with law $\hat \mu$. Let $\overset{\rightarrow}{e}_1$ be a uniformly picked directed edge of the form $(\mathbf{o}_+,\mathbf{v})$ for $\mathbf{v} \neq \mathbf{o}_-$. 
    We say that $\hat{\mu}$ is:
    \begin{itemize}
    \item stationary if $(\mathbb{G},(\mathbf{o}_-,\mathbf{o}_+))\overset{\mathcal{L}}{=}(\mathbb{G},\overset{\rightarrow}{e}_1)$
    \item revertible if $(\mathbb{G},(\mathbf{o}_-,\mathbf{o}_+))\overset{\mathcal{L}}{=}(\mathbb{G},(\mathbf{o}_+,\mathbf{o}_-))$ 
    \item unimodular if it is both revertible and stationary.
    \end{itemize}
    The subspace of unimodular laws on edge-rooted graphs will be noted $\mathcal{L}_U(\hat{\mathcal{G}}^{\star})$
\end{definition}

To simplify notation, we will say that a random rooted graph (as a random variable) is unimodular when its corresponding law is. 

\begin{remark}
    Fix $G$ a finite deterministic graph, and let $\mathbf o$ be a random vertex (resp. oriented edge). It is straightforward to check that $(G,\mathbf o)$ is unimodular iff $\mathbf o$ is uniform. Hence, unimodular graphs can be viewed as generalisations of uniformly rooted graphs, which is one of the fundamental findings of \cite{BenjaminiSchramm1}.
\end{remark}

We now present the classical transformation that maps a unimodular vertex-rooted graph to a unimodular edge-rooted graph. Heuristically, to transform a graph rooted at a uniform vertex into a graph rooted at a uniformly oriented edge, one has to pick an oriented edge starting at its root vertex, but this induces a bias by the degree of the root vertex. Indeed, a directed edge $(u,v)$ will be less likely to be the new root the greater the degree of $u$ is. This bias has to be taken into account to conserve unimodularity. This is done in the following transformations on the laws.
\begin{definition} \label{def:R}
    Take $\mu_v \in \mathcal{L}(\mathcal{G}^{\star})$ such that $0<m := \int_{\mathcal G^{\star}} \deg (o) \, \mathrm{d} \mu_v(G,o) < \infty$. Let $R(\mu_v) \in \mathcal{L}(\hat{\mathcal{G}}^{\star})$  be the unique measure such that,
    for every $f$ measurable from $\hat{\mathcal{G}}^{\star}$ to $\mathbb{R}_+$,
\[
\int_{\hat{\mathcal{G}}^{\star}}f(G,(o_-,o_+)) \, \mathrm{d} R(\mu_v)(G,(o_-,o_+)) =
\int_{\mathcal{G}*} 
 \frac{1}{m} 
\sum_{u \sim o} f(G,o,u) \, \mathrm{d} \mu_v(G,o).
\]
\end{definition}
\begin{remark}
The operator $R$ is the composition of two transformations. The first operator consists of choosing a uniformly oriented edge started at the root vertex, giving a measure $R_1(\mu_v) \in \mathcal{L}(\hat{\mathcal{G}}^{\star})$:
\[
\int_{\hat{\mathcal{G}}^{\star}} 
f(G,(o_-,o_+))\, \mathrm{d} R_1(\mu_v)(G,(o_-,o_+)))= 
\int_{\mathcal{G}*} 
 \frac{1}{\deg(o)} 
\sum_{u \sim o} f(\textbf{}G,o,u) \, \mathrm{d} \mu_v(G,o).
\]
The second step is then to cancel the bias by the degree of the vertex, giving $R(\mu_v) = R_2 \circ R_1 (\mu_v) \in \mathcal{L}(\hat{\mathcal{G}}^{\star})$:
    \[ \int_{\hat{\mathcal{G}}^{\star}}f(G,(o_-,o_+))\mathrm{d}R_2 \circ R_1 (\mu_v) (G,(o_-,o_+)) = \int_{\hat{\mathcal{G}}^{\star}} \frac{\deg(o_-)}{m} f(G,(o_-,o_+)) \, \mathrm{d}R_1 (\mu_v)(G,(o_-,o_+)).\]
    
Equivalently, to obtain an edge-rooted version of a vertex rooted random graph, one can consider the size-biased version of the original measure $\mathrm{d}\mu_v(G,o)$ (namely sample a graph with measure
$\frac{\deg(o)}{m}\mathrm{d}\mu_v(G,o)$) and then from this sampled graph, starting from the root, select an edge uniformly at random to get the edge measure $dR(\mu_v)$.
\end{remark}

\bigskip

The following proposition links vertex unimodularity and edge unimodularity.
\begin{prop}[Theorem 4.1 in \cite{aldous2018processes}]
Let $\mu_v \in \mathcal{L}(\mathcal{G}^{\star})$ be such that $0<\int_{\mathcal G^{\star}} \deg (o) \, \mathrm{d} \mu_v(G,o) < \infty$.
Then the measure $\mu_v$ is unimodular if and only if $R(\mu_v)$ is unimodular.
\end{prop}

\subsection{Unimodular Bienaymé-Galton-Watson trees}
In this subsection, we introduce Unimodular Bienaymé-Galton-Watson trees (UBGW) along with models of random graphs that converge locally in law to these trees. We will present both the vertex and the edge-rooted point of view, without weights. In either case, the weighted version with weight law $\omega$ corresponds to drawing $(w(e))_{e \in E}$ independently of law $\omega$.

Let $\pi$ be a probability measure on $\mathbb{Z}_+ = \{0,1, ... \}$ with finite expectation $m > 0$. Let $\hat{\pi}$ be the size-biased version of $\pi$, that is, $\forall k \geq 0$, $\hat{\pi} (k) =\frac{k}{m} \pi (k)$.

In the introduction, we defined the edge-rooted Unimodular Bienaymé Galton Watson Tree with reproduction law $\pi$. Let us recall the construction: take two independent copies of Bienaymé-Galton-Watson trees with offspring distribution  $\hat{\pi}$ with respective root vertex $o_-$ and $o_+$, and connect their roots by the oriented edge $(o_-,o_+)$. The resulting random tree $(\hat{\mathbb T}, (o_-,o_+))$ is an edge-rooted unimodular random graph. See Figure~\ref{fig:EUBGW} for an illustration.
\tikzset{every picture/.style={line width=0.75pt}} 
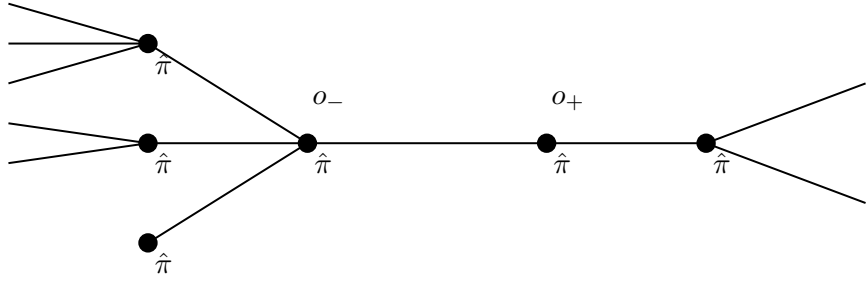
\begin{figure}[!t]
\centering
\begin{tikzpicture}[x=0.75pt,y=0.75pt,yscale=-1,xscale=1]
\draw    (250,140) -- (370,140) ;
\draw    (170,90) -- (250,140) ;
\draw    (170,140) -- (250,140) ;
\draw    (250,140) -- (170,190) ;
\draw    (370,140) -- (450,140) ;
\draw    (530,110) -- (450,140) ;
\draw    (530,170) -- (450,140) ;
\draw  [fill={rgb, 255:red, 0; green, 0; blue, 0 }  ,fill opacity=1 ] (245.65,140) .. controls (245.65,137.6) and (247.6,135.65) .. (250,135.65) .. controls (252.4,135.65) and (254.35,137.6) .. (254.35,140) .. controls (254.35,142.4) and (252.4,144.35) .. (250,144.35) .. controls (247.6,144.35) and (245.65,142.4) .. (245.65,140) -- cycle ;
\draw    (100,70) -- (170,90) ;
\draw    (100,110) -- (170,90) ;
\draw    (100,90) -- (170,90) ;
\draw    (100,130) -- (170,140) ;
\draw    (100,150) -- (170,140) ;
\draw  [fill={rgb, 255:red, 0; green, 0; blue, 0 }  ,fill opacity=1 ] (365.65,140) .. controls (365.65,137.6) and (367.6,135.65) .. (370,135.65) .. controls (372.4,135.65) and (374.35,137.6) .. (374.35,140) .. controls (374.35,142.4) and (372.4,144.35) .. (370,144.35) .. controls (367.6,144.35) and (365.65,142.4) .. (365.65,140) -- cycle ;
\draw  [fill={rgb, 255:red, 0; green, 0; blue, 0 }  ,fill opacity=1 ] (165.65,190) .. controls (165.65,187.6) and (167.6,185.65) .. (170,185.65) .. controls (172.4,185.65) and (174.35,187.6) .. (174.35,190) .. controls (174.35,192.4) and (172.4,194.35) .. (170,194.35) .. controls (167.6,194.35) and (165.65,192.4) .. (165.65,190) -- cycle ;
\draw  [fill={rgb, 255:red, 0; green, 0; blue, 0 }  ,fill opacity=1 ] (165.65,140) .. controls (165.65,137.6) and (167.6,135.65) .. (170,135.65) .. controls (172.4,135.65) and (174.35,137.6) .. (174.35,140) .. controls (174.35,142.4) and (172.4,144.35) .. (170,144.35) .. controls (167.6,144.35) and (165.65,142.4) .. (165.65,140) -- cycle ;
\draw  [fill={rgb, 255:red, 0; green, 0; blue, 0 }  ,fill opacity=1 ] (165.65,90) .. controls (165.65,87.6) and (167.6,85.65) .. (170,85.65) .. controls (172.4,85.65) and (174.35,87.6) .. (174.35,90) .. controls (174.35,92.4) and (172.4,94.35) .. (170,94.35) .. controls (167.6,94.35) and (165.65,92.4) .. (165.65,90) -- cycle ;
\draw  [fill={rgb, 255:red, 0; green, 0; blue, 0 }  ,fill opacity=1 ] (445.65,140) .. controls (445.65,137.6) and (447.6,135.65) .. (450,135.65) .. controls (452.4,135.65) and (454.35,137.6) .. (454.35,140) .. controls (454.35,142.4) and (452.4,144.35) .. (450,144.35) .. controls (447.6,144.35) and (445.65,142.4) .. (445.65,140) -- cycle ;

\draw (251,112.4) node [anchor=north west][inner sep=0.75pt]    {$o_{-}$};
\draw (371,112.4) node [anchor=north west][inner sep=0.75pt]    {$o_{+}$};
\draw (172,193.4) node [anchor=north west][inner sep=0.75pt]    {$\hat{\pi }$};
\draw (172,143.4) node [anchor=north west][inner sep=0.75pt]    {$\hat{\pi }$};
\draw (172,93.4) node [anchor=north west][inner sep=0.75pt]    {$\hat{\pi }$};
\draw (252,143.4) node [anchor=north west][inner sep=0.75pt]    {$\hat{\pi }$};
\draw (372,143.4) node [anchor=north west][inner sep=0.75pt]    {$\hat{\pi }$};
\draw (452,143.4) node [anchor=north west][inner sep=0.75pt]    {$\hat{\pi }$};
\end{tikzpicture}
\caption{A 2-Neighbourhood of an edge-rooted UBGW tree with the law of the number of children drawn on every vertex.} \label{fig:EUBGW}
\end{figure}

The vertex-rooted UBGW tree is the random tree $\mathbb T$ with the following law:
    \begin{itemize}
        \item The number of children of vertices of $\mathbb T$ are all independent.
        \item The number of children of the root $o$ is distributed according to $\pi$.
        \item Every non-root vertex has a number of children distributed according to $\hat \pi$, the sized biased version of $\pi$.
    \end{itemize}
The random tree $(\mathbb T, o)$ is a vertex-rooted unimodular random graph. See Figure~\ref{fig:VUBGW} for an illustration.
\tikzset{every picture/.style={line width=0.75pt}} 
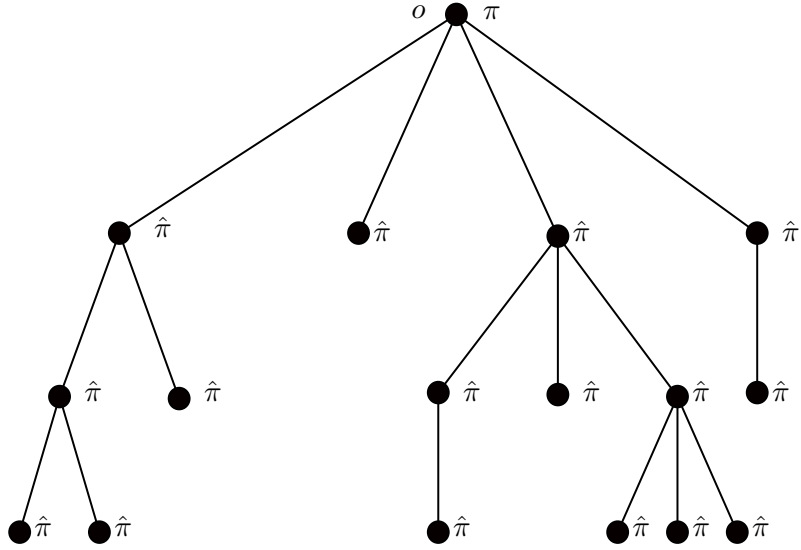
\begin{figure}[!ht]
\centering
\begin{tikzpicture}[x=0.75pt,y=0.75pt,yscale=-1,xscale=1]
\draw    (299.1,0.3) -- (130,110) ;
\draw    (299.1,0.3) -- (250,110) ;
\draw    (299.1,0.3) -- (350,111.5) ;
\draw    (299.1,0.3) -- (450,110) ;
\draw    (130,110) -- (100,192) ;
\draw    (130,110) -- (160,193) ;
\draw    (350,111.5) -- (290,190) ;
\draw    (350,111.5) -- (350,191) ;
\draw    (350,111.5) -- (410,192) ;
\draw    (450,190) -- (450,110) ;
\draw    (100,192) -- (80,260) ;
\draw    (100,192) -- (120,260) ;
\draw    (290,190) -- (290,260) ;
\draw    (410,192) -- (380,260) ;
\draw    (410,192) -- (410,260) ;
\draw    (410,192) -- (440,260) ;
\draw  [fill={rgb, 255:red, 7; green, 0; blue, 0 }  ,fill opacity=1 ] (124.77,110) .. controls (124.77,107.02) and (127.11,104.6) .. (130,104.6) .. controls (132.89,104.6) and (135.23,107.02) .. (135.23,110) .. controls (135.23,112.98) and (132.89,115.4) .. (130,115.4) .. controls (127.11,115.4) and (124.77,112.98) .. (124.77,110) -- cycle ;
\draw  [fill={rgb, 255:red, 7; green, 0; blue, 0 }  ,fill opacity=1 ] (293.88,0.3) .. controls (293.88,-2.68) and (296.21,-5.1) .. (299.1,-5.1) .. controls (301.99,-5.1) and (304.33,-2.68) .. (304.33,0.3) .. controls (304.33,3.28) and (301.99,5.7) .. (299.1,5.7) .. controls (296.21,5.7) and (293.88,3.28) .. (293.88,0.3) -- cycle ;
\draw  [fill={rgb, 255:red, 7; green, 0; blue, 0 }  ,fill opacity=1 ] (244.77,110) .. controls (244.77,107.02) and (247.11,104.6) .. (250,104.6) .. controls (252.89,104.6) and (255.23,107.02) .. (255.23,110) .. controls (255.23,112.98) and (252.89,115.4) .. (250,115.4) .. controls (247.11,115.4) and (244.77,112.98) .. (244.77,110) -- cycle ;
\draw  [fill={rgb, 255:red, 7; green, 0; blue, 0 }  ,fill opacity=1 ] (344.77,111.5) .. controls (344.77,108.52) and (347.11,106.1) .. (350,106.1) .. controls (352.89,106.1) and (355.23,108.52) .. (355.23,111.5) .. controls (355.23,114.48) and (352.89,116.9) .. (350,116.9) .. controls (347.11,116.9) and (344.77,114.48) .. (344.77,111.5) -- cycle ;
\draw  [fill={rgb, 255:red, 7; green, 0; blue, 0 }  ,fill opacity=1 ] (444.77,110) .. controls (444.77,107.02) and (447.11,104.6) .. (450,104.6) .. controls (452.89,104.6) and (455.23,107.02) .. (455.23,110) .. controls (455.23,112.98) and (452.89,115.4) .. (450,115.4) .. controls (447.11,115.4) and (444.77,112.98) .. (444.77,110) -- cycle ;
\draw  [fill={rgb, 255:red, 7; green, 0; blue, 0 }  ,fill opacity=1 ] (94.77,192) .. controls (94.77,189.02) and (97.11,186.6) .. (100,186.6) .. controls (102.89,186.6) and (105.23,189.02) .. (105.23,192) .. controls (105.23,194.98) and (102.89,197.4) .. (100,197.4) .. controls (97.11,197.4) and (94.77,194.98) .. (94.77,192) -- cycle ;
\draw  [fill={rgb, 255:red, 7; green, 0; blue, 0 }  ,fill opacity=1 ] (74.77,260) .. controls (74.77,257.02) and (77.11,254.6) .. (80,254.6) .. controls (82.89,254.6) and (85.23,257.02) .. (85.23,260) .. controls (85.23,262.98) and (82.89,265.4) .. (80,265.4) .. controls (77.11,265.4) and (74.77,262.98) .. (74.77,260) -- cycle ;
\draw  [fill={rgb, 255:red, 7; green, 0; blue, 0 }  ,fill opacity=1 ] (114.77,260) .. controls (114.77,257.02) and (117.11,254.6) .. (120,254.6) .. controls (122.89,254.6) and (125.23,257.02) .. (125.23,260) .. controls (125.23,262.98) and (122.89,265.4) .. (120,265.4) .. controls (117.11,265.4) and (114.77,262.98) .. (114.77,260) -- cycle ;
\draw  [fill={rgb, 255:red, 7; green, 0; blue, 0 }  ,fill opacity=1 ] (154.77,193) .. controls (154.77,190.02) and (157.11,187.6) .. (160,187.6) .. controls (162.89,187.6) and (165.23,190.02) .. (165.23,193) .. controls (165.23,195.98) and (162.89,198.4) .. (160,198.4) .. controls (157.11,198.4) and (154.77,195.98) .. (154.77,193) -- cycle ;
\draw  [fill={rgb, 255:red, 7; green, 0; blue, 0 }  ,fill opacity=1 ] (344.77,191) .. controls (344.77,188.02) and (347.11,185.6) .. (350,185.6) .. controls (352.89,185.6) and (355.23,188.02) .. (355.23,191) .. controls (355.23,193.98) and (352.89,196.4) .. (350,196.4) .. controls (347.11,196.4) and (344.77,193.98) .. (344.77,191) -- cycle ;
\draw  [fill={rgb, 255:red, 7; green, 0; blue, 0 }  ,fill opacity=1 ] (284.77,190) .. controls (284.77,187.02) and (287.11,184.6) .. (290,184.6) .. controls (292.89,184.6) and (295.23,187.02) .. (295.23,190) .. controls (295.23,192.98) and (292.89,195.4) .. (290,195.4) .. controls (287.11,195.4) and (284.77,192.98) .. (284.77,190) -- cycle ;
\draw  [fill={rgb, 255:red, 7; green, 0; blue, 0 }  ,fill opacity=1 ] (404.77,192) .. controls (404.77,189.02) and (407.11,186.6) .. (410,186.6) .. controls (412.89,186.6) and (415.23,189.02) .. (415.23,192) .. controls (415.23,194.98) and (412.89,197.4) .. (410,197.4) .. controls (407.11,197.4) and (404.77,194.98) .. (404.77,192) -- cycle ;
\draw  [fill={rgb, 255:red, 7; green, 0; blue, 0 }  ,fill opacity=1 ] (284.77,260) .. controls (284.77,257.02) and (287.11,254.6) .. (290,254.6) .. controls (292.89,254.6) and (295.23,257.02) .. (295.23,260) .. controls (295.23,262.98) and (292.89,265.4) .. (290,265.4) .. controls (287.11,265.4) and (284.77,262.98) .. (284.77,260) -- cycle ;
\draw  [fill={rgb, 255:red, 7; green, 0; blue, 0 }  ,fill opacity=1 ] (374.77,260) .. controls (374.77,257.02) and (377.11,254.6) .. (380,254.6) .. controls (382.89,254.6) and (385.23,257.02) .. (385.23,260) .. controls (385.23,262.98) and (382.89,265.4) .. (380,265.4) .. controls (377.11,265.4) and (374.77,262.98) .. (374.77,260) -- cycle ;
\draw  [fill={rgb, 255:red, 7; green, 0; blue, 0 }  ,fill opacity=1 ] (404.77,260) .. controls (404.77,257.02) and (407.11,254.6) .. (410,254.6) .. controls (412.89,254.6) and (415.23,257.02) .. (415.23,260) .. controls (415.23,262.98) and (412.89,265.4) .. (410,265.4) .. controls (407.11,265.4) and (404.77,262.98) .. (404.77,260) -- cycle ;
\draw  [fill={rgb, 255:red, 7; green, 0; blue, 0 }  ,fill opacity=1 ] (434.77,260) .. controls (434.77,257.02) and (437.11,254.6) .. (440,254.6) .. controls (442.89,254.6) and (445.23,257.02) .. (445.23,260) .. controls (445.23,262.98) and (442.89,265.4) .. (440,265.4) .. controls (437.11,265.4) and (434.77,262.98) .. (434.77,260) -- cycle ;
\draw  [fill={rgb, 255:red, 7; green, 0; blue, 0 }  ,fill opacity=1 ] (444.77,190) .. controls (444.77,187.02) and (447.11,184.6) .. (450,184.6) .. controls (452.89,184.6) and (455.23,187.02) .. (455.23,190) .. controls (455.23,192.98) and (452.89,195.4) .. (450,195.4) .. controls (447.11,195.4) and (444.77,192.98) .. (444.77,190) -- cycle ;
\draw (311,-5.3) node [anchor=north west][inner sep=0.75pt]    {$\pi $};
\draw (461,102.4) node [anchor=north west][inner sep=0.75pt]    {$\hat{\pi }$};
\draw (416,182.4) node [anchor=north west][inner sep=0.75pt]    {$\hat{\pi }$};
\draw (275,-5.6) node [anchor=north west][inner sep=0.75pt]    {$o$};
\draw (356,102.4) node [anchor=north west][inner sep=0.75pt]    {$\hat{\pi }$};
\draw (146,100.4) node [anchor=north west][inner sep=0.75pt]    {$\hat{\pi }$};
\draw (301,182.4) node [anchor=north west][inner sep=0.75pt]    {$\hat{\pi }$};
\draw (111,182.4) node [anchor=north west][inner sep=0.75pt]    {$\hat{\pi }$};
\draw (171,182.4) node [anchor=north west][inner sep=0.75pt]    {$\hat{\pi }$};
\draw (361,182.4) node [anchor=north west][inner sep=0.75pt]    {$\hat{\pi }$};
\draw (256,102.4) node [anchor=north west][inner sep=0.75pt]    {$\hat{\pi }$};
\draw (456,182.4) node [anchor=north west][inner sep=0.75pt]    {$\hat{\pi }$};
\draw (86,250.4) node [anchor=north west][inner sep=0.75pt]    {$\hat{\pi }$};
\draw (126,250.4) node [anchor=north west][inner sep=0.75pt]    {$\hat{\pi }$};
\draw (296,250.4) node [anchor=north west][inner sep=0.75pt]    {$\hat{\pi }$};
\draw (386,250.4) node [anchor=north west][inner sep=0.75pt]    {$\hat{\pi }$};
\draw (416,250.4) node [anchor=north west][inner sep=0.75pt]    {$\hat{\pi }$};
\draw (446,250.4) node [anchor=north west][inner sep=0.75pt]    {$\hat{\pi }$};
\end{tikzpicture}
\caption{A $3-$Neighbourhood of a vertex-rooted UBGW tree with the law of the number of children drawn on every vertex.} \label{fig:VUBGW}
\end{figure}

\bigskip

The most classical examples of random graphs converging to UBGW trees we consider are sparse Erdős–Rényi  and configuration models:
\begin{itemize}
    \item \textbf{Sparse Erdős–Rényi:} Introduced in the celebrated paper of Erdős and Rényi \cite{ErdosRenyi}, for $c>0$ and $N\geq 1$, the random graph $\mathcal G (N,\frac{c}{N})$ is defined on the vertex set $\{1,\ldots, N \}$ with independent edges between vertices with probability $\frac{c}{N}$. Once uniformly rooted, these graphs converge locally when $N$ goes to $\infty$, to a UBGW tree with reproduction law Poisson with parameter $c$.
    \item \textbf{Configuration model:} This model was introduced by Bollob\'as in 1980 \cite{B1980} and can be defined as follows.
Let $N\geq 1$ be an integer and let $d_1, \ldots, d_N \in \mathbb{Z}_+$ be such that $d_1 + \cdots +d_N$ is even. We interpret $d_i$ as a number of half-edges attached to vertex $i$. Then, the configuration model associated to the sequence $(d_i)_{1 \leq i \leq N}$ is the random multigraph with vertex set $\{1, \ldots ,N\}$ obtained by a uniform matching of these half-edges. If $d_1 + \cdots + d_N$ is odd, we change $d_N$ into $d_N+1$ and do the same construction. Now, let $\mathbf{d}^{(N)}$ be a sequence of random variables defined on the same probability space $(\Omega,\mathcal{F},\mathbb{P})$ such that for every $N\geq 1$, $\mathbf{d}^ {(N)} = (d_1^{(N)}, \ldots, d_N^{(N)}) \in \mathbb{Z}_+^N$. Furthermore, suppose that there exists $\pi$ a probability measure on $\mathbb{Z}_+$ with finite first moment such that
\[  \forall k \geq 0, \quad \frac{1}{N} \sum\limits_{ j=1}^N \mathbbm{1}_{ d_j^{(N)} = k}  \underset{ N \rightarrow +\infty}{ \longrightarrow} \pi( \{k\}).   \]
The sequence of random configuration graphs associated to the $\mathbf{d}^{(N)}$ has asymptotically a positive probability to be simple. In addition, this sequence of random graphs, when uniformly rooted, converges locally in law to the UBGW random tree with offspring distribution $\pi$, see \cite{vanderHofstadBook2} for more details.
\end{itemize}

\subsection{Matchings, optimality and unimodularity}
We start with the definition of matchings on a graph.
\begin{definition}
 For any weighted graph $G=((V,E),w)$, a matching $M=(V,E')$ on $G$ is a subgraph of $G$ such that $E' \subset E$ and  every vertex of $V$ belongs to at most one edge of $E'$. A matched graph is a pair $(G,M)$, where $M$ is a matching on $G$.
\end{definition}
We can extend the notion of unimodularity for random graphs to random matched graphs by transforming a matching into a canonical decoration. If $(G,M)=((V,E),(V,E'))$  is a matched graph then we can define the associated decoration $\mathbbm{1}_{M}:E\to\{0,1\}$ with  $\mathbbm{1}_{M}(u,v)=1$ if and only if $(u,v) \in E'$.
\begin{definition}
Let $(\mathbb{G},\mathbb M,\mathbf{o})$ be a random matched rooted graph. 
    We say that $(\mathbb{G},\mathbb{M},\mathbf{o})$ is unimodular iff $(\mathbb{G},\mathbf{o},\mathbbm{1}_{\mathbb{M}})$ is unimodular.
\end{definition}
The correspondence between $(\mathbb{G},\mathbb{M},\mathbf{o})$ and $(\mathbb{G},\mathbf{o},\mathbbm{1}_{\mathbb{M}})$ gives a representation of the space of unimodular laws on decorated matched graphs as a subspace of unimodular laws on decorated graphs $\mathcal{L}_{U,\mathcal{M}}(\mathcal{G}^{\star}) \subseteq \mathcal{L}_{U}(\mathcal{G}^{\star})$.
\begin{definition}
The subspace $\mathcal{L}_{U,\mathcal{M}}(\mathcal{G}^{\star}) \subseteq \mathcal{L}_{U}(\mathcal{G}^{\star})$ is defined as the set of 
elements $\mu$ of $\mathcal{L}_{U}(\mathcal{G}^{\star})$ such that for $(G,o,w,(f_i)_{i \in \{1,..,I\}}) \sim \mu$, we have that $f_I(u,v)$ maps $E$ to $\{0,1\}$ and the subset $\{(u,v), f_I(u,v)=1 \}$ almost surely induces a matching on $G$ in the sense that
\begin{align*}
    \mathbb{P}\left( \exists u,u',v \in V, u \neq u' \text{ and } f(v,u)=f(v,u')=1\right)=0.
\end{align*}
Similarly we can define the edge rooted version $\mathcal{L}_{U,\mathcal{M}}(\hat{\mathcal{G}}^{\star})$.
\end{definition}

We recall that the central object of interest is optimal matching, which corresponds to the classical notion of maximal weight matching in finite graphs.
When $G$ is infinite, this optimality is ill-defined. However, in the case of a unimodular random weighted graph, since the root edge is informally a typical edge, we can define optimality via its expected weight when it belongs to the matching. This leads to the following definition:
\begin{definition}
    Let $(\mathbb{G}, \mathbf o, \mathbb M )$ be a unimodular random matched rooted graph, we define the performance of $(\mathbb{G},\mathbf o,\mathbb{M})$ as:
    \begin{align*}
     \perfE(\mathbb{G},\mathbf o, \mathbb{M}) &:= \mathbb{E}\left[ w((\mathbf{o}_-,\mathbf{o}_+))\mathbbm{1}_{(\mathbf{o}_-,\mathbf{o}_+) \in \mathbb{M}} \right] \text{ in the edge-rooted setting,}  \\
    \perfV(\mathbb{G},\mathbf o,\mathbb{M}) &:=\mathbb{E}\left[ \sum_{v \sim o}w(\mathbf{o},v)\mathbbm{1}_{(\mathbf{o},v) \in \mathbb{M}} \right] \text{in the vertex-rooted setting.} 
    \end{align*}
\end{definition}

When the context is clear, we will shorten the notation to $\perf (\mathbb{M})$. 
By extension, since those quantities only depend on the law of $(\mathbb{G},\mathbf o,\mathbb{M})$, we will freely use the same notation $\perfE(\mu_{e}):=\perfE(\mathbb{G},\mathbf{o},\mathbb{M})$ for $(\mathbb{G},\mathbf{o},\mathbbm{1}_\mathbb{M}) \sim \mu_{e}$ and $\perfV(\mu_{v}):=\perfV(\mathbb{G},\mathbf{o},\mathbb{M})$ for $(\mathbb{G},\mathbf{o},\mathbbm{1}_\mathbb{M}) \sim \mu_{v}$ where $\mu_{e}\in \mathcal{L}_{U,\mathcal{M}}(\hat{\mathcal{G}}^{\star})$ and $\mu_{v} \in \mathcal{L}_{U,\mathcal{M}}(\mathcal{G}^{\star})$.

\bigskip

Let $(\mathbb G, \mathbf o)$ be a (undecorated) unimodular vertex-rooted graph. We say that $(\mathbb{G}',\mathbf{o}',\mathbb{M})$ is optimal if it is a unimodular matched vertex-rooted graph such that $\perfV((\mathbb{G}',\mathbf{o}'),\mathbb{M})$ is maximal among all unimodular rooted matched random graphs $(\mathbb{G}',\mathbf{o}',\mathbb{M}')$ such that $(\mathbb{G}',\mathbf{o}')$ has the same law as $(\mathbb{G},\mathbf{o})$.

\begin{remark}
    One could be tempted to look at performance only for objects of the form $(\mathbb{G},\mathbf{o},F(\mathbb{G},o))$ where $F$ deterministically  maps vertex-rooted weighted graphs $(V,E,o,w)$ to decorations from $E$ to $\{0,1\}$ that induces a matching. This paper does not investigate whether this restriction is relevant. 
\end{remark}
In the edge-rooted setting, optimality is defined similarly.
The next proposition shows that the operator $R$, introduced in Definition~\ref{def:R}, preserves optimality.

\begin{prop}\label{prop:chgmtpdv}
    Let $\mu_{v} \in \mathcal{L}_{U,\mathcal{M}}(\mathcal{G}^{\star})$, assume $m = \int_{\mathcal G^{\star}} \deg (o) \, \mathrm{d} \mu_v(G,o) < \infty$, then:
    \[  \perfV(\mu_{v}) =m\perfE(R(\mu_{v})).  \]
\end{prop}

    
    
    
\begin{remark}
If the graph is finite and the root is chosen uniformly either among the vertices or among the directed edges, then $\perfV$ is simply the average contribution per vertex, and $\perfE$ is the average contribution per directed edge. 
It is then clear that the two quantities are proportional. The proposition shows that it generalises to unimodular matched graphs. 
\end{remark}
\begin{proof}
We will decompose both quantities with respect to the degree of the root.  Let us write $(\mathbb{G},o,\mathbbm{1}_\mathbb{M}) \sim \mu_{v} $ and  $(\hat{\mathbb{G}}, (o_-,o_+), \mathbbm{1}_{\hat{\mathbb{M}}}) \sim R(\mu_{v})$.
In the edge-rooted case, we get:
\[\sum_{k=1}^{\infty}  \mathbb{E}[w((o_-,o_+))\mathbbm{1}_{(o_-,o_+) \in \hat{\mathbb{M}}}\mathbbm{1}_{\deg(o_-)=k}] .\]
while in the vertex-rooted case, it reads as:
\[ \sum_{k=0}^{\infty} \mathbb{E}\left[ \sum_{v \sim o} w(o,v)\mathbbm{1}_{(o,v) \in \mathbb{M}}  \mathbbm{1}_{\deg(o)=k}\right]=\sum_{k=1}^{\infty} \mathbb{E}\left[ \sum_{v \sim o} w(o,v)\mathbbm{1}_{(o,v) \in \mathbb{M}}  \mathbbm{1}_{\deg(o)=k}\right] .  \]
For each $k >0$, we shall prove the following that will give the result:
\[\mathbb{E}\left[ \sum_{v \sim o} w(o,v)\mathbbm{1}_{(o,v) \in \mathbb{M}} \mathbbm{1}_{\deg(o)=k} \right]= \mathbb{E}[\deg(o)]\mathbb{E}\left[w(o_-,o_+)\mathbbm{1}_{(o_-,o_+)\in \hat{\mathbb{M}}} \mathbbm{1}_{\deg(o_-)=k}\right]  .\]
First, let us show that 
\[\mathbb{E}\left[w(o_-,o_+)\mathbbm{1}_{(o_-,o_+)\in \hat{\mathbb{M}}} \mathbbm{1}_{\deg(o_-)=k}\right] = \frac{1}{k}\mathbb{E}\left[ \sum_{v \sim o_-} w(o_-,v)\mathbbm{1}_{(o_-,v) \in \hat{\mathbb{M}}} \mathbbm{1}_{\deg(o_-)=k} \right].  \]
For this purpose, we use reversibility to show that:
\[ \mathbb{E}\left[w(o_-,o_+)\mathbbm{1}_{(o_-,o_+)\in \hat{\mathbb{M}}} \mathbbm{1}_{\deg(o_-)=k}\right]=\mathbb{E}\left[w(o_+,o_-)\mathbbm{1}_{(o_+,o_-)\in \hat{\mathbb{M}}} \mathbbm{1}_{\deg(o_+)=k}\right]. \]
Then use the 1-step stationarity to show that:
\[ \mathbb{E}\left[w(o_+,o_-)\mathbbm{1}_{(o_+,o_-)\in \hat{\mathbb{M}}} \mathbbm{1}_{\deg(o_+)=k}\right] = \mathbb{E}\left[\frac{1}{k-1}\sum_{\substack{v \sim o_- \\ v \neq o_+}}w(o_-,v)\mathbbm{1}_{(o_-,v)\in \hat{\mathbb{M}}} \mathbbm{1}_{\deg(o_-)=k}\right],  \]
hence:
\begin{align*}
    &\mathbb{E}\left[ \sum_{v \sim o_-} w(o_-,v)\mathbbm{1}_{(o_-,v) \in \hat{\mathbb{M}}} \mathbbm{1}_{\deg(o_-)=k} \right] \\
    & \quad= \mathbb{E}\left[ \sum_{\substack{v \sim o_-\\ v \neq o_+}} w(o_-,v)\mathbbm{1}_{(o_-,v) \in \hat{\mathbb{M}}} \mathbbm{1}_{\deg(o_-)=k} \right] 
    + \mathbb{E}\left[w(o_-,o_+)\mathbbm{1}_{(o_-,o_+)\in \hat{\mathbb{M}}} \mathbbm{1}_{\deg(o_-)=k}\right] \\
    &\quad = (k-1) \mathbb{E}\left[w(o_-,o_+)\mathbbm{1}_{(o_-,o_+)\in \hat{\mathbb{M}}} \mathbbm{1}_{\deg(o_-)=k}\right] + \mathbb{E}\left[w(o_-,o_+)\mathbbm{1}_{(o_-,o_+)\in \hat{\mathbb{M}}} \mathbbm{1}_{\deg(o_-)=k}\right] \\
    &\quad = k \mathbb{E}\left[w(o_-,o_+)\mathbbm{1}_{(o_-,o_+)\in \hat{\mathbb{M}}} \mathbbm{1}_{\deg(o_-)=k}\right].
\end{align*}
Recall that the operator $R$ is the operation of taking $o_-=o$, and $o_+$ uniformly chosen among the neighbors of $o$, then biasing by $\deg(o)$.
Taking $f=\sum_{v \sim o_-} w(o_-,v)\mathbbm{1}_{(o_-,v) \in \hat{\mathbb{M}}} \mathbbm{1}_{\deg(o_-)=k}$ and applying Definition~\ref{def:R} on $f$:
\begin{align*}
    &\mathbb{E}\left[\sum_{v \sim o_-} w(o_-,v)\mathbbm{1}_{(o_-,v) \in \hat{\mathbb{M}}} \mathbbm{1}_{\deg(o_-)=k}\right] = \frac{1}{\mathbb{E}[\deg(o)]} \mathbb{E}\left[\deg(o) \sum_{v \sim o} w(o,v)\mathbbm{1}_{(o,v) \in \mathbb{M}} \mathbbm{1}_{\deg(o)=k} \right]
\end{align*}
In this second sum, $\deg(o)\mathbbm{1}_{\deg(o)=k}=k\mathbbm{1}_{\deg(o)=k}$ so it is the same as :
\[\frac{k}{\mathbb{E}[\deg(o)]} \mathbb{E}\left[ \sum_{v \sim o} w(o,v)\mathbbm{1}_{(o,v) \in \mathbb{M}} \mathbbm{1}_{\deg(o)=k} \right]. \]
Putting everything together, we have shown that
\[ \frac{k}{\mathbb{E}[\deg(o)]} \mathbb{E}\left[ \sum_{v \sim o} w(o,v)\mathbbm{1}_{(o,v) \in \mathbb{M}} \mathbbm{1}_{\deg(o)=k} \right]=  k \mathbb{E}\left[w((o_-,o_+))\mathbbm{1}_{(o_-,o_+)\in \hat{\mathbb{M}}} \mathbbm{1}_{\deg(o_-)=k}\right],\]\
which immediately entails the result.
\end{proof}
We state a useful property of unimodular graphs: events that have probability zero (resp. almost sure) at the root have probability zero (resp. almost sure) everywhere.
\begin{prop}[Lemma 2.3 in \cite{aldous2018processes}]
    Let $(\mathbf{G},\mathbf{o})$ be a unimodular graph, let $H>0$ and $f$ be a non-negative H-local function on $\mathcal{G}^{\star}$ in the sense that if $(G,o),(G',o') \in \mathcal{G}^{\star}$ such that $B_H(G,o)=B_H(G',o')$, $f(G,o)=f(G',o')$.
    Assume that 
    \[ \mathbb{E}\left[ f(\mathbf{G},\mathbf{o})\right]=0 .\]
    Then almost surely, for all $v \in \mathbf{V}$ 
    \[ f(\mathbf{G},v)=0. \]
\end{prop}

\section{Optimal unimodular matchings on UBGW trees : Theorem~\ref{unicityUBGW}} \label{sec:matchUBGW}

In this section, we construct optimal matchings on UBGW trees. As mentioned in the introduction, our construction relies on a message passing algorithm that we will present and study in depth in Section~\ref{sec:message}. We then prove that the matchings constructed are optimal in Section~\ref{sec:optimality}. Finally, we prove the uniqueness of optimal matchings in UBGW trees in Section~\ref{sec:uniqueness}.

By virtue of Proposition~\ref{prop:chgmtpdv}, it is equivalent to study optimal matchings from either vertex-rooted or edge-rooted point of view. We will change the point of view throughout the paper depending on which one is the most suitable to the situation. An indication will be given at the beginning of each subsection.

\subsection{A message passing algorithm} \label{sec:message}

In this subsection, we shall adopt the edge-rooted point of view. Before formally introducing the message passing algorithm mentioned in the introduction, we start by discussing the simpler setting of finite trees. We are looking for a dynamic program that builds the maximum matching.

\begin{figure}[t!]
\centering
\tikzset{every picture/.style={line width=0.75pt}} 

\begin{tikzpicture}[x=0.75pt,y=0.75pt,yscale=-1,xscale=1]

\draw    (170,161) -- (290,161) ;
\draw    (90,111) -- (170,161) ;
\draw    (90,161) -- (170,161) ;
\draw    (170,161) -- (90,211) ;
\draw    (290,161) -- (370,161) ;
\draw    (450,131) -- (370,161) ;
\draw    (450,191) -- (370,161) ;
\draw  [fill={rgb, 255:red, 0; green, 0; blue, 0 }  ,fill opacity=1 ] (165.65,161) .. controls (165.65,158.6) and (167.6,156.65) .. (170,156.65) .. controls (172.4,156.65) and (174.35,158.6) .. (174.35,161) .. controls (174.35,163.4) and (172.4,165.35) .. (170,165.35) .. controls (167.6,165.35) and (165.65,163.4) .. (165.65,161) -- cycle ;
\draw    (20,91) -- (90,111) ;
\draw    (20,131) -- (90,111) ;
\draw    (20,111) -- (90,111) ;
\draw    (20,151) -- (90,161) ;
\draw    (20,171) -- (90,161) ;
\draw  [fill={rgb, 255:red, 0; green, 0; blue, 0 }  ,fill opacity=1 ] (285.65,161) .. controls (285.65,158.6) and (287.6,156.65) .. (290,156.65) .. controls (292.4,156.65) and (294.35,158.6) .. (294.35,161) .. controls (294.35,163.4) and (292.4,165.35) .. (290,165.35) .. controls (287.6,165.35) and (285.65,163.4) .. (285.65,161) -- cycle ;
\draw  [fill={rgb, 255:red, 0; green, 0; blue, 0 }  ,fill opacity=1 ] (85.65,211) .. controls (85.65,208.6) and (87.6,206.65) .. (90,206.65) .. controls (92.4,206.65) and (94.35,208.6) .. (94.35,211) .. controls (94.35,213.4) and (92.4,215.35) .. (90,215.35) .. controls (87.6,215.35) and (85.65,213.4) .. (85.65,211) -- cycle ;
\draw  [fill={rgb, 255:red, 0; green, 0; blue, 0 }  ,fill opacity=1 ] (85.65,161) .. controls (85.65,158.6) and (87.6,156.65) .. (90,156.65) .. controls (92.4,156.65) and (94.35,158.6) .. (94.35,161) .. controls (94.35,163.4) and (92.4,165.35) .. (90,165.35) .. controls (87.6,165.35) and (85.65,163.4) .. (85.65,161) -- cycle ;
\draw  [fill={rgb, 255:red, 0; green, 0; blue, 0 }  ,fill opacity=1 ] (85.65,111) .. controls (85.65,108.6) and (87.6,106.65) .. (90,106.65) .. controls (92.4,106.65) and (94.35,108.6) .. (94.35,111) .. controls (94.35,113.4) and (92.4,115.35) .. (90,115.35) .. controls (87.6,115.35) and (85.65,113.4) .. (85.65,111) -- cycle ;
\draw  [fill={rgb, 255:red, 0; green, 0; blue, 0 }  ,fill opacity=1 ] (365.65,161) .. controls (365.65,158.6) and (367.6,156.65) .. (370,156.65) .. controls (372.4,156.65) and (374.35,158.6) .. (374.35,161) .. controls (374.35,163.4) and (372.4,165.35) .. (370,165.35) .. controls (367.6,165.35) and (365.65,163.4) .. (365.65,161) -- cycle ;
\draw  [draw opacity=0] (31.64,259.84) .. controls (117.13,256.93) and (185,213.34) .. (185,160) .. controls (185,106.32) and (116.26,62.51) .. (29.99,60.1) -- (22.5,160) -- cycle ; \draw  [color={rgb, 255:red, 208; green, 2; blue, 27 }  ,draw opacity=1 ] (31.64,259.84) .. controls (117.13,256.93) and (185,213.34) .. (185,160) .. controls (185,106.32) and (116.26,62.51) .. (29.99,60.1) ;  
\draw  [draw opacity=0] (30.04,248.19) .. controls (87.11,239.84) and (130.05,203.49) .. (130.03,159.93) .. controls (130.02,116.36) and (87.02,80.03) .. (29.92,71.74) -- (5.03,159.98) -- cycle ; \draw  [color={rgb, 255:red, 126; green, 211; blue, 33 }  ,draw opacity=1 ] (30.04,248.19) .. controls (87.11,239.84) and (130.05,203.49) .. (130.03,159.93) .. controls (130.02,116.36) and (87.02,80.03) .. (29.92,71.74) ;  
\draw  [draw opacity=0] (438.13,249.68) .. controls (433.8,249.89) and (429.42,250) .. (425,250) .. controls (339.4,250) and (270,209.71) .. (270,160) .. controls (270,110.29) and (339.4,70) .. (425,70) .. controls (427.64,70) and (430.26,70.04) .. (432.86,70.11) -- (425,160) -- cycle ; \draw  [color={rgb, 255:red, 144; green, 19; blue, 254 }  ,draw opacity=1 ] (438.13,249.68) .. controls (433.8,249.89) and (429.42,250) .. (425,250) .. controls (339.4,250) and (270,209.71) .. (270,160) .. controls (270,110.29) and (339.4,70) .. (425,70) .. controls (427.64,70) and (430.26,70.04) .. (432.86,70.11) ;  
\draw  [draw opacity=0] (434.75,229.82) .. controls (375.09,227.25) and (328,196.97) .. (328,160) .. controls (328,122.96) and (375.27,92.63) .. (435.1,90.16) -- (443,160) -- cycle ; \draw  [color={rgb, 255:red, 74; green, 144; blue, 226 }  ,draw opacity=1 ] (434.75,229.82) .. controls (375.09,227.25) and (328,196.97) .. (328,160) .. controls (328,122.96) and (375.27,92.63) .. (435.1,90.16) ;  

\draw (171,133.4) node [anchor=north west][inner sep=0.75pt]    {$u$};
\draw (291,133.4) node [anchor=north west][inner sep=0.75pt]    {$v$};
\draw (209,132.4) node [anchor=north west][inner sep=0.75pt]    {$w( u,v)$};
\draw (161,282.4) node [anchor=north west][inner sep=0.75pt]    {$ \begin{array}{l}
Z( u,v) =
OPT(\textcolor[rgb]{0.56,0.07,1}{T_{(u,v)}}) -OPT(\textcolor[rgb]{0.29,0.56,0.89}{T_{(u,v)} \setminus \{v\}}) .
\end{array}$};
\draw (311,52.4) node [anchor=north west][inner sep=0.75pt]  [color={rgb, 255:red, 144; green, 19; blue, 254 }  ,opacity=1 ]  {$T_{(u,v)}$};
\draw (13,192.4) node [anchor=north west][inner sep=0.75pt]  [color={rgb, 255:red, 126; green, 211; blue, 33 }  ,opacity=1 ]  {$T_{(v,u)} \setminus \{u\}$};
\draw (129,52.4) node [anchor=north west][inner sep=0.75pt]  [color={rgb, 255:red, 208; green, 2; blue, 27 }  ,opacity=1 ]  {$T_{(v,u)}$};
\draw (371,192.4) node [anchor=north west][inner sep=0.75pt]  [color={rgb, 255:red, 74; green, 144; blue, 226 }  ,opacity=1 ]  {$\textcolor[rgb]{0.29,0.56,0.89}{T_{(u,v)} \setminus \{v\}}$};
\draw (161,332.4) node [anchor=north west][inner sep=0.75pt]    {$ \begin{array}{l}
Z( v,u) = OPT(\textcolor[rgb]{0.82,0.01,0.11}{T_{(v,u)}}) -OPT(\textcolor[rgb]{0.49,0.83,0.13}{T_{(v,u)} \setminus \{u\}}) .
\end{array}$};
\draw (1,281.4) node [anchor=north west][inner sep=0.75pt]    {$ \begin{array}{l}
OPT( G) =\text{weight of} \\
\text{a maximum matching}\\
\text{on}\ G.
\end{array}$};

\end{tikzpicture}
\caption{Definitions of $Z$.}\label{fig:Zdefinitions}
\end{figure}

Fix a finite weighted rooted deterministic tree $T$ with a unique optimal matching $\mathbb{M}_{\mathrm{opt}}$. Let $\{u,v\}$ be an edge of $T$, we denote by $T_{(v,u)}$ and $T_{(u,v)}$ the two connected components of $T \setminus \{u,v\}$ containing respectively $u$ and $v$. Let us start by simple but key observations illustrated in Figure~\ref{fig:Zdefinitions}:
\begin{itemize}
    \item The maximal weight of matchings of $T$ that exclude the edge $\{u,v \}$ is merely the sum of the maximal weights of matchings of $T_{(v,u)}$ and $T_{(u,v)}$, denoted $OPT(T_{(v,u)})$ and $OPT(T_{(u,v)})$.
    \item The maximal weight of matchings that include the edge $\{u,v\}$ is the sum of the weight of $\{u,v\}$ and of the maximal weights of matchings of $T_{(v,u)} \setminus \{ u \} $ and $T_{(u,v)} \setminus \{ v \}$. Note that both $T_{(v,u)} \setminus \{ u \} $ and $T_{(u,v)} \setminus \{ v \}$ consist of collections of disjoint subtrees of $T$ issued from the children of $u$ and $v$. We denote by $OPT(T_{(v,u)} \setminus \{ u\})$ and $OPT(T_{(u,v)} \setminus \{ v \})$ the relevant maximal weights.
\end{itemize}
From this discussion, one can see that the edge $\{ u,v \}$ is in the optimal matching of $T$ iff
\[
w(u,v) > OPT(T_{(v,u)}) + OPT(T_{(u,v)}) - \Big(OPT(T_{(v,u)} \setminus \{u\}) + OPT(T_{(u,v)} \setminus \{v\}) \Big).
\]
It will be instrumental to isolate quantities depending only on $T_{(v,u)}$ and $T_{(u,v)}$ in the previous display. This leads us to introduce the following quantities
\begin{align*}
Z(u,v)&= OPT( T_{(u,v)})-OPT(T_{(u,v)} \setminus \{v\}), \\
Z(v,u)&=OPT(T_{(v,u)} )- OPT( T_{(v,u)} \setminus \{u\}), 
\end{align*}
and the criterion of $(u,v) \in \mathbb{M}_{\mathrm{opt}}$ is simply 
\begin{equation}\label{eq:defmatching}
 w(u,v)> Z(u,v)+Z(v,u) .
\end{equation}
Note that the variable $Z(u,v)$ has a neat interpretation  in terms of the matchings of $T_{(u,v)}$. Indeed, it is the marginal gain between allowing $v$ to be matched or not.

\bigskip

\begin{figure}[t!]
\tikzset{
pattern size/.store in=\mcSize, 
pattern size = 5pt,
pattern thickness/.store in=\mcThickness, 
pattern thickness = 0.3pt,
pattern radius/.store in=\mcRadius, 
pattern radius = 1pt}
\makeatletter
\pgfutil@ifundefined{pgf@pattern@name@_4nnhtjt6r}{
\pgfdeclarepatternformonly[\mcThickness,\mcSize]{_4nnhtjt6r}
{\pgfqpoint{0pt}{0pt}}
{\pgfpoint{\mcSize+\mcThickness}{\mcSize+\mcThickness}}
{\pgfpoint{\mcSize}{\mcSize}}
{
\pgfsetcolor{\tikz@pattern@color}
\pgfsetlinewidth{\mcThickness}
\pgfpathmoveto{\pgfqpoint{0pt}{0pt}}
\pgfpathlineto{\pgfpoint{\mcSize+\mcThickness}{\mcSize+\mcThickness}}
\pgfusepath{stroke}
}}
\makeatother

 
\tikzset{
pattern size/.store in=\mcSize, 
pattern size = 5pt,
pattern thickness/.store in=\mcThickness, 
pattern thickness = 0.3pt,
pattern radius/.store in=\mcRadius, 
pattern radius = 1pt}
\makeatletter
\pgfutil@ifundefined{pgf@pattern@name@_ik3oy0ks0}{
\pgfdeclarepatternformonly[\mcThickness,\mcSize]{_ik3oy0ks0}
{\pgfqpoint{0pt}{0pt}}
{\pgfpoint{\mcSize+\mcThickness}{\mcSize+\mcThickness}}
{\pgfpoint{\mcSize}{\mcSize}}
{
\pgfsetcolor{\tikz@pattern@color}
\pgfsetlinewidth{\mcThickness}
\pgfpathmoveto{\pgfqpoint{0pt}{0pt}}
\pgfpathlineto{\pgfpoint{\mcSize+\mcThickness}{\mcSize+\mcThickness}}
\pgfusepath{stroke}
}}
\makeatother
\tikzset{every picture/.style={line width=0.75pt}} 

\centering
\begin{tikzpicture}[x=0.75pt,y=0.75pt,yscale=-1,xscale=1]

\draw [color={rgb, 255:red, 144; green, 19; blue, 254 }  ,draw opacity=1 ]   (30,160) -- (130,50) ;

\draw (0,160) -- (30,160);
\draw    (30,160) -- (130,110) ;
\draw    (30,160) -- (130,160) ;
\draw    (30,160) -- (130,210) ;
\draw    (30,160) -- (130,270) ;
\draw    (130,50) -- (230,10) ;
\draw    (130,50) -- (230,50) ;
\draw    (130,50) -- (230,90) ;
\draw [pattern=_4nnhtjt6r,pattern size=6pt,pattern thickness=0.75pt,pattern radius=0pt, pattern color={rgb, 255:red, 0; green, 0; blue, 0}]   (110,100) -- (110,290) ;
\draw    (110,100) -- (200,100) ;
\draw    (110,290) -- (200,290) ;
\draw  [color={rgb, 255:red, 255; green, 255; blue, 255 }  ,draw opacity=1 ][pattern=_ik3oy0ks0,pattern size=19.049999999999997pt,pattern thickness=0.75pt,pattern radius=0pt, pattern color={rgb, 255:red, 208; green, 2; blue, 27}] (110,100) -- (200,100) -- (200,290) -- (110,290) -- cycle ;
\draw  [draw opacity=0] (235.37,99.24) .. controls (209.58,95.1) and (190,74.63) .. (190,50) .. controls (190,26.04) and (208.54,6.02) .. (233.29,1.14) -- (245,50) -- cycle ; \draw  [color={rgb, 255:red, 144; green, 19; blue, 254 }  ,draw opacity=1 ] (235.37,99.24) .. controls (209.58,95.1) and (190,74.63) .. (190,50) .. controls (190,26.04) and (208.54,6.02) .. (233.29,1.14) ;  
\draw  [draw opacity=0] (219.11,99.76) .. controls (163.47,97.27) and (120,75.94) .. (120,50) .. controls (120,23.71) and (164.63,2.16) .. (221.34,0.15) -- (230,50) -- cycle ; \draw  [color={rgb, 255:red, 74; green, 144; blue, 226 }  ,draw opacity=1 ] (219.11,99.76) .. controls (163.47,97.27) and (120,75.94) .. (120,50) .. controls (120,23.71) and (164.63,2.16) .. (221.34,0.15) ;  

\draw (111,22.4) node [anchor=north west][inner sep=0.75pt]    {$v_{i}$};
\draw (17,142.4) node [anchor=north west][inner sep=0.75pt]    {$v$};
\draw (-13,142.4) node [anchor=north west][inner sep=0.75pt]    {$u$};
\draw (119,242.4) node [anchor=north west][inner sep=0.75pt]    {$v_{j}$};
\draw (211,171) node [anchor=north west][inner sep=0.75pt]   [align=left] {$\textcolor[rgb]{0.82,0.01,0.11}{{\text{common to both } OPT(T_{(u,v)}) \text{ and }  OPT(T_{(u,v)}\setminus \{ v \})}}$};
\draw (241,30.4) node [anchor=north west][inner sep=0.75pt]    {$\textcolor[rgb]{0.56,0.07,1}{\text{in } OPT( T_{(u,v)})}$};
\draw (241,60.4) node [anchor=north west][inner sep=0.75pt]    {$\textcolor[rgb]{0.29,0.56,0.89}{\text{in} \ OPT( T_{(u,v)} \ \setminus \{v\})}$};

\end{tikzpicture}
\caption{Illustration of the deduction of the recursive equation assuming $v_i$ is matched to $v$.}\label{fig:Zrecursion}
\end{figure}
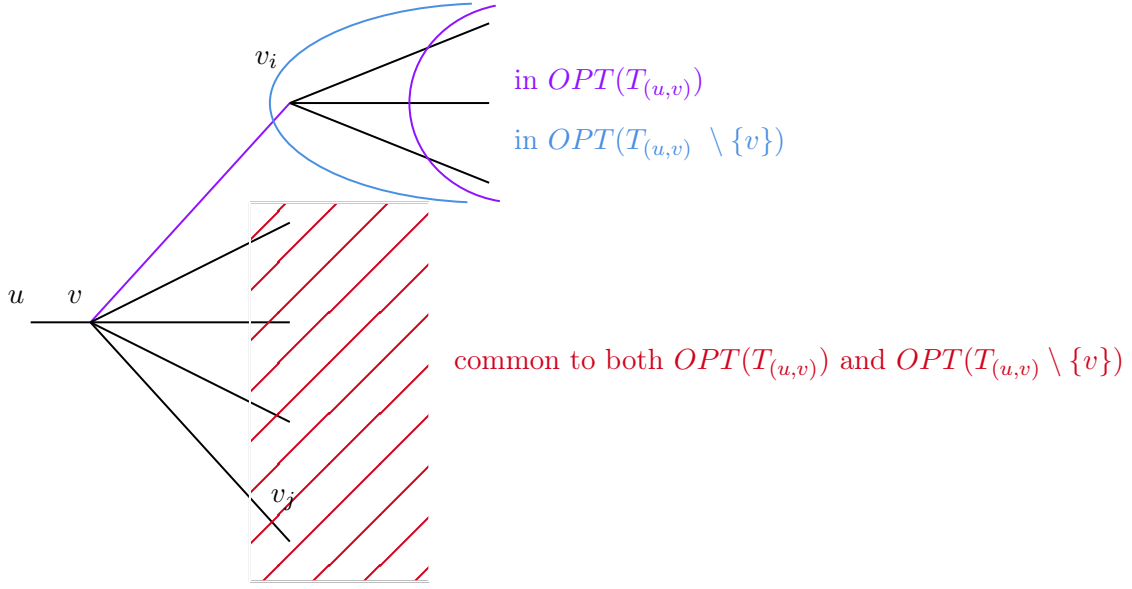

The variables $Z$ have the nice property of satisfying a recursive equation. We describe this recursion for $Z(u,v)$ and $T_{(u,v)}$. See Figure~\ref{fig:Zrecursion} for an illustration.

Listing $v_1,...v_k$ the children of $v$ in $T_{(u,v)}$, assume that the maximum matching of $T_{(u,v)}$ matches $v$ with $v_i$. 
In that situation, the maximum matching of $T_{(u,v)}$ and the maximum matching of $T_{(u,v)} \setminus \{v\}$ coincide on the subtrees $T_{(v,v_j)}$ for $j \neq i$.
On $T_{(v,v_i)}$, this maximum matching has matched $v_i$ with $v$, so $v_i$ is not matched to other vertices. Therefore, our maximum matching on $T_{(u,v)}$ restricted to $T_{(v,v_i)}$ is the union of $\{v,v_i\}$ with the maximum matching of
$T_{(v,v_i)} \setminus \{v_i\}$. In that case, the weight of the maximal matching is given by:
\[
w(v,v_i) + OPT(T_{(v,v_i)} \setminus \{v_i\}) + \sum_{j \neq i} OPT(T_{(v,v_j)}).
\]

On the other hand, if a vertex $v$ is not matched in the maximum matching of $T_{(u,v)}$, then the later is also the maximum matching in $T_{(u,v)} \setminus \{v\}$. Thus, inside each sub-tree $T_{(v,v_i)}$, it coincides with the maximal matching of $T_{(v,v_i)}$.
In that case, the weight of the maximal matching is given by:
\[ \sum_{j} OPT(T_{(v,v_j)}).  \]

Putting all the different cases together, we have the identity
\[
OPT(T_{(u,v)}) = \max \left\{ \sum_{j} OPT(T_{(v,v_j)}) \, , \max_{i\in \{ 1, \ldots , k \} } \left\{ w(v,v_i) + OPT(T_{(v,v_i)} \setminus \{v_i\}) + \sum_{j \neq i} OPT(T_{(v,v_j)}) \right\} \right\}.
\]
Recalling the definition of $Z(u,v)$, we get:
\begin{align*}
Z(u,v)&= 
OPT(T_{(u,v)})
-\sum_{j} OPT(T_{(v,v_j)}) \\
&= \max\left\{0, \max_{i\in \{ 1, \ldots , k \} } \left\{w(v,v_i)-(OPT(T_{(v,v_i)})-OPT(T_{(v,v_i)} \setminus \{v_i\}) \right\} \right\} \\
&= \max\left\{0,\max_{i\in \{ 1, \ldots , k \} } w(v,v_i)-Z(v,v_i)\right\}.
\end{align*}
In conclusion, we have obtained the recursive Equation \eqref{eq:Z_hrecursion} given in the introduction.

Note that, since $T$ is a finite tree, it is possible to calculate $Z(u,v)$ for all vertices $u \sim v \in T$ by starting when $v$ is a leaf, in which case $Z(u,v) = 0$. By construction, our decision rule \eqref{eq:decisionrule} constructs the optimal matching on $T$ from the values of $Z$.

\bigskip

Now, let us discuss how to extend the construction of variables $Z(u,v)$ when the underlying tree is a (possibly infinite) UBGW tree $\mathcal{T}$. This is not trivial since the previous deterministic construction starting from leaves is impossible when the tree is infinite. Moreover, it is not always clear that a solution of the recursive system \eqref{eq:Z_hrecursion} on a given random weighted UBGW tree $\mathcal T$ can be constructed measurably.

We want to construct a random couple $(\mathcal{T},Z(u,v)_{(u,v) \in \overset{\rightarrow}{E}})$ such that the decorated tree is unimodular and the variables $Z$ satisfy the recursive system \eqref{eq:Z_hrecursion}. The unimodularity of the pair requires that the variables $Z$ have the same distribution. In addition, we will restrict ourselves to variables such that, for any vertex $v$ with neighbours $v_1, \ldots , v_k$, the variables $(Z(v,v_i))_{1 \leq i \leq k}$ are independent. This assumption comes from the fact that, in the finite setting, the variables $Z(v,v_i)$ are computed from the disjoint subtrees $T_{(v,v_i)}$.
As consequence, the law of the variables $Z$ must satisfy a recursive distributional equation. Lemma \ref{lem:ZRDE} below guarantees the existence of a solution to this equation.



\begin{lemma}\label{lem:ZRDE}
    Let $N$ be a random variable with law $\hat{\pi}$ and $(w_k)_{k \in \mathbb{N}}$ be a sequence of i.i.d. random variables with law $\omega$, independent of $N$. 
    Then there exists a law $\zeta$ such that for all sequence $(Z_k)_{k \geq 0}
    $ i.i.d of law $\zeta$ and independent of $N$ and of the sequence $(w_k)_{k \in \mathbb{N}}$, the following equality in law holds:
    \[Z_0\overset{\mathrm{(law)}}{=}\max(0,\max_{1 \leq i \leq N} (w_i-Z_i)) .\]
\end{lemma}
\begin{remark} 
We do not claim at this stage that the invariant law $\zeta$ is unique, this will be the purpose of Section~\ref{sec:thuniqueness}.
\end{remark}

\begin{proof}
Let us express the invariance of $\zeta$ in terms of its cumulative distribution function $h$. Let us call $\hat{\phi}$ the generating function of $\hat{\pi}, \hat{\phi}(x)=\sum_{k=0}^{\infty}\mathbb{P}( N=k)x^k=\sum_{k=0}^{\infty}\hat{\pi}_kx^k $. 

Let $t \in \mathbb{R}$, then by definition of $h$ and $Z_i$,
\begin{align*}
h(t)&= \mathbb{P}(Z_0\leq t)= \mathbb{P} \left( \max(0,\max_{1 \leq i \leq N} (w_i-Z_i)) \leq t\right) \\
&=\mathbbm{1}_{t \geq 0} \mathbb{P}\left( \forall i \leq N, w_i - Z_i \leq t  \right)  \\
&= \mathbbm{1}_{t \geq 0} \sum_{k=0}^{\infty} \hat{\pi}_k \mathbb{P}( Z_i\geq w_i-t)^k \\
&= \mathbbm{1}_{t \geq 0} \sum_{k=0}^{\infty} \hat{\pi_k} (1-\mathbb{E}[h(w_1-t) ])^k \\
&= \mathbbm{1}_{t \geq 0 } \hat{\phi}(1- \mathbb{E}[h(w_1-t)]).
\end{align*}
The invariance of $\zeta$ is equivalent to:
\begin{equation} \label{eq:hequation}
\forall t\in\mathbb{R}, \quad h(t)=    \mathbbm{1}_{t \geq 0 } \hat{\phi}(1- \mathbb{E}[h(w_1-t)]).
\end{equation}
The objective is to apply Schauder's fixed point theorem on the subspace $X$ of non-decreasing functions of $\mathbb{R} \mapsto [0,1]$ equipped with the product topology of $\mathbb{R}^{\mathbb{R}}$. The space $X$ is closed and contained in the compact space  $[0,1]^{\mathbb{R}}$ (by Tychonoff's theorem) so $X$ is compact. It is also convex, hence a compact convex subspace of $\mathbb{R}^\mathbb{R}$ that is a Hausdorff topological vector space. 
Let us define the map 
\begin{alignat*}{2}
\mathrm{F}\colon X&\rightarrow X\\
           f&\mapsto \mathrm{F}(f)\colon {}&\mathbb{R} &\rightarrow [0,1] \\
               &&t & \mapsto  {} \mathbbm{1}_{t \geq 0 } \hat{\phi}(1- \mathbb{E}[f(w_1-t)]) . 
\end{alignat*}
By dominated convergence theorem, the map $\mathrm{F}$ is a continuous map from $X$ into itself. 
Hence, Schauder's fixed point theorem implies that there exists a function $f_0 \in X$ such that $\mathrm{F}(f_0)=f_0$. To conclude we just need to verify that the image $\mathrm{F}(X)$ is contained in the subspace of cadlag functions.  
Since elements of $X$ are non-decreasing, the fact that they admit left limits is obvious, we just need to show that any element of $\mathrm{F}(X)$ is right-continuous, that is for any increasing map $f$, the map 
$
t \mapsto   \mathbbm{1}_{t \geq 0 } \hat{\phi}(1- \mathbb{E}[f(w_1-t)])$ is right continuous. 
It is just a consequence of $w_1$ having no atoms, increasing functions having at most countable discontinuities and the dominated convergence theorem.
\end{proof}

\bigskip

Now that we know that solutions to Equation~\eqref{eq:Zlaw} exist, we want to construct unimodular pairs $(\mathcal T,Z)$ for which the variables $Z$ have the same stationary law and satisfy the propagation Equation~\eqref{eq:Z_hrecursion}. More precisely, we will see in the next proposition that for any law $\zeta$ with cumulative distribution function $h$ given by the previous lemma,  a pair $(\mathcal T',Z_h)$ satisfying the following properties can be resampled:
 \begin{itemize}
 \item[i)] The family $Z_h$ is a solution of Equation \eqref{eq:Z_hrecursion} on $\mathcal{T}'$,
 \item[ii)] For all $(u,v) \in \overset{\rightarrow}{E'}$, the variable $Z_h(u,v)$ is distributed according to $\zeta$, in average over $T_{(u,v)}$,
 \item[iii)] The first marginal $\mathcal T'$ has the same law as $\mathcal T$.
 \end{itemize}
Actually, our next proposition states this fact in the more general setting of unimodular decorated weighted UBGW tree. Indeed, this more general setting will be important later.
\begin{prop}\label{prop:constructionMZ_h}
    Let $h$ be a solution to Equation~\eqref{eq:hequation} so that the law $\zeta_h$ with cumulative distribution $h$ is a solution to Equation \eqref{eq:Zlaw}.   Let $\mathbb{T}$ be a unimodular decorated weighted UBGW tree with reproduction law $\pi$ and weights law $\omega$ (for example, $\mathbb{T}=(\mathcal{T},\mathbbm{1}_\mathbb{M})$,
    where $\mathbb M$ is a matching of $\mathcal T$).
\begin{enumerate}[label=(\roman*)]
    \item There exists a unimodular random decorated tree $(\mathbb{T}',Z_h(u,v)_{(u,v) \in \overset{\rightarrow}{E'}})$ such that:
    \begin{itemize}
        \item The law of $\mathbb{T}'$ is the law of $\mathbb{T}$.
        \item  Each $Z_h(u,v)$, in average over $T_{(u,v)}$, has law $\zeta_h$.
        \item For every $(u,v) \in \overset{\rightarrow}{E'}$, one has
    \[ Z_h(u,v)=\max \left\{0,\max_{\substack{u' \sim v \\ u' \neq u}} \left\{  w(v,u')-Z_h(v,u') \right\} \right\} . \]
    \end{itemize}
    \item The pair $(\mathbb{T}',Z_h(u,v)_{(u,v) \in \overset{\rightarrow}{E'}})$ defines a unimodular matching $\mathbb M_h$ on $\mathbb{T}'$ via the rule 
\[
\{u,v\} \in \mathbb{M}_h
\Leftrightarrow 
 Z_h(u,v)+Z_h(v,u) < w(u,v). \]

Furthermore, this "edge-rule" is equivalent to the following "vertex-rule": $v$ is the neighbour of $u$ in $\mathbb{M}_h$ if and only if
\[
v=\argmax_{v' \sim u } (w(u,v')-Z_h(u,v')) \text{ and } w(u,v)-Z_h(u,v)>0. 
\]
In particular, almost surely,
\begin{align*}
    u \text{ is not matched by }\mathbb{M}_h \Leftrightarrow \max_{v' \sim u } (w(u,v')-Z_h(u,v')) \leq 0.
\end{align*}
\end{enumerate}
\end{prop}
\begin{remark}
The proof of the existence of random variables $(\mathbb{T},Z)$ is not constructive and we do not know of a construction for the variables $Z$ that is measurable in terms of the tree and that satisfies the recursion \eqref{eq:Z_hrecursion}, except in specific examples (e.g. finite trees).
\end{remark}

We now turn to the proof of the proposition:

\begin{proof}[Proof of Proposition \ref{prop:constructionMZ_h}]
i)  For an illustration of this proof, we refer to Figure~\ref{fig:constructionZh}. Let $H \in \mathbb{N}^{\star}$, recall $B_H(\mathbb{T})$ the $H-$neighbourhood of the root edge of $\mathbb{T}$. 
We first construct the restriction of $(Z_h(u,v))_{(u,v)\in \overset{\rightarrow}{E}}$ to $B_H(\mathbb{T})$, then observe consistency over $H$, and then apply Kolmogorov's extension theorem.
\begin{figure}
\centering

\includegraphics[scale=1]{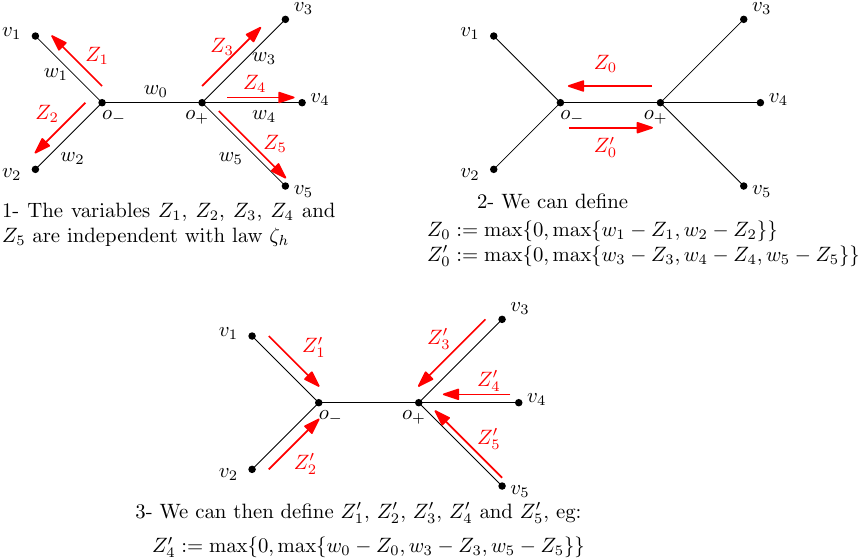}

\caption{Construction of the law of $Z_h$ on a $1-$neighbourhood.}\label{fig:constructionZh}
\end{figure}

We call the depth of a non-oriented edge $\{u,v\}$ its distance to the root edge. 
For $k \in \mathbb{N}$, let $E_k$ be the set of directed edges  $(u,v) \in \overset{\rightarrow}{E}$ such that the depth of $\{u,v\}$ is $k$ and $(u,v)$ is pointing away from the root. We set $(Z_h(u,v))_{(u,v) \in E_H}$ to be independent variables with law $\zeta_h$ as defined in the previous lemma. 
\bigskip

We can then use Recursion~\eqref{eq:Z_hrecursion}: $$Z_h(u,v) = \max(0,\underset{\substack{u' \sim v \\ u' \neq u}}{\max}(w(v,u')-Z_h(v,u')))$$ to define $(Z_h(u,v)))_{(u,v) \in E_{H-1}}$. 
By induction, we define every $Z_h(u,v)$ for $(u,v)$ pointing away from the root and then on the edge-root and its symmetric.
We can then define $Z_h(u,v)$ on the set of directed edges pointing  towards the root of depth $1$ that we call $E_{-1}$, then by induction, we can define them on the set of directed edges pointing towards the root of depth $k$ that we call $E_{-k}$, for $k$ running from $1$ to $H$. 

In this way, we have defined $Z_h(u,v)$ on $B_H(\mathbb{T})$. We can see that because $\zeta_h$ is an invariant law for the RDE, the restriction of $Z_h(u,v)$ to $B_{H-1}(\mathbb{T})$ has the same law as if we defined it directly on $B_{H-1}(\mathbb{T})$. 
By Kolmogorov's extension theorem, we deduce that there exists a process $(\mathbb{T}', Z_h)$ such that the first marginal's law is $\mathbb{T}$ and $Z_h$ satisfies Recursion~\eqref{eq:Z_hrecursion} on $\mathbb{T}'$ with the prescribed law.

It remains to show unimodularity. First, we want to show that the law of $Z_h$ restricted to a $H-$neighbourhood of $(o_+,v)$ with $v$ chosen uniformly among the children of $o_+$ is the same as the law of $Z_h$ restricted to the $H-$neighbourhood of $(o_-,o_+)$.
Since Recursion~\eqref{eq:Z_hrecursion} is preserved, one only needs to show that the exterior variables of  $Z_h$ on the $H-$boundary of $(o_+,v)$ are i.i.d variables of law $\zeta_h$.
The $H-$neighbourhood of $(o_+,v)$ is included in the $(H+1)-$neighbourhood of $(o_-,o_+)$ and thus all the message variables $Z_h$ inside the $H-$neighbourhood of $(o_+,v)$ can be calculated from those of the 
$(H+1)-$neighbourhood of $(o-_,o_+)$. One can check by applying Recursion~\eqref{eq:Z_hrecursion} up to twice from the $(H+1)-$boundary of $(o_-,o_+)$, that we effectively recover variables $Z$ on the $H-$boundary of $(o_+,v)$ that are i.i.d of law $\zeta_h$. The decorated tree
$(\mathbb{T},Z_h)$ is thus stationary. 

Clearly, the law of the variables $Z_h(u,v)$ does not depend on the direction of the root since its finite dimensional statistics are symmetrical on every $B_H(\mathbb{T})$, so changing $(o_-,o_+)$ into $(o_+,o_-)$ does not change their law. Hence, the decorated tree $(\mathbb{T},Z_h)$ is also reversible.

ii) For an illustration of this proof, we refer the reader to Figure~\ref{fig:Zbehaviour}. We need to show that for every vertex $u \in V$, there is at most one neighbouring vertex $v \in V$ satisfying Rule~\eqref{eq:decisionrule}
\[  w(u,v) > Z_h(u,v) + Z_h(v,u).\]
To this end, we will first show the following implication for any $\{u,v\}\in E$
\[ w(u,v) > Z_h(u,v) + Z_h(v,u) \Rightarrow \{v\}= \argmax_{v' \sim u } (w(u,v')-Z_h(u,v')). \]
       
\begin{figure}[t!]
\centering

\tikzset{every picture/.style={line width=0.75pt}} 

\begin{tikzpicture}[x=0.75pt,y=0.75pt,yscale=-1,xscale=1]

\draw    (0,30) -- (120,120) ;
\draw    (230,40) -- (120,120) ;
\draw    (120,120) -- (120,220) ;
\draw [color={rgb, 255:red, 208; green, 2; blue, 27 }  ,draw opacity=1 ]   (90,120) -- (1.58,51.23) ;
\draw [shift={(0,50)}, rotate = 37.87] [color={rgb, 255:red, 208; green, 2; blue, 27 }  ,draw opacity=1 ][line width=0.75]    (10.93,-3.29) .. controls (6.95,-1.4) and (3.31,-0.3) .. (0,0) .. controls (3.31,0.3) and (6.95,1.4) .. (10.93,3.29)   ;
\draw [color={rgb, 255:red, 208; green, 2; blue, 27 }  ,draw opacity=1 ]   (120,110) -- (228.38,31.18) ;
\draw [shift={(230,30)}, rotate = 143.97] [color={rgb, 255:red, 208; green, 2; blue, 27 }  ,draw opacity=1 ][line width=0.75]    (10.93,-3.29) .. controls (6.95,-1.4) and (3.31,-0.3) .. (0,0) .. controls (3.31,0.3) and (6.95,1.4) .. (10.93,3.29)   ;
\draw [color={rgb, 255:red, 208; green, 2; blue, 27 }  ,draw opacity=1 ]   (140,130) -- (140,218) ;
\draw [shift={(140,220)}, rotate = 270] [color={rgb, 255:red, 208; green, 2; blue, 27 }  ,draw opacity=1 ][line width=0.75]    (10.93,-3.29) .. controls (6.95,-1.4) and (3.31,-0.3) .. (0,0) .. controls (3.31,0.3) and (6.95,1.4) .. (10.93,3.29)   ;
\draw [color={rgb, 255:red, 245; green, 166; blue, 35 }  ,draw opacity=1 ]   (90,220) -- (90,132) ;
\draw [shift={(90,130)}, rotate = 90] [color={rgb, 255:red, 245; green, 166; blue, 35 }  ,draw opacity=1 ][line width=0.75]    (10.93,-3.29) .. controls (6.95,-1.4) and (3.31,-0.3) .. (0,0) .. controls (3.31,0.3) and (6.95,1.4) .. (10.93,3.29)   ;
\draw [color={rgb, 255:red, 245; green, 166; blue, 35 }  ,draw opacity=1 ]   (0,20) -- (118.4,108.8) ;
\draw [shift={(120,110)}, rotate = 216.87] [color={rgb, 255:red, 245; green, 166; blue, 35 }  ,draw opacity=1 ][line width=0.75]    (10.93,-3.29) .. controls (6.95,-1.4) and (3.31,-0.3) .. (0,0) .. controls (3.31,0.3) and (6.95,1.4) .. (10.93,3.29)   ;
\draw [color={rgb, 255:red, 245; green, 166; blue, 35 }  ,draw opacity=1 ]   (230,70) -- (151.6,128.8) ;
\draw [shift={(150,130)}, rotate = 323.13] [color={rgb, 255:red, 245; green, 166; blue, 35 }  ,draw opacity=1 ][line width=0.75]    (10.93,-3.29) .. controls (6.95,-1.4) and (3.31,-0.3) .. (0,0) .. controls (3.31,0.3) and (6.95,1.4) .. (10.93,3.29)   ;
\draw  [fill={rgb, 255:red, 0; green, 0; blue, 0 }  ,fill opacity=1 ] (115,120) .. controls (115,117.24) and (117.24,115) .. (120,115) .. controls (122.76,115) and (125,117.24) .. (125,120) .. controls (125,122.76) and (122.76,125) .. (120,125) .. controls (117.24,125) and (115,122.76) .. (115,120) -- cycle ;
\draw  [fill={rgb, 255:red, 0; green, 0; blue, 0 }  ,fill opacity=1 ] (-5,30) .. controls (-5,27.24) and (-2.76,25) .. (0,25) .. controls (2.76,25) and (5,27.24) .. (5,30) .. controls (5,32.76) and (2.76,35) .. (0,35) .. controls (-2.76,35) and (-5,32.76) .. (-5,30) -- cycle ;
\draw  [fill={rgb, 255:red, 0; green, 0; blue, 0 }  ,fill opacity=1 ] (115,220) .. controls (115,217.24) and (117.24,215) .. (120,215) .. controls (122.76,215) and (125,217.24) .. (125,220) .. controls (125,222.76) and (122.76,225) .. (120,225) .. controls (117.24,225) and (115,222.76) .. (115,220) -- cycle ;
\draw  [fill={rgb, 255:red, 0; green, 0; blue, 0 }  ,fill opacity=1 ] (225,40) .. controls (225,37.24) and (227.24,35) .. (230,35) .. controls (232.76,35) and (235,37.24) .. (235,40) .. controls (235,42.76) and (232.76,45) .. (230,45) .. controls (227.24,45) and (225,42.76) .. (225,40) -- cycle ;

\draw (122,123.4) node [anchor=north west][inner sep=0.75pt]    {$u$};
\draw (48,77.4) node [anchor=north west][inner sep=0.75pt]    {$w_{1}$};
\draw (171,80.4) node [anchor=north west][inner sep=0.75pt]    {$w_{2}$};
\draw (98,162.4) node [anchor=north west][inner sep=0.75pt]    {$w_{3}$};
\draw (12,-2.6) node [anchor=north west][inner sep=0.75pt]    {$\textcolor[rgb]{0.96,0.65,0.14}{Z}\textcolor[rgb]{0.96,0.65,0.14}{_{1}^{'}} =\max( 0,w_{2} -\textcolor[rgb]{0.82,0.01,0.11}{Z}\textcolor[rgb]{0.82,0.01,0.11}{_{2}})$};
\draw (173,110) node [anchor=north west][inner sep=0.75pt]    {$ \begin{array}{l}
\textcolor[rgb]{0.96,0.65,0.14}{Z}\textcolor[rgb]{0.96,0.65,0.14}{_{2}^{'}} =\max(0, w_{1} -\textcolor[rgb]{0.82,0.01,0.11}{Z}\textcolor[rgb]{0.82,0.01,0.11}{_{1}})
\end{array}$};
\draw (-79,167.4) node [anchor=north west][inner sep=0.75pt]    {$\textcolor[rgb]{0.96,0.65,0.14}{Z}\textcolor[rgb]{0.96,0.65,0.14}{_{3}^{'}} =\max( 0,w_{1} -\textcolor[rgb]{0.82,0.01,0.11}{Z}\textcolor[rgb]{0.82,0.01,0.11}{_{1}})$};
\draw (170,140.4) node [anchor=north west][inner sep=0.75pt]{
$\begin{array}{l}
\text{Suppose } w_{1} -\textcolor[rgb]{0.82,0.01,0.11}{Z}\textcolor[rgb]{0.82,0.01,0.11}{_{1}} \  >\ w_{2} -\textcolor[rgb]{0.82,0.01,0.11}{Z}\textcolor[rgb]{0.82,0.01,0.11}{_{2}}  >w_{3} -\textcolor[rgb]{0.82,0.01,0.11}{Z}\textcolor[rgb]{0.82,0.01,0.11}{_{3}}. \\
\text{We can see that the} \ \textcolor[rgb]{0.96,0.65,0.14}{Z_{i}^{'}} \ \text{pointing towards} \ u\\
\text{are all equal to }
\max_{v\sim u}( 0, w( u,v)-\textcolor[rgb]{0.82,0.01,0.11}{Z}( u,v))\\
\text{except for the } v\ \text{hitting that maximum}\\
\text{in which case the value is the second largest.}
\end{array}$};
\draw (21,102.4) node [anchor=north west][inner sep=0.75pt]    {$\textcolor[rgb]{0.82,0.01,0.11}{Z_{1}}$};
\draw (159,32.4) node [anchor=north west][inner sep=0.75pt]    {$\textcolor[rgb]{0.82,0.01,0.11}{Z}\textcolor[rgb]{0.82,0.01,0.11}{_{2}}$};
\draw (149,170.4) node [anchor=north west][inner sep=0.75pt]    {$\textcolor[rgb]{0.82,0.01,0.11}{Z}\textcolor[rgb]{0.82,0.01,0.11}{_{3}}$};
\draw (-18,2.4) node [anchor=north west][inner sep=0.75pt]    {$v_{1}$};
\draw (239,22.4) node [anchor=north west][inner sep=0.75pt]    {$v_{2}$};
\draw (102,222.4) node [anchor=north west][inner sep=0.75pt]    {$v_{3}$};

\end{tikzpicture}

\caption{Illustration of the behaviour of the $Z_h$ around a vertex $u$.}\label{fig:Zbehaviour}

\end{figure}
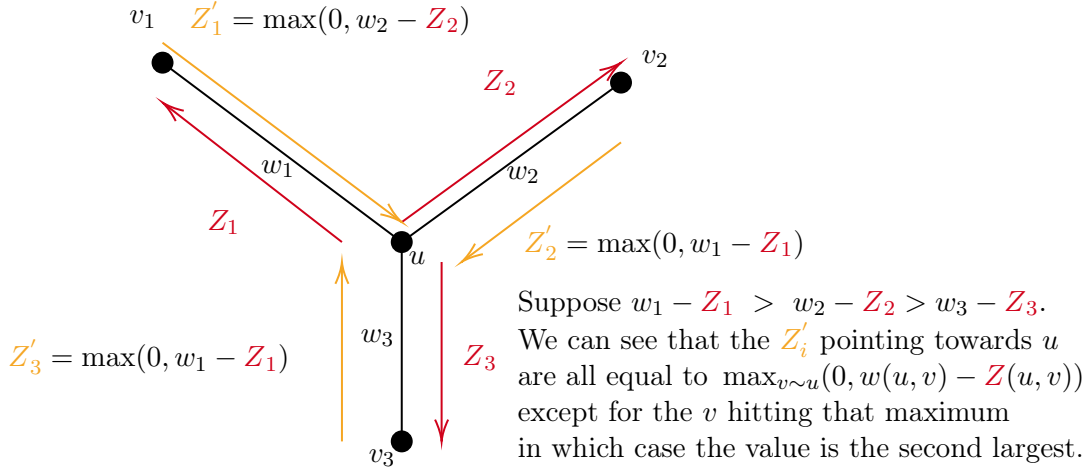
Indeed:
\begin{align*}
 & \, w(u,v)> Z_h(u,v)+Z_h(v,u)  \\
 \Rightarrow &\,  w(u,v)-Z_h(u,v) > Z_h(v,u)= \max(0,\max_{\substack{ v' \sim u \\ v' \neq v}} (w(u,v')-Z_h(u,v'))) .
\end{align*}
This implies that for every neighbour of $u$, $v'$ different from $v$, $w(u,v)-Z_h(u,v) > w(u,v')-Z_h(u,v')$, which in turn implies that $v$ is the unique element of $\argmax_{v' \sim u } (w(u,v')-Z_h(u,v'))$.
Hence the vertex $u$ is matched to at most one neighbour.

To show the equivalent vertex rule, we only need to show that $u$ is not matched if and only if $\max_{v' \sim u } (w(u,v')-Z_h(u,v')) \leq 0$.

For the if part, we have that any incident $Z_h(v,u)$ for $v \sim u$ evaluates to zero by applying Recursion~\eqref{eq:Z_hrecursion}. So for any $v \sim u$, $Z_h(v,u)+Z_h(u,v)=Z_h(u,v)$.
Since $\max_{v' \sim u } (w(u,v')-Z_h(u,v')) \leq 0$, we deduce that for any $v' \sim u$, $w(u,v')\leq Z_h(u,v')$. But we have just seen that $Z_h(v',u)=0$, hence $w(u,v')\leq Z_h(u,v')+Z_h(v',u)$ for any $v'\sim u$. So $v'$ is not matched to $u$.

Reciprocally, if $u$ is not matched by $\mathbb{M}_h$, then for any $v' \sim u$ \[ \max(0,\max_{\substack{v' \sim u \\ v' \neq v}}(w(u,v')-Z_h(u,v')))=Z_h(v',u)\geq w(u,v')-Z_h(u,v'), \]  so the maximum of $\argmax_{v' \sim u } (w(u,v')-Z_h(u,v'))$ is either not reached by any $v' \sim u$ or there are several $v'$ in the $\argmax$. Since $\omega$ is non-atomic, there is almost surely no ties and the $\max$ is not reached. So $\max_{v' \sim u }(w(u,v')-Z_h(u,v')\leq 0$ which concludes.
\end{proof}

From now on, despite resampling $\mathbb{T}$ when applying Proposition~\ref{prop:constructionMZ_h}, we will forgo the notation $\mathbb{T}'$ for $\mathbb{T}$ to ease notations.

\bigskip

We conclude this section with a technical, but useful, result. It gives a property of the probability distributions that are solutions of Equation~\eqref{eq:Zlaw}.  

\begin{lemma}\label{lem:Zatomic}
     Any law $\zeta_h$ defined as a solution of Equation~\eqref{eq:Zlaw} is atomic at zero and only at zero.
\end{lemma}
\begin{proof}
Let $h$ be any solution to Equation~\eqref{eq:hequation}.
First, let us show that $\zeta_h$ is atomic at zero.
Evaluate Equation~$\eqref{eq:hequation}$ at $0$:
\[ h(0)=\hat{\phi}(1-\mathbb{E}_{W \sim \omega}[h(W)]  )  .\]

The result if immediate if $\hat{\phi}(0)>0$, hence we assume $\hat{\phi}(0)=0$ and $h(0)=0$.  Since $\hat{\phi}>0$ on $]0,1]$, we deduce that $\mathbb{E}_{W \sim \omega}[h(W)]=1$. 
This implies in particular that $\inf(\mathrm{supp}(\omega)) \geq \sup(\mathrm{supp}(\zeta_h))$. Since $\mathrm{supp}(\zeta_h) \subset \mathbb{R}_+$, this implies that $\mathrm{supp}(\omega) \subset \mathbb{R}_+$ as well.

Recall equation~\eqref{eq:hequation}:
\begin{equation*}
    h(t)=\mathbbm{1}_{t \geq 0} \hat{\phi}\left(1-\mathbb{E}_{W\sim\omega}\left[h(W-t)\right] \right).
\end{equation*}
Evaluate this equation at some variable $W'$ of law $\omega$, this yields:
\[ h(W')=\mathbbm{1}_{W' \geq 0} \hat{\phi}\left(1-\mathbb{E}_{W\sim\omega}\left[h(W-W')\right] \right)   \]
We have shown previously that $W'\geq 0 $ almost surely so we can discard $\mathbbm{1}_{W' \geq 0}$. Now take expectation with respect to $W'$, this yields:
\[ 1= \mathbb{E}_{W' \sim \omega} \left[\hat{\phi}\left(1-\mathbb{E}_{W\sim\omega}\left[h(W-W')\right]\right) \right] . \]
But since $\hat{\phi}<1$ outside of $1$, this implies that almost surely on $W'$, 
\[ \mathbb{E}_{W\sim \omega}[h(W-W')]=0 .\]
Let $W',W''$ independent variables of law $\omega$. Evaluate equation~\eqref{eq:hequation} on $W'-W''$:
\[  h(W'-W'')=\mathbbm{1}_{W'-W'' \geq 0} \hat{\phi}\left(1-\mathbb{E}_{W\sim\omega}\left[h(W+W''-W')\right] \right)  . \]
Taking the expectation on both sides, 
\[0= \mathbb{E}\left[ \mathbbm{1}_{W'-W'' \geq 0} \hat{\phi}\left(1-\mathbb{E}_{W\sim\omega}\left[h(W-(W'-W''))\right] \right)\right].\]
So almost surely, either $W'-W''<0$ or $\mathbb{E}_{W\sim\omega}[h(W-(W'-W''))]=1$.
Conditionally on $W'-W''\geq 0$, $\mathbb{E}_{W\sim\omega}[h(W-(W'-W''))]=1$.  But the function $h$ is non-decreasing so this statement stays true for $W'-W''<0$.

So now, taking expectations on $W'-W''$, we obtain:
\[ \mathbb{E}(h(W+W''-W'))=1 .  \]

We can now iterate this procedure. Let $(W_{i})_{i \in \mathbb{N}^{\star}}$ be i.i.d variables with law $\omega$, we will show the following statement by induction:
\begin{equation}\label{eq:Zatomezero}
    \forall n \in \mathbb{N}^{\star}, \mathbb{E}\left[h\left(\sum_{i=1}^{n}W_{2i-1}-\sum_{i=1}^{n-1}W_{2i}\right)  \right]=1.
\end{equation}
The statement holds for $n=2$. Assume it holds for some $n \in \mathbb{N}^{\star}$.
Evaluating $h$ at $\sum_{i=1}^{n}W_{2i-1}-\sum_{i=1}^{n-1}W_{2i}$, we get
\[ h\left(\sum_{i=1}^{n}W_{2i-1}-\sum_{i=1}^{n-1}W_{2i}\right)=\mathbbm{1}_{\sum_{i=1}^{n}W_{2i-1}-\sum_{i=1}^{n-1}W_{2i} \geq 0} \hat{\phi}\left(1-\mathbb{E}\left[h\left(\sum_{i=1}^{n}W_{2i}-\sum_{i=1}^{n}W_{2i-1}\right)\right] \right) .  \]
where the expectation is taken only over $W_{2n}$.

Now take expectation over $(W_i)_{i \leq 2n-1}$, by the induction hypothesis the left-hand side will evaluate to one:
\[1= \mathbb{E}\left[\mathbbm{1}_{\sum_{i=1}^{n}W_{2i-1}-\sum_{i=1}^{n-1}W_{2i} \geq 0} \hat{\phi}\left(1-\mathbb{E}\left[h\left(\sum_{i=1}^{n}W_{2i}-\sum_{i=1}^{n}W_{2i-1}\right)\right] \right) \right].   \]
We deduce that almost surely on $(W_i)_{i \leq 2n-1}$,
\[\mathbb{E}_{W_{2n} \sim \omega}\left[ h\left(\sum_{i=1}^{n}W_{2i}-\sum_{i=1}^{n}W_{2i-1}\right)\right]=0. \]
Now evaluate $h$ on $\sum_{i=1}^{n}W_{2i}-\sum_{i=1}^{n}W_{2i-1}$, this yields:
\begin{align*}
h &\left(\sum_{i=1}^{n}W_{2i}-\sum_{i=1}^{n}W_{2i-1}\right) \\
&\quad=\mathbbm{1}_{\sum_{i=1}^{n}W_{2i}-\sum_{i=1}^{n}W_{2i-1} \geq 0}\hat{\phi}\left(1-\mathbb{E}\left[h\left(W_{2n+1}-\left(\sum_{i=1}^{n}W_{2i}-\sum_{i=1}^{n}W_{2i-1}\right)\right)\right] \right)  ,
\end{align*}
where the expectation is over $W_{2n+1}$.
Now take expectation over $(W_i)_{i \leq 2n}$, the left-hand side evaluates to 0 as we have just shown:
\[0= \mathbb{E}\left[\mathbbm{1}_{\sum_{i=1}^{n}W_{2i}-\sum_{i=1}^{n}W_{2i-1} \geq 0}\hat{\phi}\left(1-\mathbb{E}\left[h\left(W_{2n+1}-\left(\sum_{i=1}^{n}W_{2i}-\sum_{i=1}^{n}W_{2i-1}\right)\right)\right] \right) \right]. \]
We deduce that:
\[ \mathbb{P}\left( \mathbb{E}_{W_{2n+1}\sim\omega}\left[h\left(W_{2n+1}-\left(\sum_{i=1}^{n}W_{2i}-\sum_{i=1}^{n}W_{2i-1}\right)  \right) \right]=1 \middle| \sum_{i=1}^{n}W_{2i}-\sum_{i=1}^{n}W_{2i-1} \geq 0    \right)=1.    \]
Again, since $h$ is non-decreasing, we can discard the conditioning to obtain:
\[ \mathbb{P}\left( \mathbb{E}_{W_{2n+1}\sim\omega}\left[h\left( \sum_{i=1}^{n+1}W_{2i-1}-\sum_{i=1}^{n}W_{2i}  \right) \right]=1 \right)=1.    \]
So taking expectation over $(W_{i})_{i \leq 2n+1}$, we have proven that:
\[\mathbb{E}\left[h\left( \sum_{i=1}^{n+1}W_{2i-1}-\sum_{i=1}^{n}W_{2i}  \right) \right]=1, \]
which completes the induction.
Now, introducing a variable $Z_h$ of law $\zeta_h$ independent from the $(W_i)_{1 \leq i \leq 2n+1}$, Equation~\eqref{eq:Zatomezero} can be rewritten as:
\[ \forall n \in \mathbb{N}^*,   \mathbb{P}\left( Z_h \leq \sum_{i=1}^{n+1} W_{2i-1}  - \sum_{i=1}^{n}W_i \right) =1    .\]
This implies that:
\[ \forall n \in \mathbb{N}^{\star}, \inf(\mathrm{supp}(\omega)) - n\left[\sup(\mathrm{supp}(\omega))-\inf(\mathrm{supp}(\omega))\right] \geq \sup(\mathrm{supp}(\zeta_h))  .\]
But since $\omega$ is atomless, $\mathrm{supp}(\omega)$ is not reduced to a single point so $\sup(\mathrm{supp}(\omega))-\inf(\mathrm{supp}(\omega))>0$. Taking $n \rightarrow \infty$ then shows that $-\infty \geq \sup(\mathrm{supp}(\zeta_h))$. This is clearly impossible as $\zeta_h$ is a law on $\mathbb{R}_+$. We then conclude that $h(0)>0$ which exactly means that $\zeta_h$ is atomic at zero.

Now let us observe that $\zeta$ is non-atomic outside of zero. This is essentially a consequence of the RDE, for $Z,Z_i \sim \zeta$ and $w_i\sim \omega$, all mutually independent:
\begin{equation*}
    Z\overset{\mathcal{L}}{=}\max\left(0,\max_{1\leq i \leq N}(w_i-Z_i)\right).
\end{equation*}
As $w_i$ is non-atomic, $w_i-Z_i$ is non-atomic as well, so the right-hand side is non-atomic when it doesn't evaluate to zero, which concludes.
\end{proof}

\subsection{Adding self-loops}\label{sec:selfloops}
In the remainder of this section, we will adopt the vertex-rooted point of view. Dealing with partial matchings rather than perfect matchings will be very cumbersome in future proofs. For example, when comparing two matchings, we would need to discuss several cases depending on whether a vertex is matched or not for the two matchings. A simple solution to deal with this is to add self loops, in which case the graphs always have perfect matchings, and any partial matching on the graph without self loops can be augmented into a perfect matching on the graph with self loops. Reciprocally, a perfect matching on a graph with self loops can be restricted to a partial matching on the corresponding graph without self loops. 

We need to choose adequate weights for the self loops so that a partial matching on the graph without self loops is optimal among all partial matchings of the original graph if and only if the corresponding perfect matching on the augmented graph with self loops is optimal among all perfect matchings of the augmented graph.

Starting with a system of variables $Z(u,v)$ for $(u,v) \in \overset{\rightarrow}{E}$ that satisfies recursion \eqref{eq:Z_hrecursion}, we want to extend it to variables $Z(v,v)$ for each self loop so that we get a similar recursion as \eqref{eq:Z_hrecursion} and the decision rule \eqref{eq:decisionrule} constructs the perfect matching of the augmented graph associated to the matching defined by the original variables.
Recall that, for $u \neq v$ and $u \sim v$, the original recursion \eqref{eq:Z_hrecursion} on the variables $Z$ reads
\[
Z(u,v)= \max(0 , \max_{\substack{u' \sim v \\ u' \neq u,v}}(w(v,u')-Z(v,u'))).
\]
To take into account the self loop $(v,v)$, we want instead
\[
Z(u,v)= \max(0 , \max_{\substack{u' \sim v \\ u' \neq u}}(w(v,u')-Z(v,u'))).
\]
Choosing
\[
w(v,v) = Z(v,v)
\]
ensures that this recursion is still valid and reads
\[
Z(u,v)= \max_{\substack{u' \sim v \\ u' \neq u}}(w(v,u')-Z(v,u')).
\]
For $Z(v,v)$ to satisfy the same recursion, we set
\begin{align*}
w(v,v) = Z(v,v) := \max_{\substack{u' \sim v \\ u' \neq v}}(w(v,u')-Z(v,u')).
\end{align*}
This definition of $Z(v,v)$ has a simple interpretation in the finite setting: it is the (signed) gain if one forces $v$ to be matched to one of its neighbours.

\bigskip

\begin{figure}[t!]
\centering
\tikzset{every picture/.style={line width=0.75pt}} 
\begin{tikzpicture}[x=0.75pt,y=0.75pt,yscale=-1,xscale=1]

\draw    (110,120) -- (110,220) ;
\draw    (110,20) -- (110,120) ;
\draw    (110,120) -- (10,120) ;
\draw    (210,120) -- (110,120) ;
\draw [color={rgb, 255:red, 208; green, 2; blue, 27 }  ,draw opacity=1 ]   (120,90) -- (120,22) ;
\draw [shift={(120,20)}, rotate = 90] [color={rgb, 255:red, 208; green, 2; blue, 27 }  ,draw opacity=1 ][line width=0.75]    (10.93,-3.29) .. controls (6.95,-1.4) and (3.31,-0.3) .. (0,0) .. controls (3.31,0.3) and (6.95,1.4) .. (10.93,3.29)   ;
\draw [color={rgb, 255:red, 208; green, 2; blue, 27 }  ,draw opacity=1 ]   (140,130) -- (208,130) ;
\draw [shift={(210,130)}, rotate = 180] [color={rgb, 255:red, 208; green, 2; blue, 27 }  ,draw opacity=1 ][line width=0.75]    (10.93,-3.29) .. controls (6.95,-1.4) and (3.31,-0.3) .. (0,0) .. controls (3.31,0.3) and (6.95,1.4) .. (10.93,3.29)   ;
\draw [color={rgb, 255:red, 208; green, 2; blue, 27 }  ,draw opacity=1 ]   (100,150) -- (100,218) ;
\draw [shift={(100,220)}, rotate = 270] [color={rgb, 255:red, 208; green, 2; blue, 27 }  ,draw opacity=1 ][line width=0.75]    (10.93,-3.29) .. controls (6.95,-1.4) and (3.31,-0.3) .. (0,0) .. controls (3.31,0.3) and (6.95,1.4) .. (10.93,3.29)   ;
\draw [color={rgb, 255:red, 208; green, 2; blue, 27 }  ,draw opacity=1 ]   (80,110) -- (12,110) ;
\draw [shift={(10,110)}, rotate = 360] [color={rgb, 255:red, 208; green, 2; blue, 27 }  ,draw opacity=1 ][line width=0.75]    (10.93,-3.29) .. controls (6.95,-1.4) and (3.31,-0.3) .. (0,0) .. controls (3.31,0.3) and (6.95,1.4) .. (10.93,3.29)   ;
\draw    (140,90) .. controls (170.45,59.8) and (174.45,118.8) .. (110,120) ;
\draw    (140,90) -- (111.41,118.59) ;
\draw [shift={(110,120)}, rotate = 315] [color={rgb, 255:red, 0; green, 0; blue, 0 }  ][line width=0.75]    (10.93,-3.29) .. controls (6.95,-1.4) and (3.31,-0.3) .. (0,0) .. controls (3.31,0.3) and (6.95,1.4) .. (10.93,3.29)   ;

\draw (161,140.4) node [anchor=north west][inner sep=0.75pt]    {$w_{1} ,\textcolor[rgb]{0.82,0.01,0.11}{Z_{1}}$};
\draw (41,200.4) node [anchor=north west][inner sep=0.75pt]    {$w_{2} ,\textcolor[rgb]{0.82,0.01,0.11}{Z}\textcolor[rgb]{0.82,0.01,0.11}{_{2}}$};
\draw (11,82.4) node [anchor=north west][inner sep=0.75pt]    {$w_{3} ,\textcolor[rgb]{0.82,0.01,0.11}{Z}\textcolor[rgb]{0.82,0.01,0.11}{_{3}}$};
\draw (131,22.4) node [anchor=north west][inner sep=0.75pt]    {$w_{4} ,\textcolor[rgb]{0.82,0.01,0.11}{Z}\textcolor[rgb]{0.82,0.01,0.11}{_{4}}$};
\draw (166,52.4) node [anchor=north west][inner sep=0.75pt]    {$ \begin{array}{l}
w^{s}_h( v,v) :=Z_{h}^{s}( v,v) \ 
:=\underset{i}{\max}( w_{i} -\textcolor[rgb]{0.82,0.01,0.11}{Z}_{\textcolor[rgb]{0.82,0.01,0.11}{i}}).
\end{array}$};
\draw (112,123.4) node [anchor=north west][inner sep=0.75pt]    {$v$};

\end{tikzpicture}
\caption{Extension to self-loops.}\label{fig:selfloops}
\end{figure}
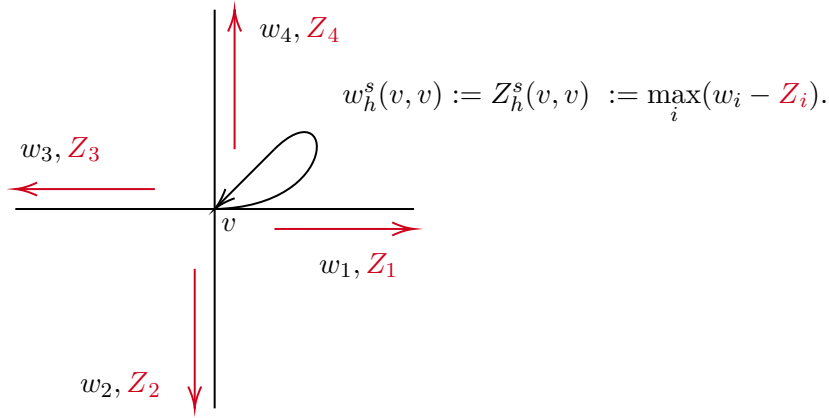

Formally, we define the extension
\begin{equation} \label{eq:Tsh}
(\mathbb{T}^{s}_h,Z_h^{s},(f_i^s)_{i \in \{0,...,I\}})=((V^{s},E^{s}),o,w^{s}_h,Z_h^{s},(f_i^s)_{i \in \{0,....,I\}})
\end{equation}
deterministically on every outcome of $(\mathbb{T},Z_h)=((V,E),o,w,Z_h, (f_i)_{i \in \{0,...,I\}})$, refer to Figure~\ref{fig:selfloops} for an illustration.
The vertex set, the root and the decorations remain unchanged, $V^s=V$ $o^s=o$ and $f_i^s$ restricted to $\overset{\rightarrow}{E}$ is $f_i$ and $+\infty$ otherwise.  We will enrich the edge set by adding self-loops $(v,v)$ for every $v \in V$ to construct $E^s$.
If $e=(u,v)$ with $u \neq v$, we do not change $w$ nor $Z_h$: we set $w^{s}_h(u,v)=w(u,v)$ and $Z_h^{s}(u,v)=Z_h(u,v)$. 

If $e=(v,v)$, we set
\begin{align*}
Z_h^{s}(v,v)&= \max_{\substack{u' \sim v \\ u' \neq v}}(w(v,u')-Z_h(v,u')),\\ 
w^{s}_h(v,v)&=Z_h^{s}(v,v) .
\end{align*}
The following recursive equation holds for $(\mathbb{T}^s_h,Z_h^{s})$: 
\begin{equation}\label{eq:recursionZ_hS}
 \forall (u,v) \in \overset{\rightarrow}{E^s}, \quad Z_h^{s}(u,v)= \max_{\substack{u' \sim v \\ u' \neq u}}(w^{s}_h(v,u')-Z_h^{s}(v,u'))  . 
\end{equation}

\begin{remark}
This is the same recursion that has been studied by Aldous \cite{aldous2000zeta2} in the setting of the complete graph and the Poisson Weighted Infinite Tree.
\end{remark}

Now, let us look at the analogous decision to the decision rule~\eqref{eq:decisionrule}:

\[\{u,v\} \in \mathbb{M}_h^s \Leftrightarrow w^s_h(u,v)> Z_h^s(u,v)+Z_h^s(v,u).
\]
and show that $\mathbb M_h^s$ is the perfect matching of $\mathbb T^s_h$ corresponding to $\mathbb M_h$.

If $(u,v)$ is not a self loop, the rule for $\mathbb M_h^s$ is the same as for $\mathbb M_h$ since we have not modified the variables appearing in the equation. Therefore, the only extra edges of $\mathbb{M}_h^s$ are self loops.

Now we need to show that a self loop $(v,v) \in \mathbb M_h^s$ iff the vertex $v$ is not matched in $\mathbb M_h$.
By Proposition~\ref{prop:constructionMZ_h} ii), $v$ is unmatched by $\mathbb{M}_h$ if and only if $\max(0,\max_{v' \sim v}w(v,v')-Z_h(v,v'))=0$. Since we have defined $Z_h^s(v,v)$ as $\max_{v' \sim v}w(v,v')-Z_h(v,v')$, we deduce that $v$ is unmatched by $\mathbb{M}_h$ if and only if $w^s_h(v,v)=Z_h^s(v,v)\leq0$.
Since the $w$'s are non-atomic, we can replace the previous inequality with a strict one. On the other hand, the rule for self loops reads
\[\{v,v\} \in \mathbb{M}_h^s \Leftrightarrow w^s_h(v,v)> Z_h^s(v,v)+Z_h^s(v,v) = 2 \, w^s_h(v,v),
\]
which is equivalent to $w^s_h(v,v) < 0$ and to $v$ not being matched in $\mathbb M_h$.

\bigskip

Furthermore we have also found the equivalent vertex-rule for $\mathbb{M}_h^s$:
\begin{equation}\label{eq:critereZhS}
\{u,v\} \in \mathbb{M}_h^{s} \Leftrightarrow w^{s}_h(u,v)> Z_h^{s}(u,v)+Z_h^{s}(v,u) \Leftrightarrow v=\argmax_{v' \sim v } (w^s_h(u,v')-Z_h^s(u,v')) .
\end{equation}
The following property we have just shown will be key later:
\begin{equation}\label{eq:critereselfloops}
(v,v) \in \mathbb{M}_h^{s} \Leftrightarrow Z_h^{s}(v,v)=w^{s}_h(v,v)<0. 
\end{equation}
Intuitively, this corresponds to the fact that $u$ is matched to one of its neighbours if and only if the gain of forcing $u$ to be matched to one of its neighbours is positive. In other words, $u$ is matched by $\mathbb{M}_h$ if and only if $w^s_h(u,u)=Z^s(u,u)\geq 0$. However, this does not give information about who is the chosen neighbour of $u$. 

\bigskip

We will later use the fact that introducing self loops with this specific choice of weights conserves unimodularity:

\begin{lemma} \label{lem:Tsunimodular}
    Let $(\mathbb{T},(f_i)_{i \in \{0,...,I\}})$ be a decorated weighted (vertex rooted) UBGW random tree.
    Fix $h$ a solution to Equation~\eqref{eq:hequation} and recall the definition of $\mathbb T^s_h$ in Equation~\eqref{eq:Tsh}. Then $(\mathbb{T}^{s}_h,Z_h^{s},(f_i^s)_{i \in \{0,...,I\}})$ is unimodular.
\end{lemma}

\begin{proof}
    We need to show that the mass-transport principle still holds for $(\mathbb{T}^{s}_h,(f_i^s)_{i \in \{0,...,I\}})$. The idea is that $(\mathbb{T}^{s}_h,(f_i^s)_{i \in \{0,...,I\}})$ is the image of $(\mathbb{T},Z_h,f_i)$ through a deterministic bijective operator that doesn't modify vertices. For every $f$ that is a measurable function over the space of doubly pointed decorated trees, we want to show that:
\[ \mathbb{E}_{(T,Z,f_i) \sim (\mathbb{T}^s_h,Z_h^s,f_i^s)}\left[\sum_{v\in V^{s}} f(T,Z_h,(f_i),o,v) \right]
=\mathbb{E}_{(T,Z,f_i) \sim (\mathbb{T}^s_h,Z_h^s,f_i^s)}\left[\sum_{v\in V^{s}} f(T,Z_h,(f_i),v,o) \right]    .\]
To this end we will write $(T^s,Z_h^s,f_i^s) \sim s(\hat{T},Z_h,f_i))$ where $(\hat{T},\hat{Z_h},(\hat{f_i^s})) \sim (\mathbb{T},Z_h)$, and notice that neither the vertex set nor the root have changed so the required equality rewrites as:
 \[ \mathbb{E}_{(T,Z_h,f_i) \sim (\mathbb{T},Z_h,(f_i)}\left[\sum_{v\in V} f(s(T,Z_h,(f_i)),o,v) \right]
=\mathbb{E}_{(T,Z_h,f_i) \sim (\mathbb{T},Z_h,(f_i)}\left[\sum_{v\in V} f(s(T,Z_h,(f_i)),v,o) \right].\]
which is true by unimodularity of $(\mathbb{T},Z_h,(f_i))$.
\end{proof}

\bigskip

Finally, we define a new performance on the space of matchings on $\mathbb{T}^{s}_h$, we only consider perfect matchings $\mathbb{M}^{s}$ where every vertex has exactly one neighbour.   
\begin{equation}
\perf^{s}_V(\mathbb{T}^s_h,\mathbb{M}^{s}) =  \mathbb{E}\left[ \sum_{v \sim o}w^{s}_h(o,v)\mathbbm{1}_{(o,v) \in\mathbb{M}^s}  \right] .
\end{equation}
The aim of the next section is to first prove that $\mathbb{M}_h^s$ always maximises $\perf^{s}_V$ among unimodular perfect matchings of $(\mathbb T^s_h,w^s_h)$. As a consequence, as can be guessed from \eqref{eq:critereselfloops}, we will then show that $\mathbb{M}_h$ is optimal among unimodular matchings of $(\mathbb T,w)$.


\subsection{Optimality}\label{sec:optimality}

The goal of this subsection is to prove the following statement:

\begin{prop}\label{prop:optimal}
   For any  solution $h$ to Equation~\eqref{eq:hequation}, the matching  $(\mathbb{T},\mathbb{M}_h)$ is optimal on $\mathbb{T}$. Furthermore, almost surely,  any unimodular optimal matching $\mathbb{M}$ on $(\mathbb{T},Z_h,\mathbb{M}_h)$ shares the same set of unmatched vertices.
\end{prop}
Before stating the proof, let us stress the fact that given any unimodular matching $\mathbb{M}$ on $\mathbb{T}$, Proposition \ref{prop:constructionMZ_h} ii) states that  $\mathbb{M}_h$ can be resampled on top of $(\mathbb{T},\mathbb{M})$ such that $(\mathbb{T},\mathbb{M},\mathbb{M}_h)$ is unimodular, and both matchings co-exist on the tree while retaining unimodularity. 

We will proceed by proving the following lemma:
\begin{lemma}\label{lem:perfS}
     Fix $h$ a solution to Equation~\eqref{eq:hequation} and recall the definition of $\mathbb T^s_h$ in Equation~\eqref{eq:Tsh}. Given $(\mathbb T, Z_h)$, the matching $\mathbb{M}_h^{s}$ maximises $\perf^{s}_V$ among all unimodular perfect matchings of $\mathbb T^s_h$. 
\end{lemma}
\begin{proof}[Proof of Lemma \ref{lem:perfS}]
Fix $\mathbb{M}^{s}$ any unimodular matching on $\mathbb{T}^s_h$. And sample $(\mathbb{T}^s_h,\mathbb{M}^s,\mathbb{M}_h^s)$ as in Proposition \ref{prop:constructionMZ_h} ii).

We define the neighbour function $n$ associated with $\mathbb{M}^s$ for which $n(u)$ is the vertex matched with $u$ in $\mathbb{M}^s$. 
Similarly, we define the neighbour function $n_h$ associated with $\mathbb{M}_h^s$. 

Set $v_0=o$, $v_{1}=n(o)$, $v_{-1}=n_h(o)$, $v_{-2}=n(v_{-1})=n(n_h(o))$, see Figure~\ref{fig:constructionMZ_h} for an illustration. 
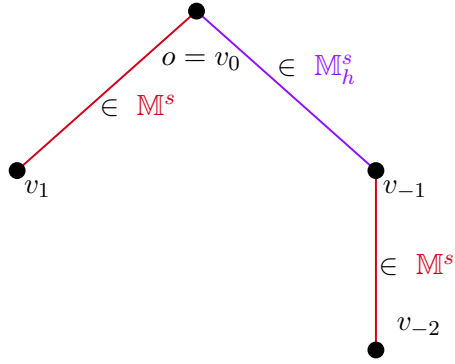
\begin{figure}[t!]
\centering    
\tikzset{every picture/.style={line width=0.75pt}} 

\begin{tikzpicture}[x=0.75pt,y=0.75pt,yscale=-1,xscale=1]

\draw [color={rgb, 255:red, 208; green, 2; blue, 27 }  ,draw opacity=1 ]   (100,0) -- (10,80) ;
\draw [color={rgb, 255:red, 144; green, 19; blue, 254 }  ,draw opacity=1 ]   (100,0) -- (190,80) ;
\draw [color={rgb, 255:red, 208; green, 2; blue, 27 }  ,draw opacity=1 ]   (190,80) -- (190,170) ;
\draw  [fill={rgb, 255:red, 0; green, 0; blue, 0 }  ,fill opacity=1 ] (103.73,0) .. controls (103.73,-2.06) and (102.06,-3.73) .. (100,-3.73) .. controls (97.94,-3.73) and (96.27,-2.06) .. (96.27,0) .. controls (96.27,2.06) and (97.94,3.73) .. (100,3.73) .. controls (102.06,3.73) and (103.73,2.06) .. (103.73,0) -- cycle ;
\draw  [fill={rgb, 255:red, 0; green, 0; blue, 0 }  ,fill opacity=1 ] (13.73,80) .. controls (13.73,77.94) and (12.06,76.27) .. (10,76.27) .. controls (7.94,76.27) and (6.27,77.94) .. (6.27,80) .. controls (6.27,82.06) and (7.94,83.73) .. (10,83.73) .. controls (12.06,83.73) and (13.73,82.06) .. (13.73,80) -- cycle ;
\draw  [fill={rgb, 255:red, 0; green, 0; blue, 0 }  ,fill opacity=1 ] (193.73,80) .. controls (193.73,77.94) and (192.06,76.27) .. (190,76.27) .. controls (187.94,76.27) and (186.27,77.94) .. (186.27,80) .. controls (186.27,82.06) and (187.94,83.73) .. (190,83.73) .. controls (192.06,83.73) and (193.73,82.06) .. (193.73,80) -- cycle ;
\draw  [fill={rgb, 255:red, 0; green, 0; blue, 0 }  ,fill opacity=1 ] (193.73,170) .. controls (193.73,167.94) and (192.06,166.27) .. (190,166.27) .. controls (187.94,166.27) and (186.27,167.94) .. (186.27,170) .. controls (186.27,172.06) and (187.94,173.73) .. (190,173.73) .. controls (192.06,173.73) and (193.73,172.06) .. (193.73,170) -- cycle ;

\draw (50,40.4) node [anchor=north west][inner sep=0.75pt]    {$\in \ \textcolor[rgb]{0.82,0.01,0.11}{\mathbb{M}^{s}}$};
\draw (190,120.4) node [anchor=north west][inner sep=0.75pt]    {$\in \ \mathbb{\textcolor[rgb]{0.82,0.01,0.11}{M}}\textcolor[rgb]{0.82,0.01,0.11}{^{s}}$};
\draw (139,20.4) node [anchor=north west][inner sep=0.75pt]    {$\in \ \textcolor[rgb]{0.56,0.07,1}{\mathbb{M}_{h}^{s}}$};
\draw (192,83.4) node [anchor=north west][inner sep=0.75pt]    {$v_{-1}$};
\draw (199,152.4) node [anchor=north west][inner sep=0.75pt]    {$v_{-2}$};
\draw (12,83.4) node [anchor=north west][inner sep=0.75pt]    {$v_{1}$};
\draw (81,19.4) node [anchor=north west][inner sep=0.75pt]    {$o=v_{0}$};

\end{tikzpicture}
\caption{Illustration of the configuration in the proof of Proposition \ref{prop:constructionMZ_h}.}\label{fig:constructionMZ_h}
\end{figure}
We want to show that
\[ \mathbb{E}[ w^s_h(v_0,v_{1})-w^{s}_h(v_0,v_{-1})] \leq 0. \]
To this end, we will use the recursion on $Z_h^s$ to bound the expectation by an expression involving $Z_h^{s}$. 

The expectation is zero outside the event $A=\{n_h(o)\neq n(o)\}=\{v_{-1} \neq v_1\}$. Let us work on this event from now on. 
By definition of $Z_h^{s}$, we have that \[Z_h^{s}(v_{-2},{v_{-1}}) = \max_{y \sim v_{-1}, y \neq v_{-2}} (w^{s}_h(v_{-1},y)-Z_h^{s}(v_{-1},y)). \] 
However, because $(v_0,v_{-1})$ is in $\mathbb{M}_h$, we know that \[v_{0}=\argmax_{y \sim v_{-1}} (w^{s}_h(v_{-1},y)-Z_h^{s}(v_{-1},y))\] and that $v_{-2}$ does not satisfy that maximum (else, $v_0=v_{-2}$ which contradicts being in $A$).

Hence
\begin{equation}\label{eq:aldous1}
Z_h^{s}(v_{-2},{v_{-1}}) = w^{s}_h(v_{-1},v_{0})-Z_h^{s}(v_{-1},v_{0}).  
\end{equation}
Now, by definition of the $Z_h^s$, we also have 
\begin{equation}\label{eq:aldous2}
    Z_h^{s}(v_{-1},v_0) = \max_{\substack{y \sim v_0 \\ y \neq v_{-1}}} ((w^{s}_h(v_{0},y)-Z_h^{s}(v_{0},y)) \geq w^{s}_h(v_{0},v_1)-Z_h^{s}(v_{0},v_1)
\end{equation}    
as we condition on $v_{-1} \neq v_{1}$. 

Combining \eqref{eq:aldous1} and \eqref{eq:aldous2}, we get that:
\begin{align}\label{eq:aldous4}
w^{s}_h(v_0,v_1)-w^{s}_h(v_0,v_{-1})&= w^{s}_h(v_0,v_1)-Z_h^{s}(v_{-1},v_0)-Z_h^{s}(v_{-2},v_{-1}) \\
&= \left[w^{s}_h(v_0,v_1)-Z_h^{s}(v_0,v_1)-Z_h^{s}(v_{-1},v_0) \right] +Z_h^{s}(v_0,v_1)-Z_h^{s}(v_{-2},v_{-1}) \\
&\leq Z_h^{s}(v_0,v_1)-Z_h^{s}(v_{-2},v_{-1}). 
\end{align}
Taking expectation and remembering that $v_1=n(o)$ and $v_{-1}=n_h(o)$, we get that:
\begin{align*}\label{eq:aldous3}
\perfV^s(\mathbb{T}^s_h,\mathbb{M}^s) - \perfV^s(\mathbb{T}^s_h, \mathbb{M}_h^s) 
& = \mathbb{E}[w^s_h(o, n(o))] - \mathbb{E}[w^s_h(o,n_h(o))] \\
& = \mathbb{E}\left[( w^{s}_h(v_0,v_1)-w^{s}_h(v_0,v_{-1}))\mathbbm{1}_{A}     \right] \\
& \leq \mathbb{E}\left[  ( Z_h^{s}(v_0,v_1)-Z_h^{s}(v_{-2},v_{-1}  ))\mathbbm{1}_{A}  \right].
\end{align*}
It remains to show that the final expectation is zero, to this end, we will use unimodularity to prove that conditionally on $A$, $Z_h^{s}(v_0,v_1)$ and $Z_h^{s}(v_{-2},v_{-1})$ have the same law.

Let $B$ some Borel set of $\mathbb{R}$, we want to show that:
\begin{equation} \label{eq:ZhonB}
\mathbb{P}\left(  Z_h^{s}(v_0,v_1)\in B, A \right)=\mathbb{P}\left(  Z_h^{s}(v_{-2},v_{-1})\in B, A \right).
\end{equation}
To prove this, we will use a chain of intermediary equalities between the two events by applying the mass-transport principle, which is possible thanks to Lemma~\ref{lem:Tsunimodular}. 

\begin{figure}[t!]
    \centering

\tikzset{every picture/.style={line width=0.75pt}} 

\begin{tikzpicture}[x=0.75pt,y=0.75pt,yscale=-1,xscale=1]

\draw [color={rgb, 255:red, 208; green, 2; blue, 27 }  ,draw opacity=1 ]   (100,32) -- (10,112) ;
\draw [color={rgb, 255:red, 144; green, 19; blue, 254 }  ,draw opacity=1 ]   (100,32) -- (190,112) ;
\draw [color={rgb, 255:red, 208; green, 2; blue, 27 }  ,draw opacity=1 ]   (190,112) -- (190,202) ;
\draw  [fill={rgb, 255:red, 0; green, 0; blue, 0 }  ,fill opacity=1 ] (103.73,32) .. controls (103.73,29.94) and (102.06,28.27) .. (100,28.27) .. controls (97.94,28.27) and (96.27,29.94) .. (96.27,32) .. controls (96.27,34.06) and (97.94,35.73) .. (100,35.73) .. controls (102.06,35.73) and (103.73,34.06) .. (103.73,32) -- cycle ;
\draw  [fill={rgb, 255:red, 0; green, 0; blue, 0 }  ,fill opacity=1 ] (13.73,112) .. controls (13.73,109.94) and (12.06,108.27) .. (10,108.27) .. controls (7.94,108.27) and (6.27,109.94) .. (6.27,112) .. controls (6.27,114.06) and (7.94,115.73) .. (10,115.73) .. controls (12.06,115.73) and (13.73,114.06) .. (13.73,112) -- cycle ;
\draw    (100,20) -- (11.49,98.67) ;
\draw [shift={(10,100)}, rotate = 318.37] [color={rgb, 255:red, 0; green, 0; blue, 0 }  ][line width=0.75]    (10.93,-3.29) .. controls (6.95,-1.4) and (3.31,-0.3) .. (0,0) .. controls (3.31,0.3) and (6.95,1.4) .. (10.93,3.29)   ;
\draw  [fill={rgb, 255:red, 0; green, 0; blue, 0 }  ,fill opacity=1 ] (193.73,112) .. controls (193.73,109.94) and (192.06,108.27) .. (190,108.27) .. controls (187.94,108.27) and (186.27,109.94) .. (186.27,112) .. controls (186.27,114.06) and (187.94,115.73) .. (190,115.73) .. controls (192.06,115.73) and (193.73,114.06) .. (193.73,112) -- cycle ;
\draw  [fill={rgb, 255:red, 0; green, 0; blue, 0 }  ,fill opacity=1 ] (193.73,202) .. controls (193.73,199.94) and (192.06,198.27) .. (190,198.27) .. controls (187.94,198.27) and (186.27,199.94) .. (186.27,202) .. controls (186.27,204.06) and (187.94,205.73) .. (190,205.73) .. controls (192.06,205.73) and (193.73,204.06) .. (193.73,202) -- cycle ;
\draw [color={rgb, 255:red, 208; green, 2; blue, 27 }  ,draw opacity=1 ]   (480,30) -- (390,110) ;
\draw [color={rgb, 255:red, 144; green, 19; blue, 254 }  ,draw opacity=1 ]   (480,30) -- (570,110) ;
\draw [color={rgb, 255:red, 208; green, 2; blue, 27 }  ,draw opacity=1 ]   (570,110) -- (570,200) ;
\draw    (480,20) -- (391.49,98.67) ;
\draw [shift={(390,100)}, rotate = 318.37] [color={rgb, 255:red, 0; green, 0; blue, 0 }  ][line width=0.75]    (10.93,-3.29) .. controls (6.95,-1.4) and (3.31,-0.3) .. (0,0) .. controls (3.31,0.3) and (6.95,1.4) .. (10.93,3.29)   ;
\draw  [fill={rgb, 255:red, 0; green, 0; blue, 0 }  ,fill opacity=1 ] (393.73,110) .. controls (393.73,107.94) and (392.06,106.27) .. (390,106.27) .. controls (387.94,106.27) and (386.27,107.94) .. (386.27,110) .. controls (386.27,112.06) and (387.94,113.73) .. (390,113.73) .. controls (392.06,113.73) and (393.73,112.06) .. (393.73,110) -- cycle ;
\draw  [fill={rgb, 255:red, 0; green, 0; blue, 0 }  ,fill opacity=1 ] (483.73,30) .. controls (483.73,27.94) and (482.06,26.27) .. (480,26.27) .. controls (477.94,26.27) and (476.27,27.94) .. (476.27,30) .. controls (476.27,32.06) and (477.94,33.73) .. (480,33.73) .. controls (482.06,33.73) and (483.73,32.06) .. (483.73,30) -- cycle ;
\draw  [fill={rgb, 255:red, 0; green, 0; blue, 0 }  ,fill opacity=1 ] (573.73,110) .. controls (573.73,107.94) and (572.06,106.27) .. (570,106.27) .. controls (567.94,106.27) and (566.27,107.94) .. (566.27,110) .. controls (566.27,112.06) and (567.94,113.73) .. (570,113.73) .. controls (572.06,113.73) and (573.73,112.06) .. (573.73,110) -- cycle ;
\draw  [fill={rgb, 255:red, 0; green, 0; blue, 0 }  ,fill opacity=1 ] (573.73,203.73) .. controls (573.73,201.67) and (572.06,200) .. (570,200) .. controls (567.94,200) and (566.27,201.67) .. (566.27,203.73) .. controls (566.27,205.78) and (567.94,207.45) .. (570,207.45) .. controls (572.06,207.45) and (573.73,205.78) .. (573.73,203.73) -- cycle ;
\draw   (270,120) -- (287.5,100) -- (287.5,110) -- (322.5,110) -- (322.5,100) -- (340,120) -- (322.5,140) -- (322.5,130) -- (287.5,130) -- (287.5,140) -- cycle ;

\draw (192,83.4) node [anchor=north west][inner sep=0.75pt]    {$v_{-1}$};
\draw (199,152.4) node [anchor=north west][inner sep=0.75pt]    {$v_{-2}$};
\draw (11,122.4) node [anchor=north west][inner sep=0.75pt]    {$v_{1}$};
\draw (111,10.4) node [anchor=north west][inner sep=0.75pt]    {$v_{0}$};
\draw (-19,32.4) node [anchor=north west][inner sep=0.75pt]    {$Z_h( u,v) \in B$};
\draw (122,137) node [anchor=north west][inner sep=0.75pt]   [align=left] {};
\draw (571,82.4) node [anchor=north west][inner sep=0.75pt]    {$v_{0}$};
\draw (581,152.4) node [anchor=north west][inner sep=0.75pt]    {$v_{1}$};
\draw (381,120.4) node [anchor=north west][inner sep=0.75pt]    {$v_{-2}$};
\draw (491,2.4) node [anchor=north west][inner sep=0.75pt]    {$v_{-1}$};
\draw (351,32.4) node [anchor=north west][inner sep=0.75pt]    {$Z_h( u,v) \in B$};
\draw (507,171.27) node [anchor=north west][inner sep=0.75pt]   [align=left] {};

\end{tikzpicture}

    \caption{The second idea is that the two events have the same probability seen either from $v_0$ or $v_{-1}$.}
    \label{fig:optimalityproof}
\end{figure}
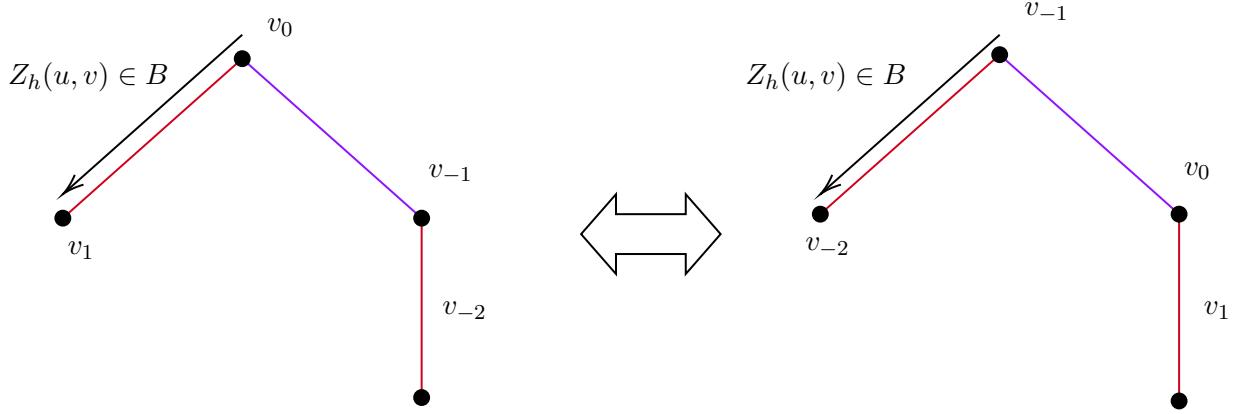

Let us recall the definition of $v_k$ so the desired equality rewrites as 
\[ \mathbb{P}\left(  Z_h^{s}(o,n(o))\in B, n(o)\neq n_h(o) \right)=\mathbb{P}\left(  Z_h^{s}(n(n_h(o)),n_h(o))\in B, n(o) \neq n_h(o) \right).\]
Define the measurable function $f$ on the space of double rooted  decorated trees
\[f(T,M,Z_h^s,a,b)=\mathbbm{1}_{ Z_h^s(a,b)\in B, n(a)\neq n_h(a), b=n(a)   } .  \]
Applying mass-transport principle to $f$ we get
\[ \mathbb{E} \left[  \sum_{v\in V} f(\mathbb{T}^s_h,\mathbb{M}^s,Z_h^s,o,v) \right] = \mathbb{E} \left[  \sum_{v\in V} f(\mathbb{T}^s_h,\mathbb{M}^s,Z_h^s,v,o) \right]  \]
Computing the first expectation yields
\begin{align*}
\mathbb{E}\left[ \sum_{v\in V}\mathbbm{1}_{Z_h^s(o,v) \in B, n(o) \neq n_h(o), v=n(o)}  \right] &= \mathbb{E}\left[ \mathbbm{1}_{Z_h^s(o,n(o))\in B, n(o) \neq n_h(o)} \right] \\
 &= \mathbb{P}\left( Z_h^s(o,n(o)) \in B , n(o) \neq n_h(o)  \right) .
 \end{align*}
 Computing the second expectation yields
\begin{align*}
\mathbb{E}\left[ \sum_{v\in V}\mathbbm{1}_{Z_h^s(v,o) \in B, n(v) \neq n_h(v), o=n(v)}  \right] &= \mathbb{E}\left[ \mathbbm{1}_{Z_h^s(n(o),o)\in B, n(n(o)) \neq n_h(n(o))} \right] \\
 &= \mathbb{P}\left( Z_h^s(n(o),o) \in B , n(o) \neq n_h(o)  \right) 
 \end{align*}
where we used that $o=n(v)$ is equivalent to $v=n(o)$ and the fact that $n(n(o))\neq n_h(n(o))$ is equivalent to $n(o) \neq n_h(o)$.

We now define another measurable function $f'$ on the same space:
\[f'(T,M,Z_h^s,a,b)= \mathbbm{1}_{Z_h^s(n(a),a)\in B, n(a) \neq n_h(a), b=n_h(a)} .\]
Once again, applying the mass-transport principle to $f'$ we get
\[ \mathbb{E} \left[  \sum_{v\in V} f'(\mathbb{T}^s_h,\mathbb{M}^s,Z_h^s,o,v) \right] = \mathbb{E} \left[  \sum_{v\in V} f'(\mathbb{T}^s_h,\mathbb{M}^s,Z_h^s,v,o) \right]  .\]
Computing the first expectation yields
\begin{align*}
\mathbb{E}\left[ \sum_{v\in V}\mathbbm{1}_{Z_h^s(n(o),o) \in B, n(o) \neq n_h(o), v=n_h(o)}  \right] &= \mathbb{E}\left[ \mathbbm{1}_{Z_h^s(n(o),o)\in B, n(o) \neq n_h(o)} \right] \\
 &= \mathbb{P}\left( Z_h^s(n(o),o) \in B , n(o) \neq n_h(o)  \right) .
 \end{align*}
 Computing the second expectation yields
 \begin{align*}
\mathbb{E}\left[ \sum_{v\in V}\mathbbm{1}_{Z_h^s(n(v),v) \in B, n(v) \neq n_h(v), o=n_h(v)}  \right] &= \mathbb{E}\left[ \mathbbm{1}_{Z_h^s(n(n_h(o)),n_h(o))\in B, n(n_h(o)) \neq n_h(n_h(o))} \right] \\
 &= \mathbb{P}\left( Z_h^s(n(n_h(o)),n_h(o)) \in B , n(o) \neq n_h(o)  \right) .
 \end{align*}
where we used that $o=n_h(v)$ is equivalent to $v=n_h(o)$ and the fact that $n(n_h(o)) \neq (n_h(n_h(o)))$ is equivalent to $n(o) \neq n_h(o)$. 

In conclusion, we showed
\[ \mathbb{P}\left(  Z_h^s(o,n(o))\in B, n(o)\neq n_h(o) \right)=\mathbb{P}\left( Z_h^s(n(o),o) \in B , n(o) \neq n_h(o)  \right)\]
\[=\mathbb{P}\left(  Z_h^s(n(n_h(o)),n_h(o))\in B, n(o) \neq n_h(o) \right),\]
yielding \eqref{eq:ZhonB} and the Lemma.
\end{proof}

We are now ready to prove our optimality result for unimodular matchings.

\begin{proof}[Proof of Proposition \ref{prop:optimal}]
Let $(\mathbb{T},\mathbb{M})$ an optimal matching on $\mathbb{T}$ and $n$ its associated neighbour function.
Let $(\mathbb{T},\mathbb{M},Z_h)$ such that the law of the marginals $(\mathbb{T},\mathbb{M})$ and $(\mathbb{T},Z_h)$ are preserved.
The main idea is that by adding self-loops with the rules given in the previous subsection, the matching $\mathbb M_h^s$ is the most penalised perfect matching of $\mathbb T^s_h$. This allows us to compare $\mathbb M$ and $\mathbb M_h$ by comparing their counterparts on $\mathbb T^s_h$.
The fact that $\mathbb M$ is optimal translates into the fact that
\[\perfV(\mathbb{M})=\mathbb{E}\left[ w^{s}_h(o,n(o))\mathbbm{1}_{o \neq n(o)}  \right] \]
is maximal among unimodular perfect matchings of $\mathbb T^s_h$.

Let us decompose $w^{s}_h(o,n(o))$ and $w^{s}_h(o,n_h(o))$ into their positive and negative parts. We denote by $\mathbb M^s$ the perfect matching of $\mathbb T_h^s$ corresponding to $\mathbb M$.
\begin{align*}
    \perfV^s(\mathbb{M}^{s}) &= \mathbb{E}\left[ w^{s}_h(o,n(o))\mathbbm{1}_{w^{s}_h(o,n(o))\geq 0}\right]+\mathbb{E}\left[ w^{s}_h(o,n(o))\mathbbm{1}_{w^{s}_h(o,n(o))< 0}\right] \\
    \perfV^s(\mathbb{M}_h^{s}) &= \mathbb{E}\left[ w^{s}_h(o,n_h(o))\mathbbm{1}_{w^{s}_h(o,n_h(o))\geq 0}\right]+\mathbb{E}\left[ w^{s}_h(o,n_h(o))\mathbbm{1}_{w^{s}_h(o,n_h(o))< 0}\right].
\end{align*}
But remembering \eqref{eq:critereselfloops}, $\mathbb{M}_h^{s}$ is precisely the matching that includes all the strictly negative weight self-loops and only them.

This implies that the previous decomposition can be written as
\[ \perfV^s(\mathbb{M}_h^s) = \perfV(\mathbb{M}_h) + \mathbb{E}[w^s_h(o,o)\mathbbm{1}_{w^{s}_h(o,o)<0} ]. \]
Now we can lower bound the positive part of $\perfV^s(\mathbb{M}^s)$ by $\perfV(\mathbb{M})$ by omitting the positive weight self-loops and we can lower bound the negative part by $w^s_h(o,o)\mathbbm{1}_{w^{s}_h(o,o)<0} .$ because $\mathbb{M}^{s}$ can at most include every negative self-loop: 
\[ \perfV^s(\mathbb{M}^{s}) \geq \perfV(\mathbb{M}) + \mathbb{E}[w^{s}_h(o,o))\mathbbm{1}_{w^{s}_h(o,o)<0}]. \]
But $\mathbb{M}$ is optimal so $\perfV(\mathbb{M}) \geq \perfV(\mathbb{M}_h)$, so we proved
\[ \perfV^s(\mathbb{M}^{s}) \geq \perfV(\mathbb{M}_h) + \mathbb{E}[w^{s}_h(o,o)\mathbbm{1}_{w^{s}_h(o,o)<0})] = \perfV^s(\mathbb{M}_h^s) .\]
However, we proved in Lemma \ref{lem:perfS} that $\mathbb{M}_h^{s}$ is optimal for $\perfV^{s}$ so
\[ \perfV^s(\mathbb{M}_h^s) \geq  \perfV^s(\mathbb{M}^{s})  .\]
So equality holds and all previous inequalities were equalities:
\begin{align*}
&\perfV(\mathbb{\mathbb{M}})=\perfV(\mathbb{M}_h),  \\
&\perfV(\mathbb{\mathbb{M}})=\mathbb{E}\left[ w^{s}_h(o,n(o))\mathbbm{1}_{w^{s}_h(o,n(o))\geq 0}\right], \\
&\mathbb{E}\left[ w^{s}_h(o,n(o))\mathbbm{1}_{w^{s}_h(o,n(o))< 0}\right]=\mathbb{E}[w^{s}_h(o,o)\mathbbm{1}_{w^{s}_h(o,o)<0})].
\end{align*}
So we have shown that $\mathbb{M}^{s}$ maximizes $\perfV^s$ and that $\mathbb{M}_h$ is optimal on $\mathbb{T}$. 
Furthermore, by unimodularity, we get that the set of unmatched vertices by $\mathbb{M}$ is almost surely the same as the set of unmatched vertices by $\mathbb{M}_h$, and it is almost surely the set of vertices $v$ such that $w^{s}_h(v,v)<0$.
\end{proof}

\begin{remark} \label{rmk:perf=perfs}
    If $\mathbb M$ is optimal, we have $ \perfV^s(\mathbb{M}^s) =  \perfV^s(\mathbb{M}_h^s)$ in addition to $ \perfV(\mathbb{M}) =  \perfV(\mathbb{M}_h)$.
\end{remark}

\subsection{Uniqueness}\label{sec:uniqueness}
We will keep using the vertex-rooted point of view in this section. The goal of this section is to prove that $\mathbb{M}_h$ is also the unique optimal matching in law. We will actually prove the following stronger result:

\begin{prop}\label{prop:optimalunique}
    Let $(\mathbb{T},\mathbb{M})$ be a unimodular optimal matched UBGW tree, sample $(\mathbb{T},\mathbb{M},Z_h)$ with Proposition~\ref{prop:constructionMZ_h}.
    Almost surely, $\mathbb{M}=\mathbb{M}_h$.

\end{prop}
\begin{proof}
We will show that almost surely, $\mathbb{M}^s=\mathbb{M}_h^s$ and deduce that almost surely, $\mathbb{M}=\mathbb{M}_h$ through the projection that forgets self-loops. Recall that $A$ is the event when $n_h(o)\neq n(o)$ and assume $\mathbb P(A) \neq 0$.





Going back into the proof of optimality of $\mathbb{M}_h^s$.
Taking expectation in Equation~\eqref{eq:aldous4} we get
\begin{align*}
\perfV^s(\mathbb{M}^s)-\perfV^s(\mathbb{M}_h^s)&= \mathbb{E}\left[(w^{s}_h(v_0,v_1)-Z_h^{s}(v_0,v_1)-Z_h^{s}(v_{-1},v_0))\mathbbm{1}_{v_1 \neq v_{-1}} \right] 
\end{align*}
Equality between $\perfV^s(\mathbb{M}_h^s)$ and $\perfV^s(\mathbb{M}^s)$ (see Remark~\ref{rmk:perf=perfs}) implies that the right-hand side is equal to $0$. Furthermore, Equation~\eqref{eq:aldous2} says that the variable in the expectation of the right-hand side is non negative almost surely. We conclude that, almost surely,
\[Z_h^{s}(v_{-1},o) =\max_{\substack {y \sim o \\ y \neq v_{-1}}}(w^{s}_h(o,y)-Z_h^{s}(o,y))    =w^{s}_h((o,v_1)-Z_h^{s}(o,v_1)). \]
By definition of $n_h$,
\[ v_{-1}=n_h(o) = \argmax_{y \sim o }(w^{s}_h(o,y)-Z_h^{s}(o,y)). \]
So $v_1$ achieves the maximum among the list of $(w^{s}_h(o,y)-Z_h^{s}(o,y))$ stripped of its maximum. 
We deduce that $v_1=n(o)$ achieves the second largest among the $(w^{s}_h(o,y)-Z_h^{s}(o,y))$ that we will write as $\overset{[2]}{\underset{y \sim o}{\argmax}}( w^{s}_h(o,y)-Z_h^{s}(o,y))$. 

Thus:
\[\mathbb{P}\left( n(o)=\underset{y \sim o}{\argmax}(w^{s}_h(o,y)-Z_h^{s}(o,y)) \text{ or } \underset{y \sim o}{\argmaxx}(w^{s}_h(o,y)-Z_h^{s}(o,y))  \right)=1. \]
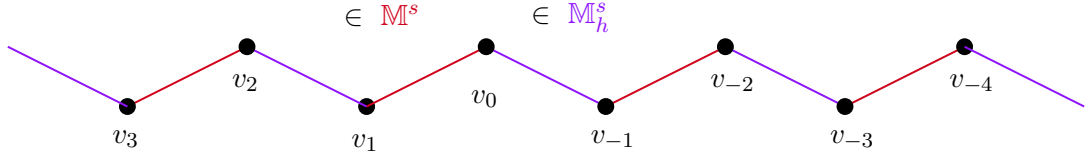
\begin{figure}[t!]
    \centering

\tikzset{every picture/.style={line width=0.75pt}} 

\begin{tikzpicture}[x=0.75pt,y=0.75pt,yscale=-1,xscale=1]

\draw  [fill={rgb, 255:red, 0; green, 0; blue, 0 }  ,fill opacity=1 ] (183.73,100) .. controls (183.73,97.94) and (182.06,96.27) .. (180,96.27) .. controls (177.94,96.27) and (176.27,97.94) .. (176.27,100) .. controls (176.27,102.06) and (177.94,103.73) .. (180,103.73) .. controls (182.06,103.73) and (183.73,102.06) .. (183.73,100) -- cycle ;
\draw [color={rgb, 255:red, 208; green, 2; blue, 27 }  ,draw opacity=1 ]   (300,100) -- (360,70) ;
\draw [color={rgb, 255:red, 144; green, 19; blue, 254 }  ,draw opacity=1 ]   (300,100) -- (240,70) ;
\draw [color={rgb, 255:red, 208; green, 2; blue, 27 }  ,draw opacity=1 ]   (180,100) -- (240,70) ;
\draw [color={rgb, 255:red, 208; green, 2; blue, 27 }  ,draw opacity=1 ]   (60,100) -- (120,70) ;
\draw [color={rgb, 255:red, 208; green, 2; blue, 27 }  ,draw opacity=1 ]   (420,100) -- (480,70) ;
\draw [color={rgb, 255:red, 144; green, 19; blue, 254 }  ,draw opacity=1 ]   (180,100) -- (120,70) ;
\draw [color={rgb, 255:red, 144; green, 19; blue, 254 }  ,draw opacity=1 ]   (420,100) -- (360,70) ;
\draw  [fill={rgb, 255:red, 0; green, 0; blue, 0 }  ,fill opacity=1 ] (483.73,70) .. controls (483.73,67.94) and (482.06,66.27) .. (480,66.27) .. controls (477.94,66.27) and (476.27,67.94) .. (476.27,70) .. controls (476.27,72.06) and (477.94,73.73) .. (480,73.73) .. controls (482.06,73.73) and (483.73,72.06) .. (483.73,70) -- cycle ;
\draw  [fill={rgb, 255:red, 0; green, 0; blue, 0 }  ,fill opacity=1 ] (423.73,100) .. controls (423.73,97.94) and (422.06,96.27) .. (420,96.27) .. controls (417.94,96.27) and (416.27,97.94) .. (416.27,100) .. controls (416.27,102.06) and (417.94,103.73) .. (420,103.73) .. controls (422.06,103.73) and (423.73,102.06) .. (423.73,100) -- cycle ;
\draw  [fill={rgb, 255:red, 0; green, 0; blue, 0 }  ,fill opacity=1 ] (363.73,70) .. controls (363.73,67.94) and (362.06,66.27) .. (360,66.27) .. controls (357.94,66.27) and (356.27,67.94) .. (356.27,70) .. controls (356.27,72.06) and (357.94,73.73) .. (360,73.73) .. controls (362.06,73.73) and (363.73,72.06) .. (363.73,70) -- cycle ;
\draw  [fill={rgb, 255:red, 0; green, 0; blue, 0 }  ,fill opacity=1 ] (303.73,100) .. controls (303.73,97.94) and (302.06,96.27) .. (300,96.27) .. controls (297.94,96.27) and (296.27,97.94) .. (296.27,100) .. controls (296.27,102.06) and (297.94,103.73) .. (300,103.73) .. controls (302.06,103.73) and (303.73,102.06) .. (303.73,100) -- cycle ;
\draw  [fill={rgb, 255:red, 0; green, 0; blue, 0 }  ,fill opacity=1 ] (243.73,70) .. controls (243.73,67.94) and (242.06,66.27) .. (240,66.27) .. controls (237.94,66.27) and (236.27,67.94) .. (236.27,70) .. controls (236.27,72.06) and (237.94,73.73) .. (240,73.73) .. controls (242.06,73.73) and (243.73,72.06) .. (243.73,70) -- cycle ;
\draw  [fill={rgb, 255:red, 0; green, 0; blue, 0 }  ,fill opacity=1 ] (123.73,70) .. controls (123.73,67.94) and (122.06,66.27) .. (120,66.27) .. controls (117.94,66.27) and (116.27,67.94) .. (116.27,70) .. controls (116.27,72.06) and (117.94,73.73) .. (120,73.73) .. controls (122.06,73.73) and (123.73,72.06) .. (123.73,70) -- cycle ;
\draw  [fill={rgb, 255:red, 0; green, 0; blue, 0 }  ,fill opacity=1 ] (63.73,100) .. controls (63.73,97.94) and (62.06,96.27) .. (60,96.27) .. controls (57.94,96.27) and (56.27,97.94) .. (56.27,100) .. controls (56.27,102.06) and (57.94,103.73) .. (60,103.73) .. controls (62.06,103.73) and (63.73,102.06) .. (63.73,100) -- cycle ;
\draw [color={rgb, 255:red, 144; green, 19; blue, 254 }  ,draw opacity=1 ]   (540,100) -- (480,70) ;
\draw [color={rgb, 255:red, 144; green, 19; blue, 254 }  ,draw opacity=1 ]   (60,100) -- (0,70) ;

\draw (167,47.4) node [anchor=north west][inner sep=0.75pt]    {$\in \ \mathbb{\textcolor[rgb]{0.82,0.01,0.11}{M}}\textcolor[rgb]{0.82,0.01,0.11}{^{s}}$};
\draw (261,46.4) node [anchor=north west][inner sep=0.75pt]    {$\in \ \mathbb{\textcolor[rgb]{0.56,0.07,1}{M}}\textcolor[rgb]{0.56,0.07,1}{_{h}^{s}}$};
\draw (171,112.4) node [anchor=north west][inner sep=0.75pt]    {$v_{1}$};
\draw (231,90.4) node [anchor=north west][inner sep=0.75pt]    {$v_{0}$};
\draw (291,110.4) node [anchor=north west][inner sep=0.75pt]    {$v_{-1}$};
\draw (111,82.4) node [anchor=north west][inner sep=0.75pt]    {$v_{2}$};
\draw (51,110.4) node [anchor=north west][inner sep=0.75pt]    {$v_{3}$};
\draw (351,82.4) node [anchor=north west][inner sep=0.75pt]    {$v_{-2}$};
\draw (411,110.4) node [anchor=north west][inner sep=0.75pt]    {$v_{-3}$};
\draw (471,82.4) node [anchor=north west][inner sep=0.75pt]    {$v_{-4}$};

\end{tikzpicture}

    \caption{Alternating path from $o$.}
    \label{fig:uniquenesspath}
\end{figure}
Now let us set: 
\begin{align*}
v_0&=o ,\\
v_1&=\underset{y \sim o}{\argmaxx}(w^{s}_h(o,y)-Z_h^{s}(o,y)) ,\\
v_{-1}&=\underset{y \sim o}{\argmax}(w^{s}_h(o,y)-Z_h^{s}(o,y)) ,\\
v_{-2k}&=\underset{y \sim w_{-2k+1}}{\argmaxx}(w^{s}_h(v_{-2k+1},y)-Z_h^{s}(v_{-2k+1},y))  &&\forall k \in \mathbb{N}^{\star},\\
v_{-2k-1}&=\underset{y \sim w_{-2k}}{\argmax}(w^{s}_h(v_{-2k},y)-Z_h^{s}(v_{-2k},y))  &&\forall k \in \mathbb{N}^{\star} , \\
v_{2k}&=\underset{y \sim w_{2k-1}}{\argmax}(w^{s}_h(v_{2k-1},y)-Z_h^{s}(v_{2k-1},y))  &&\forall k \in \mathbb{N}^{\star} , \\
v_{2k+1}&=\underset{y \sim w_{2k}}{\argmaxx}(w^{s}_h(v_{2k},y)-Z_h^{s}(v_{2k},y))  &&\forall k \in \mathbb{N}^{\star}.
\end{align*}
We always have by construction that $(v_{-2k},v_{-2k-1})\in \mathbb{M}_h$ for any $k \in \mathbb{Z}$. 

Conditionally on the event $A=\{n_h(o)\neq n(o)\}$, we get that:
\[\mathbb{P}\left(\bigcup_{k \in \mathbb{Z}}\left\{ (v_{-2k-1},v_{-2k-2})\in \mathbb{M}\right\} \middle| A \right)=1 \]
Furthermore, by the second part of Proposition~\ref{prop:optimal}, almost surely, for any $v\in V$, $n_h(v)=v$ if and only if $n(v)=v$, from this we deduce that on $A$ we have $v_{-1}\neq v_0$ and, by induction, $  \forall k \in \mathbb{N}, v_{-k} \neq v_{-k-1}$ (no loops in the path) and also that $\forall k \in \mathbb{N}, v_{-k-2} \neq v_{-k}$ (the path cannot go back up). 

Thus, $\mathbb{P}\left( k \mapsto v_k \text{ is into } | A   \right)=1.$ Let us write $C=\{ k \in \mathbb{Z}\mapsto v_k \text{ is into}$\}.  Define for $(N,N')\in \overline{\mathbb{Z}}^2$, $N< N'$:
 \[C_{N,N'}=\{k \mapsto v_k \text{ is into for } N \leq k \leq N' \}.\]
We have that $\mathbb{P}(A)>0$ and $\mathbb{P}(C|A)=1$ so $\mathbb{P}(C)\geq \mathbb{P}(A)>0$. Furthermore, $C \subset C_{-\infty,0}$ so $\mathbb{P}(C_{-\infty,0})>0$. We are going to show that this is not possible by establishing a contradiction for the value on the probability
\begin{equation} \label{eq:contradiction}
    \mathbb{P}\left( \overline{C_{0,2}} |  C_{-\infty, 0} \right).
\end{equation}

\bigskip

First, we will use unimodularity to show that $\mathbb{P}\left( \overline{C_{0,2}} |  C_{-\infty, 0} \right)=0$.
 We will rely on the following Lemma whose proof is postponed to the end of this Section.
\begin{lemma}\label{lem:doublyinfinitepath}
For any $N \in \mathbb{N}$,
\[ \mathbb{P}\left(C_{-2N,2} \right)=\mathbb{P}\left(C_{-2N-2,0} \right)   . \]    
\end{lemma}

\bigskip

From this Lemma we deduce that:
\[ \mathbb{P}(C_{-\infty,2})=\mathbb{P}(C_{-\infty,0}). \]
As long as $\deg(o)\geq 1$, there is one regular edge with $o$ as an endpoint and the self-loop $\{o,o\}$ the $\argmax$ and $\argmaxx$ around $o$ are thus automatically different. We can then decompose $C_{-\infty,2}$ as $C_{-\infty,0}\cap C_{0,2}\cap \{\deg(o)\geq 1\}$
\[ \mathbb{P}\left( C_{0,2}\cap \{\deg(o)\geq 1\}| C_{-\infty,0} \right)=1 . \]
Hence $\mathbb{P}\left( \overline{C_{0,2}} \cup \{\deg(o)< 1\}| C_{-\infty,0} \right) = 0 $ and

\[ \mathbb{P}\left( \overline{C_{0,2}} |  C_{-\infty, 0} \right) =0 .\]

\bigskip

Second, by exhibiting explicit configurations of positive total mass where $\overline{C_{0,2}}$ happens, we will show that $\mathbb{P}\left( \overline{C_{0,2}} |  C_{-\infty, 0} \right)>0$.

Set $T_{v_{-1}}$ the subtree rooted at $v_{-1}$. 

Referring to Figure~\ref{fig:uniquenessdependency}, the idea is that if $v_{0} \neq v_{-1}$, the only dependence the $\{v_{-k}, k \geq 1 \}$ has from $T \setminus T_{v_{-1}}$ comes from $Z_h(v_{0},v_{-1})$ and $w(v_{-1},v_0)$, so we will condition on the values of $Z_h(v_0,v_{-1}),w(v_0,v_{-1})$. Let $\mathcal{P}$ be the law of $Z_h(v_0,v_{-1}),w(v_0,v_{-1})$ conditionally on $C_{-\infty,0}$, we have then shown that:
\begin{align*}
    &\mathbb{P}\left( \overline{C_{0,2}} |  C_{-\infty, 0} \right) \\
    \geq& \iiint \mathbb{P}\left( \overline{C_{0,2}} | Z_h(v_0,v_{-1}),w(v_0,v_{-1}), C_{-\infty,0}  \right) \mathrm{d}\mathcal{P}. \\
    =& \iiint \mathbb{P}\left( \overline{C_{0,2}} | Z_h(v_0,v_{-1}),w(v_0,v_{-1})\right) \mathrm{d}\mathcal{P}.
\end{align*}
It suffices to show that
\[\mathbb{P}\left( \overline{C_{0,2}} | Z_h(v_0,v_{-1}),w(v_0,v_{-1})\right) >0 \] $\mathcal{P}$-almost everywhere and this will imply a contradiction.
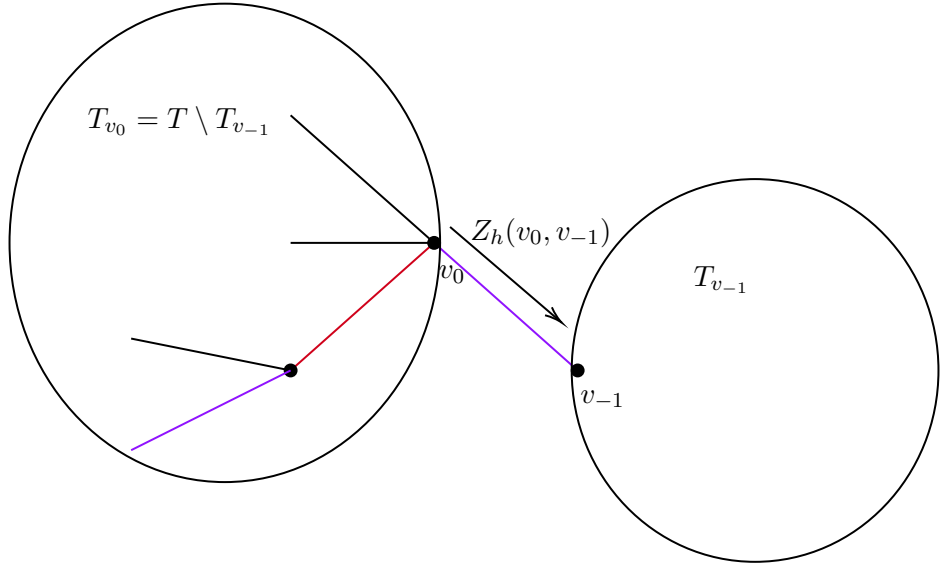
\begin{figure}[t!]
    \centering

\tikzset{every picture/.style={line width=0.75pt}} 

\begin{tikzpicture}[x=0.75pt,y=0.75pt,yscale=-0.8,xscale=0.8]

\draw [color={rgb, 255:red, 208; green, 2; blue, 27 }  ,draw opacity=1 ]   (270,100) -- (180,180) ;
\draw [color={rgb, 255:red, 144; green, 19; blue, 254 }  ,draw opacity=1 ]   (270,100) -- (360,180) ;
\draw  [fill={rgb, 255:red, 0; green, 0; blue, 0 }  ,fill opacity=1 ] (273.73,100) .. controls (273.73,97.94) and (272.06,96.27) .. (270,96.27) .. controls (267.94,96.27) and (266.27,97.94) .. (266.27,100) .. controls (266.27,102.06) and (267.94,103.73) .. (270,103.73) .. controls (272.06,103.73) and (273.73,102.06) .. (273.73,100) -- cycle ;
\draw  [fill={rgb, 255:red, 0; green, 0; blue, 0 }  ,fill opacity=1 ] (183.73,180) .. controls (183.73,177.94) and (182.06,176.27) .. (180,176.27) .. controls (177.94,176.27) and (176.27,177.94) .. (176.27,180) .. controls (176.27,182.06) and (177.94,183.73) .. (180,183.73) .. controls (182.06,183.73) and (183.73,182.06) .. (183.73,180) -- cycle ;
\draw  [fill={rgb, 255:red, 0; green, 0; blue, 0 }  ,fill opacity=1 ] (363.73,180) .. controls (363.73,177.94) and (362.06,176.27) .. (360,176.27) .. controls (357.94,176.27) and (356.27,177.94) .. (356.27,180) .. controls (356.27,182.06) and (357.94,183.73) .. (360,183.73) .. controls (362.06,183.73) and (363.73,182.06) .. (363.73,180) -- cycle ;
\draw   (356.27,180) .. controls (356.27,113.73) and (407.76,60) .. (471.27,60) .. controls (534.79,60) and (586.27,113.73) .. (586.27,180) .. controls (586.27,246.27) and (534.79,300) .. (471.27,300) .. controls (407.76,300) and (356.27,246.27) .. (356.27,180) -- cycle ;
\draw    (270,100) -- (180,100) ;
\draw    (270,100) -- (180,20) ;
\draw    (180,180) -- (80,160) ;
\draw [color={rgb, 255:red, 144; green, 19; blue, 254 }  ,draw opacity=1 ]   (180,180) -- (80,230) ;
\draw   (3.73,100) .. controls (3.73,17.16) and (64.17,-50) .. (138.73,-50) .. controls (213.28,-50) and (273.73,17.16) .. (273.73,100) .. controls (273.73,182.84) and (213.28,250) .. (138.73,250) .. controls (64.17,250) and (3.73,182.84) .. (3.73,100) -- cycle ;
\draw    (280,90) -- (348.48,148.7) ;
\draw [shift={(350,150)}, rotate = 220.6] [color={rgb, 255:red, 0; green, 0; blue, 0 }  ][line width=0.75]    (10.93,-3.29) .. controls (6.95,-1.4) and (3.31,-0.3) .. (0,0) .. controls (3.31,0.3) and (6.95,1.4) .. (10.93,3.29)   ;

\draw (360,190.4) node [anchor=north west][inner sep=0.75pt]    {$v_{-1}$};
\draw (271,112.4) node [anchor=north west][inner sep=0.75pt]    {$v_{0}$};
\draw (51,12.4) node [anchor=north west][inner sep=0.75pt]    {$T_{v_{0}} =T\setminus T_{v_{-1}}{}$};
\draw (431,112.4) node [anchor=north west][inner sep=0.75pt]    {$T_{v_{-1}}$};
\draw (291,82.4) node [anchor=north west][inner sep=0.75pt]    {$Z_h( v_{0} ,v_{-1})$};

\end{tikzpicture}

    \caption{The variables attached to $T_{v_0}$ only depend on $T_{v_{-1}} \cup \{ v_0, v_{-1}\}$ and its attached variables through $Z(v_0,v_{-1})$ and $w(v_0,v_{-1})$.}
    \label{fig:uniquenessdependency}
\end{figure}

The next Lemma shows an explicit lower bound as required:
\begin{lemma}\label{lem:Mhstop}
    Let $(W_i)_{i \geq 1}$ be i.i.d variables of law $\omega$.

    Let $\Tilde{\pi}_k=\mathbb{P}(\deg(o)=k| C_{-\infty,0})$.

    Let $\mathbf{X}$ any Borel subset of $\mathrm{supp}(\omega) \times\mathrm{supp}(\zeta_h)$:
    \[\mathbb{P}\left(\overline{C_{0,2}} |   \left(w(v_0,v_{-1}),Z_h(v_0,v_{-1})\right)\in \mathbf{X} \right) \geq \sum_{k>1}\Tilde{\pi}_k\mathbb{P}_{Z \sim \zeta_h}(Z=0)^{k-1} .\]
\end{lemma}
Lemma~\ref{lem:Zatomic} implies that $\mathbb{P}_{Z \sim \zeta_h}(Z=0)>0$.  This uniform lower bound does not depend on the values of $w(v_0,v_{-1}),Z_h(v_0,v_{-1})$,  so it is $\mathcal{P}$ almost-sure.
\end{proof}

We now turn to the proof of our two Lemmas.

\begin{proof}[Proof of Lemma \ref{lem:doublyinfinitepath}]
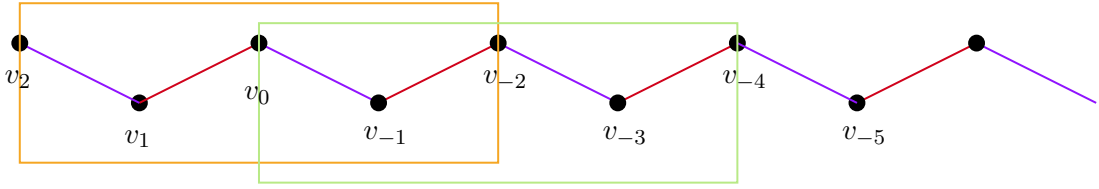
\begin{figure}[h!]
    \centering

\tikzset{every picture/.style={line width=0.75pt}} 

\begin{tikzpicture}[x=0.75pt,y=0.75pt,yscale=-1,xscale=1]

\draw  [fill={rgb, 255:red, 0; green, 0; blue, 0 }  ,fill opacity=1 ] (183.73,100) .. controls (183.73,97.94) and (182.06,96.27) .. (180,96.27) .. controls (177.94,96.27) and (176.27,97.94) .. (176.27,100) .. controls (176.27,102.06) and (177.94,103.73) .. (180,103.73) .. controls (182.06,103.73) and (183.73,102.06) .. (183.73,100) -- cycle ;
\draw [color={rgb, 255:red, 208; green, 2; blue, 27 }  ,draw opacity=1 ]   (300,100) -- (360,70) ;
\draw [color={rgb, 255:red, 144; green, 19; blue, 254 }  ,draw opacity=1 ]   (300,100) -- (240,70) ;
\draw [color={rgb, 255:red, 208; green, 2; blue, 27 }  ,draw opacity=1 ]   (180,100) -- (240,70) ;
\draw [color={rgb, 255:red, 208; green, 2; blue, 27 }  ,draw opacity=1 ]   (540,100) -- (600,70) ;
\draw [color={rgb, 255:red, 208; green, 2; blue, 27 }  ,draw opacity=1 ]   (420,100) -- (480,70) ;
\draw [color={rgb, 255:red, 144; green, 19; blue, 254 }  ,draw opacity=1 ]   (180,100) -- (120,70) ;
\draw [color={rgb, 255:red, 144; green, 19; blue, 254 }  ,draw opacity=1 ]   (420,100) -- (360,70) ;
\draw  [fill={rgb, 255:red, 0; green, 0; blue, 0 }  ,fill opacity=1 ] (483.73,70) .. controls (483.73,67.94) and (482.06,66.27) .. (480,66.27) .. controls (477.94,66.27) and (476.27,67.94) .. (476.27,70) .. controls (476.27,72.06) and (477.94,73.73) .. (480,73.73) .. controls (482.06,73.73) and (483.73,72.06) .. (483.73,70) -- cycle ;
\draw  [fill={rgb, 255:red, 0; green, 0; blue, 0 }  ,fill opacity=1 ] (423.73,100) .. controls (423.73,97.94) and (422.06,96.27) .. (420,96.27) .. controls (417.94,96.27) and (416.27,97.94) .. (416.27,100) .. controls (416.27,102.06) and (417.94,103.73) .. (420,103.73) .. controls (422.06,103.73) and (423.73,102.06) .. (423.73,100) -- cycle ;
\draw  [fill={rgb, 255:red, 0; green, 0; blue, 0 }  ,fill opacity=1 ] (363.73,70) .. controls (363.73,67.94) and (362.06,66.27) .. (360,66.27) .. controls (357.94,66.27) and (356.27,67.94) .. (356.27,70) .. controls (356.27,72.06) and (357.94,73.73) .. (360,73.73) .. controls (362.06,73.73) and (363.73,72.06) .. (363.73,70) -- cycle ;
\draw  [fill={rgb, 255:red, 0; green, 0; blue, 0 }  ,fill opacity=1 ] (303.73,100) .. controls (303.73,97.94) and (302.06,96.27) .. (300,96.27) .. controls (297.94,96.27) and (296.27,97.94) .. (296.27,100) .. controls (296.27,102.06) and (297.94,103.73) .. (300,103.73) .. controls (302.06,103.73) and (303.73,102.06) .. (303.73,100) -- cycle ;
\draw  [fill={rgb, 255:red, 0; green, 0; blue, 0 }  ,fill opacity=1 ] (243.73,70) .. controls (243.73,67.94) and (242.06,66.27) .. (240,66.27) .. controls (237.94,66.27) and (236.27,67.94) .. (236.27,70) .. controls (236.27,72.06) and (237.94,73.73) .. (240,73.73) .. controls (242.06,73.73) and (243.73,72.06) .. (243.73,70) -- cycle ;
\draw  [fill={rgb, 255:red, 0; green, 0; blue, 0 }  ,fill opacity=1 ] (123.73,70) .. controls (123.73,67.94) and (122.06,66.27) .. (120,66.27) .. controls (117.94,66.27) and (116.27,67.94) .. (116.27,70) .. controls (116.27,72.06) and (117.94,73.73) .. (120,73.73) .. controls (122.06,73.73) and (123.73,72.06) .. (123.73,70) -- cycle ;
\draw  [fill={rgb, 255:red, 0; green, 0; blue, 0 }  ,fill opacity=1 ] (543.73,100) .. controls (543.73,97.94) and (542.06,96.27) .. (540,96.27) .. controls (537.94,96.27) and (536.27,97.94) .. (536.27,100) .. controls (536.27,102.06) and (537.94,103.73) .. (540,103.73) .. controls (542.06,103.73) and (543.73,102.06) .. (543.73,100) -- cycle ;
\draw [color={rgb, 255:red, 144; green, 19; blue, 254 }  ,draw opacity=1 ]   (540,100) -- (480,70) ;
\draw [color={rgb, 255:red, 144; green, 19; blue, 254 }  ,draw opacity=1 ]   (660,100) -- (600,70) ;
\draw  [fill={rgb, 255:red, 0; green, 0; blue, 0 }  ,fill opacity=1 ] (603.73,70) .. controls (603.73,67.94) and (602.06,66.27) .. (600,66.27) .. controls (597.94,66.27) and (596.27,67.94) .. (596.27,70) .. controls (596.27,72.06) and (597.94,73.73) .. (600,73.73) .. controls (602.06,73.73) and (603.73,72.06) .. (603.73,70) -- cycle ;
\draw  [color={rgb, 255:red, 245; green, 166; blue, 35 }  ,draw opacity=1 ] (120,50) -- (360,50) -- (360,130) -- (120,130) -- cycle ;
\draw  [color={rgb, 255:red, 184; green, 233; blue, 134 }  ,draw opacity=1 ] (240,60) -- (480,60) -- (480,140) -- (240,140) -- cycle ;

\draw (171,112.4) node [anchor=north west][inner sep=0.75pt]    {$v_{1}$};
\draw (231,90.4) node [anchor=north west][inner sep=0.75pt]    {$v_{0}$};
\draw (291,110.4) node [anchor=north west][inner sep=0.75pt]    {$v_{-1}$};
\draw (111,82.4) node [anchor=north west][inner sep=0.75pt]    {$v_{2}$};
\draw (531,110.4) node [anchor=north west][inner sep=0.75pt]    {$v_{-5}$};
\draw (351,82.4) node [anchor=north west][inner sep=0.75pt]    {$v_{-2}$};
\draw (411,110.4) node [anchor=north west][inner sep=0.75pt]    {$v_{-3}$};
\draw (471,82.4) node [anchor=north west][inner sep=0.75pt]    {$v_{-4}$};

\end{tikzpicture}

    \caption{We see the same event from the perspectives of $v_2$ inside the orange box and $v_0$ inside the green box.}
    \label{fig:uniquenessreroot}
\end{figure}
    Refer to Figure~\ref{fig:uniquenessreroot} for an illustration.
    Define the alternating path starting from a vertex $u$ as:
    \[ u_1(u)=\argmax_{y \sim u}(w^{s}_h(u,y)-Z_h^{s}(u,y)), u_2(u)=\argmaxx_{y \sim u_1}(w^{s}_h(u_1,y)-Z_h^{s}(u_1,y)),  \]
    and for $n \geq 1,$
    \[ u_{2n+1}(u)=\argmax_{y \sim u_{2n}}(w^{s}_h(u_{2n},y)-Z_h^{s}(u_{2n},y)), u_{2n+2}(u)=\argmaxx_{y \sim u_{2n+1}}(w^{s}_h(u_{2n+1},y)-Z_h^{s}(u_{2n+1},y)).  \]
    By convention, set $u_0=\text{Id}_{V}$, in other words, $u_0(u)=u$.
    Let $N\in \mathbb{N}$, define events:
    \[ C_N(T,Z_h,a,b)=  \left\{ n \mapsto u_n(a) \text{ is into for } 0 \leq n \leq 2N    \right\} \cap \left\{ n \mapsto u_n(b) \text{ is into for } 0 \leq n \leq 2N    \right\} \]
    
    Let us define $f_N(T,M,Z_h,a,b):=\mathbbm{1}_{C_N(T,Z_h,a,b) \cap \{b=u_{2}(a)\}  }$.
    
    On one hand,
    \[ \mathbb{E}\left[\sum_{v \in V}f_N(T,M,Z_h,o,v)\right]=\mathbb{P}\left(C_{-N,2}\right) ,\]
    on the other hand
    \[ \mathbb{E}\left[\sum_{v \in V}f_N(T,M,Z_h,v,o)\right]=\mathbb{P}\left(C_{-N-2,0} \right)  . \]
    Applying mass-transport principle to $f_N$ yields the result.
\end{proof}

\begin{proof}[Proof of Lemma \ref{lem:Mhstop}]
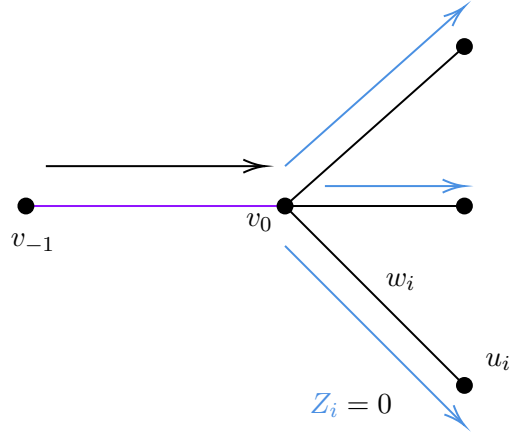
\begin{figure}[t!]

\centering

\tikzset{every picture/.style={line width=0.75pt}} 

\begin{tikzpicture}[x=0.75pt,y=0.75pt,yscale=-1,xscale=1]

\draw [color={rgb, 255:red, 144; green, 19; blue, 254 }  ,draw opacity=1 ]   (160,120) -- (30,120) ;
\draw    (160,120) -- (250,40) ;
\draw    (160,120) -- (250,120) ;
\draw    (160,120) -- (250,210) ;
\draw [color={rgb, 255:red, 74; green, 144; blue, 226 }  ,draw opacity=1 ]   (160,100) -- (248.51,21.33) ;
\draw [shift={(250,20)}, rotate = 138.37] [color={rgb, 255:red, 74; green, 144; blue, 226 }  ,draw opacity=1 ][line width=0.75]    (10.93,-3.29) .. controls (6.95,-1.4) and (3.31,-0.3) .. (0,0) .. controls (3.31,0.3) and (6.95,1.4) .. (10.93,3.29)   ;
\draw [color={rgb, 255:red, 74; green, 144; blue, 226 }  ,draw opacity=1 ]   (180,110) -- (248,110) ;
\draw [shift={(250,110)}, rotate = 180] [color={rgb, 255:red, 74; green, 144; blue, 226 }  ,draw opacity=1 ][line width=0.75]    (10.93,-3.29) .. controls (6.95,-1.4) and (3.31,-0.3) .. (0,0) .. controls (3.31,0.3) and (6.95,1.4) .. (10.93,3.29)   ;
\draw [color={rgb, 255:red, 74; green, 144; blue, 226 }  ,draw opacity=1 ]   (160,140) -- (248.59,228.59) ;
\draw [shift={(250,230)}, rotate = 225] [color={rgb, 255:red, 74; green, 144; blue, 226 }  ,draw opacity=1 ][line width=0.75]    (10.93,-3.29) .. controls (6.95,-1.4) and (3.31,-0.3) .. (0,0) .. controls (3.31,0.3) and (6.95,1.4) .. (10.93,3.29)   ;
\draw  [fill={rgb, 255:red, 0; green, 0; blue, 0 }  ,fill opacity=1 ] (33.73,120) .. controls (33.73,117.94) and (32.06,116.27) .. (30,116.27) .. controls (27.94,116.27) and (26.27,117.94) .. (26.27,120) .. controls (26.27,122.06) and (27.94,123.73) .. (30,123.73) .. controls (32.06,123.73) and (33.73,122.06) .. (33.73,120) -- cycle ;
\draw  [fill={rgb, 255:red, 0; green, 0; blue, 0 }  ,fill opacity=1 ] (163.73,120) .. controls (163.73,117.94) and (162.06,116.27) .. (160,116.27) .. controls (157.94,116.27) and (156.27,117.94) .. (156.27,120) .. controls (156.27,122.06) and (157.94,123.73) .. (160,123.73) .. controls (162.06,123.73) and (163.73,122.06) .. (163.73,120) -- cycle ;
\draw  [fill={rgb, 255:red, 0; green, 0; blue, 0 }  ,fill opacity=1 ] (253.73,210) .. controls (253.73,207.94) and (252.06,206.27) .. (250,206.27) .. controls (247.94,206.27) and (246.27,207.94) .. (246.27,210) .. controls (246.27,212.06) and (247.94,213.73) .. (250,213.73) .. controls (252.06,213.73) and (253.73,212.06) .. (253.73,210) -- cycle ;
\draw  [fill={rgb, 255:red, 0; green, 0; blue, 0 }  ,fill opacity=1 ] (253.73,120) .. controls (253.73,117.94) and (252.06,116.27) .. (250,116.27) .. controls (247.94,116.27) and (246.27,117.94) .. (246.27,120) .. controls (246.27,122.06) and (247.94,123.73) .. (250,123.73) .. controls (252.06,123.73) and (253.73,122.06) .. (253.73,120) -- cycle ;
\draw  [fill={rgb, 255:red, 0; green, 0; blue, 0 }  ,fill opacity=1 ] (253.73,40) .. controls (253.73,37.94) and (252.06,36.27) .. (250,36.27) .. controls (247.94,36.27) and (246.27,37.94) .. (246.27,40) .. controls (246.27,42.06) and (247.94,43.73) .. (250,43.73) .. controls (252.06,43.73) and (253.73,42.06) .. (253.73,40) -- cycle ;
\draw    (40,100) -- (148,100) ;
\draw [shift={(150,100)}, rotate = 180] [color={rgb, 255:red, 0; green, 0; blue, 0 }  ][line width=0.75]    (10.93,-3.29) .. controls (6.95,-1.4) and (3.31,-0.3) .. (0,0) .. controls (3.31,0.3) and (6.95,1.4) .. (10.93,3.29)   ;

\draw (139,122.4) node [anchor=north west][inner sep=0.75pt]    {$v_{0}$};
\draw (21,132.4) node [anchor=north west][inner sep=0.75pt]    {$v_{-1}$};
\draw (171,212.4) node [anchor=north west][inner sep=0.75pt]    {$\textcolor[rgb]{0.29,0.56,0.89}{Z_{i}} =0$};
\draw (259,192.4) node [anchor=north west][inner sep=0.75pt]    {$u_{i}$};
\draw (209,152.4) node [anchor=north west][inner sep=0.75pt]    {$w_{i}$};

\end{tikzpicture}

    \caption{For $\mathbb{M}_h$ to stop, it suffices that all the $Z_i=0$.}
    \label{fig:uniquenesslemma}
\end{figure}

Let $u_i$ be the children of $o$ as vertices for $i \leq \deg(o)$, $Z_i=Z_h(o,u_i)$.
This follows from seeing that:
\begin{align*}
&\left\{  \left(W(v_0,v_{-1}),Z_h(v_0,v_{-1})\right) \in \mathbf{X}  ,\overline{C_{0,2}} \right\}   \\ 
\subseteq&\left\{  \left(W(v_0,v_{-1}),Z_h(v_0,v_{-1})\right) \in \mathbf{X},n_h(v_1)=v_1\right\} \\
\subseteq&\left\{ \left(W(v_0,v_{-1}),Z_h(v_0,v_{-1})\right) \in \mathbf{X}, \forall i \leq \deg(o), Z_i=0\right\} .\\
\end{align*}
Refer to Figure~\ref{fig:uniquenesslemma} for an illustration. 
Every variable appearing in this event is independent of all others, so we can sum over $\deg(o)$ to obtain the bound:

\begin{align*}
    &\mathbb{P}\left(\overline{C_{0,2}} |  \left(w(v_0,v_{-1}),Z_h(v_0,v_{-1})\right)\in \mathbf{X} \right) \\
    &\geq \frac{\mathbb{P}\left(\left(w(v_0,v_{-1}),Z_h(v_0,v_{-1})\right) \in \mathbf{X} \right) \sum_{k>1}\Tilde{\pi}_k\mathbb{P}_{Z \sim \zeta_h}(Z=0)^{k-1}}{\mathbb{P} \left(w(v_0,v_{-1}),Z_h(v_0,v_{-1})\right)\in \mathbf{X}) }   .
\end{align*}
The factors $\mathbb{P}\left(\left(w(v_0,v_{-1}),Z_h(v_0,v_{-1})\right)\in \mathbf{X}\right)$ cancel out, so we get the desired expression.
\end{proof}

\section{From infinite to finite : Theorem~\ref{maintheorem}}
\label{sec:infinitetofinite}

In this entire section, unless stated otherwise, we will adopt the edge-rooted point of view. The aim of this section is to prove Theorem~\ref{maintheorem}, namely that the optimal matching on a random graph converges locally to the unimodular optimal matching on the limiting UBGW tree. Theorem~\ref{unicityUBGW}, proved in Section~\ref{sec:matchUBGW}, states that there is a unique unimodular optimal matching on a UBGW tree with iid weights. A natural approach to prove Theorem~\ref{maintheorem} is then to show that every subsequential limit of optimal matchings is a unimodular optimal matching on the limiting tree. More precisely, we will prove the following statement:
\begin{prop} \label{prop:infinifini}
    Assume $(G_n,o_n)_{n\geq1}$ is a sequence of finite graphs which converges locally to the unimodular BGW tree $(\mathbb{T},o)$ with reproduction law $\pi$ and weights $\omega$. 
    
    Assume $\mathbb{E}_{N \sim \pi}[N]< \infty$, $\mathbb{E}_{W \sim \omega}[W]< \infty$, $\omega$ is atomless, and that $\frac{|E_n|}{|V_n|}$ converges in probability to $\frac{\mathbb{E}_{N\sim \pi}[N]}{2}$.
    
    Let $M_{\mathrm{opt}}(G_n)$ be the optimal matching on $G_n$, then every subsequence of $(G_n,o_n,M_{\mathrm{opt}}(G_n))$ has a subsequence that converges locally to $(\mathbb{T},o,\mathbb{M})$ where $\mathbb{M}$ is an optimal unimodular matching on $\mathbb{T}$. 
\end{prop}
The proof can be broken down into two parts, one relatively easy, and one more challenging:
\begin{itemize}
    \item Every subsequence of $(G_n,o_n,M_{\mathrm{opt}}(G_n))$ has a subsequence that converges locally to some $(\mathbb{T},o,\mathbb{M})$ where $\mathbb{M}$ is a unimodular matching on $\mathbb{T}$.
    \item The previously found $\mathbb{M}$ must be optimal.
\end{itemize}

The complexity lies in  showing that the limit of any convergent subsequence is necessarily optimal. We will do this by constructing an almost optimal matching on $G_n$ from the optimal matching $(\mathbb{T}, \mathbb{M}_h)$. Indeed, as a consequence, this will show that $M_{\mathrm{opt}}(G_n)$ does asymptotically better than $\perfV(\mathbb{M}_h)$.
\subsection{Local convergence of subsequences}\label{sec:tightness}
In this section we will carry out the easy part using the following lemma:
\begin{lemma}\label{lem:tightness}
    Let $(G_n,o_n)$ be a sequence in the space of rooted graphs, $G_n=(V_n,E_n)$. Let $(f_n)_{n \in \mathbb{N}}$ be any (random) functions on $E_n$ such that the sequence $(f_n)$ is uniformly bounded. Assume that $G_n$ converges locally to an almost surely locally finite graph $G$. Then for any $H\in \mathbb{N}$,  the sequence of decorated graphs $B_H(G_n,o_n,f_n)$ is tight. 
\end{lemma}
\begin{proof}
    Let $H \in \mathbb{N}$ and $\varepsilon >0$. Denote by 
    $\tau_{H,N}$ the set of rooted graphs of depth $H$ and degree bounded by $N$. Since $G$ is locally finite, for $N$ large enough, we have $\mathbb{P}(B_H(G,o) \in \tau_{H,N}) > 1- \frac{\epsilon}{2}$.
    By local convergence, $\mathbb{P}(B_H(G_n,o_n) \in \tau_{H,N})  \underset{n \rightarrow \infty}{\longrightarrow} \mathbb{P}(B_H(G,o) \in \tau_{H,N}) $ hence for $n$ large enough,
    \[ \mathbb{P}(B_H(G_n,o_n) \in \tau_{H,N}) > 1-\epsilon \]
    Let $A=\sup_{n \geq 0} \|f_n\|_{\infty} $, denote by $\tau_{H,N,A}$ the set of edge-decorated graphs consisting of elements of $\tau_{H,N}$ decorated by a function bounded by $A$. Then $\tau_{H,N,A}$ is compact and for $n$ large enough:
    \[ \mathbb{P}(B_H(G_n,o_n,f_n)\in \tau_{H,N,A}) >1-\epsilon,\]
    $B_H(G_n,o,f_n)$ is thus tight.
\end{proof}
We can now go back to the proof of Proposition $\ref{prop:infinifini}$.
Take $f_n=\mathbbm{1}_{M_\mathrm{opt}(G_n)}$, it is uniformly bounded by $1$.
Apply Lemma~\ref{lem:tightness} to $B_0(G_n,o_n,f_n)$ to obtain a subsequence such that $B_0(G_{k_n},o_{k_n},f_{k_n})$ converges. 

But $(G_{k_n},o_{k_n})$ still converges locally to $(G,o)$, so apply Lemma~\ref{lem:tightness} again to find a subsequence $k_n'$ of $k_n$ such that $B_1(G_{k'_n},o_{k'_n},f_{k'_n})$ converges locally. Repeat the process and use a diagonal argument to obtain a subsequence $l_n$ such that for all $H \in \mathbb{N}$, $B_H(G_{l_n},o_{l_n},f_{l_n})$ converges weakly to some limit $B_H(G,o,g_{H})$ for some decoration $g_H$.

Furthermore, since the decorations $g_H$ are obtained from a diagonal extraction, the family $(g_H)_{H \geq 0}$ is consistent and therefore defines a limiting decoration $g$. To prove that this decoration is a matching, we just have to take  successively limits in $n \rightarrow \infty$ and $H \rightarrow \infty$ in the equalities 
\[ \mathbb{P}\left(\sum_{u \sim v} f_{l_n}(u,v) \leq 1, (u,v) \in V(B_H(G_n,o_n))\right)=1 ,\mathbb{P} \left(f_{l_n}(e) \in \{0,1\}, e \in E(B_H(G_n,o_n))\right)=1. \]

Now we just need to check that the limiting object is unimodular, it is a consequence that any local limit of unimodular graphs is unimodular, we refer to \cite{aldous2018processes}.

The sequence of graphs $(G_n,o_n)$ converges weakly to $\mathbb{T}$. Take $j_n$ a subsequence such that $\perfV(G_{j_n},\mathbb{M}_{\mathrm{opt}}(G_{j_n}))$ converges to $\limsup_{n \rightarrow \infty} \perfV(G_n,\mathbb{M}_{\mathrm{opt}}(G_n))$.

Now do the previous procedure on the sequence of $(G_{j_n},o_{j_n},f_{j_n})$ to find a subsequence $l_n$ such that $(G_{l_n},o_{l_n},f_{l_n})$ converges weakly to some $(\mathbb{T},g)$.
Set $\mathbb{M}$ the matching on $\mathbb{T}$ such that $\mathbbm{1}_\mathbb{M}=g$, then $\mathbb{M}$ is at most optimal by definition so it has less performance than $\mathbb{M}_h$:

\[ \limsup_{n \rightarrow \infty} \perfV(G_n,\mathbb{M}_{\mathrm{opt}}(G_n))= \perfV(\mathbb{T},\mathbb{M})  \leq  \perfV(\mathbb{T},\mathbb{M}_h).   \]
We have thus shown that every subsequence of $(G_n,o_n,M_{\mathrm{opt}}(G_n))$ has a locally convergent subsequence to some $(\mathbb{T},\mathbb{M})$ where $\mathbb{M}$ is an at most optimal matching on $\mathbb{T}$.

\subsection{Optimality of limits of subsequences}\label{sec:edgebenchmark}
In this section, we carry out the hard part of the proof of Proposition $\ref{prop:infinifini}$. 
We want to reconstruct a quasi matching on $G_n=(V_n,E_n)$ from $(\mathbb{T},\mathbb{M}_h)$ in the sense that the quasi matching is a subgraph that differs from a true matching with $o(|E_n|)$ amount of edges . The crucial idea is, for every $(i,j) \in V_n^2$, to compute the probability that $(i,j) \in \mathbb{M}_h$ conditionally on the $H$-neighbourhood of $o$ in $\mathbb{T}$ being the corresponding one of $(i,j)$ in $G_n$. 
Informally, this constructs a "score" matrix where each score is the likelihood of matching $(i,j)$ by looking up to depth $H$. 

The first hurdle is that the $H$- neighbourhood of $(i,j)$ in $G_n$ may not be a tree, so we will consider the universal cover of $G_n$ rooted at $(i,j)$ instead (see below for a definition), and control the resulting error. 

Second, we will show that, by completing the diagonal with the probability that $i$ is unmatched, we create, with a small error, a symmetric stochastic matrix. 

Finally, this means that with a small error as $n \rightarrow \infty$, we created a random symmetric stochastic matrix on $(V_n)^2$ that we decompose into a convex combination of involution matrices with the Birkhoff-Von Neumann theorem, which we reinterpret as a random matching. 

By construction, this random matching on $G_n$ will have the same performance as $(\mathbb{T},\mathbb{M}_h)$ with asymptotically small error, which yields the result. 

\subsubsection{Score function on edges of finite graphs}

\begin{figure}
 \centerline{\includegraphics[scale=0.8]{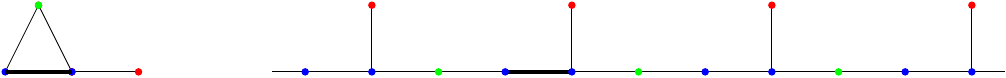}}
\caption{Example of the universal cover of an edge-rooted graph.  \label{fig:cover}}
\end{figure}

Let $(i,j) \in V_n^2$, we define $G_{n,\infty}^{(i,j)}$ as $\emptyset$ if $\{i,j\}$ is not an edge of $G_n$, and as the universal cover $G_{n,\infty}^{(i,j)}$ of $G_n$ rooted at the oriented edge $(i,j)$ if $\{i,j\}$ is an edge of $G_n$. This cover is the tree of non-backtracking trajectories from the edge $(i,j)$. Namely, when $\{i,j\}$ is an edge of $G_n$, $G_{n,\infty}^{(i,j)}$ is a tree rooted at the oriented edge $(i,j)$, children of $i$ (resp. $j$) are the neighbours of $i$ (resp. $j$) in $G_n$, with $j$ (resp. $i$) excluded. The children of any vertex $v$ with parent $v'$ are then the neighbours of $v$ in $G_n$, with $v'$ excluded. See Figure~\ref{fig:cover} for an illustration. For $H \geq 0$, we set $G_{n,H}^{(i,j)}:=B_H(G_{n,\infty}^{(i,j)})$.

We will need the following statement that translates the fact that the universal cover of $G_n$ converges locally to the same limit as $G_n$. 
\begin{lemma}\label{lem:BFSconverge}
    Let $(i,j)$ be chosen uniformly in $E_n$, then:
    \[ G_{n,\infty}^{(i,j)}  \underset{n \rightarrow \infty}{\overset{loc}{\longrightarrow}} \mathbb{T} .   \]
\end{lemma}
\begin{proof}[Proof of Lemma~\ref{lem:BFSconverge}]
    Fix $H>0$, $G_n$ converges locally to $\mathbb{T}$ so:
    \[ \mathbb{P}\left(B_H(G_n,(i,j)) \text{ is a tree} \right) \underset{n \rightarrow \infty}{\longrightarrow} \mathbb{P}\left(B_H(\mathbb{T})\text{ is a tree} \right)=1. \]
    So w.h.p $B_H((G_n,(i,j)))$ is a tree, in which case the universal cover of $G_n$ rooted on $o_n$ up to height $H$ coincides with $G_n$. The conclusion follows.
\end{proof}

Let $\hat{\lambda}$ be the law of the UBGW tree of i.i.d weights $\omega$ and reproduction law $\pi$, let $\hat{\lambda}_n$ be the measure arising from $G_{n,\infty}^{(i,j)}$ with $(i,j)$ picked uniformly in $E_n$. 
We also denote by $\hat{\lambda}_H$ and $\hat{\lambda}_{n,H}$ the corresponding measures when restricted to the space of trees of height $H$. Lemma~\ref{lem:BFSconverge} states that:
\[ \forall H > 0, \hat{\lambda}_{n,H} \overset{weak}{\underset{n \rightarrow \infty}{\longrightarrow}} \hat{\lambda}_H . \]
From $(\mathbb{T},\mathbb{M}_h)$, we define a function on  the space of rooted graphs $\hat{\mathcal{G}}^{\star}$ in the following way:
\begin{equation}\label{eq:defg}    
 g(T)= \mathbb{P}(o \in \mathbb{M} | \mathbb{T}=T ),  \forall T \in \hat{\mathcal{G}}^{\star},
\end{equation}
where it is seen as a Radon-Nikodym density. 

Integrating with respect to the second marginal, we can rewrite $\perfE(\mathbb{T},\mathbb{M}_h)$ as
\begin{equation}\label{treemarginal}
 \perfE(\mathbb{T},\mathbb{M}_h) = \mathbb{E}\left[ w(o)\mathbbm{1}_{o \in \mathbb{M}_h} \right]= \int_{\hat{\mathcal{G}}^{\star}}w(o)g(T)\hat{\lambda}(\mathrm{d}T),  
 \end{equation}
For $H>0$, we define the restricted $g_{H}$ as $g_H(\emptyset) = 0$ and 
\begin{equation}\label{eq:defgH}
 g_{H}(T)=\mathbb{P}(o \in \mathbb{M}_h| B_H(\mathbb{T})=B_H(T)), \forall T \in \hat{\mathcal{G}}^{\star}. 
\end{equation}
The function $g_{H}$, which can be seen as a conditional expectation if we only looked at $T$ up to depth $H$, is exactly what we wanted at the beginning: it returns a matching score for the root if $\mathbb{M}_h$ "sees" $T$ up to depth $H$.
Furthermore, it can be computed explicitly from the law $\zeta_h$, on $B_H(T)$ with an identical construction as in Proposition~\ref{prop:constructionMZ_h} as follows.
Given $B_H(T)$, draw i.i.d variables of law $\zeta_h$ on the $H-$boundary pointing away from the root, then use recursion~\eqref{eq:Z_hrecursion} to define it inside $B_H(T)$, then, writing $o=(o_-,o_+)$, $g_H(T)$ is the probability that $Z(o_-,o_+)+Z(o_+,o_-)\leq w(o)$.

For integrability reasons, we also exclude large weights, for $x \in \mathbb{R}$, define $g_{H,x}=g_H \mathbbm{1}_{w(o)\leq x}$. 
\begin{figure}[t!]
\centering
\includegraphics[scale=0.8]{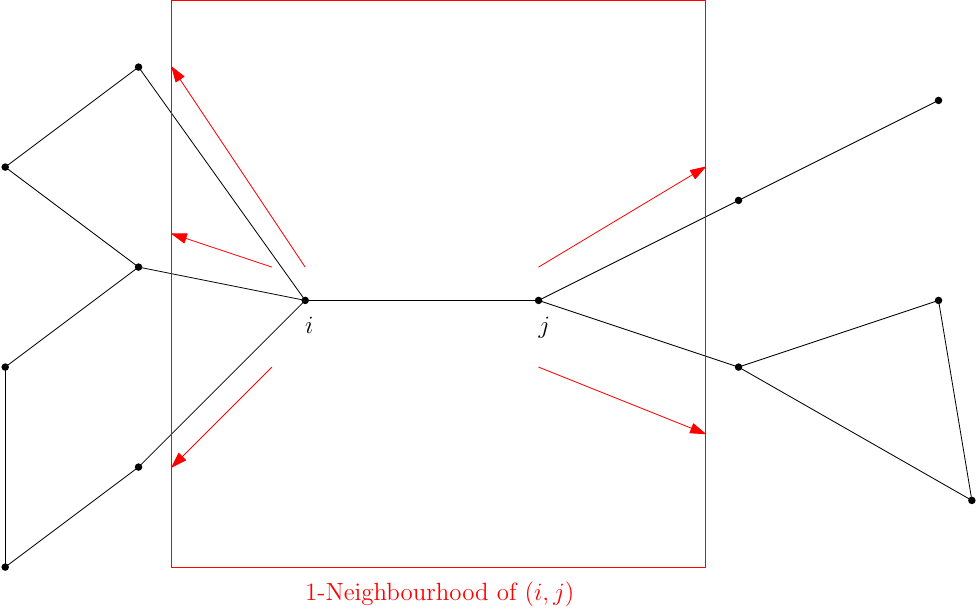}
\caption{Example of a score $q_{i,j,H,x}$ when the $H-$Boundary is a tree. Here, the score $q_{i,j,1,x}$ is the probability that $Z(i,j)+Z(j,i)< w(i,j)\leq x$ where $Z(i,j)$ and $Z(j,i)$ are computed from independent exterior variables $Z$ with law $\zeta_h$ (drawn as red arrows) on the $1-$boundary.}\label{fig:scorematrix}
\end{figure}
We can finally define a random "score" matrix:

\begin{definition}\label{def:scorematrix}
For all $(i,j) \in V_n^2, i \neq j$, define 
\begin{equation}
q_{i,j,H,n,x} := g_{H,x}(G_{n,\infty}^{(i,j)}),
\end{equation}
and 
\begin{equation}
 q_{j,j,H,n,x}:=1-\sum_{i\neq j} q_{j,i,H,n,x}.
\end{equation}
The score matrix is the $|V_n| \times |V_n|$ random matrix defined by  
\begin{equation}
Q_{H,n,x}(G_n):=\left(q_{i,j,H,n,x}\right)_{(i,j) \in V_n^2}.
\end{equation}
\end{definition}
Once again, $q_{i,j,H,n,x}$ can be computed from the distribution $\zeta_h$, see Figure~\ref{fig:scorematrix}.

\bigskip

We will use several intermediary results which proofs are postponed to Subsection~\ref{sec:lemmas}. Let $w_{i,j}$ be the weights of $(i,j)$ in $(V_n)^2$ (we set it as $0$ if the edge is not present).
\begin{lemma}\label{lem:asymptoticperf}
The matrices $Q_{H,n,x}$ are asymptotically optimal:
\begin{equation}\label{eq:eqasymptotic}
  \forall H>0 ,  \lim_{x \rightarrow \infty} \lim_{n \rightarrow \infty} \mathbb{E} \left[\frac{1}{|V_n|}  \sum_{(i,j) \in V_n^2} q_{i,j,H,n,x} w_{i,j}  \right] =  \perfV(\mathbb{T},\mathbb{M}_h).
\end{equation}
\end{lemma}

\begin{lemma}\label{lem:intersectionconv}
The matrices $Q_{H,n,x}$ are asymptotically non-negative:
\begin{equation}\label{eq:stoc}
\forall x \in \mathbb{R}, \lim_{H \rightarrow \infty}\limsup_{n \rightarrow \infty}  \mathbb{E}\left[\frac{1}{|V_n|} \sum_{j \in V_n} \left| q_{j,j,H,n,x}  \right|_{-}  \right] = 0.  
\end{equation}
\end{lemma}
Let $\text{Birk}_d$ be the Birkhoff polytope on $M_d(\mathbb{R})$, $\text{Sym}_d$ the subspace of symmetric matrices, $\text{RStoc}_d$ the right stochastic matrices for which every row sums to 1 and $\text{LStoc}_d$ the left stochastic matrices for which every column sums to 1. 

Lemma~\ref{lem:intersectionconv} controls the distance of the matrix $Q_{H,n,x}$ to $\text{LStoc}_{|V_n|}$. Since $Q_{H,n,x}$ is symmetric, it also controls the distance to $\text{RStoc}_{|V_n|}$. Now, we would like to control the distance to $\text{Birk}_{|V_n|}=\text{LStoc}_{|V_n|} \cap \text{RStoc}_{|v_n|}$.

For this purpose, we need some geometric result.
Write $\|M\|_1=\sum_{i,j}|m_{i,j}|$ for $ M\in M_d(\mathbb{R})$. We will use the following proposition that states that if a matrix $M$ is $o(d)$ close to $\text{LStoc}_d$ and $o(d)$ close to $\text{RStoc}_d$, then it is $o(d)$ close to $\text{Birk}_d$:
\begin{prop}\label{prop:algorithmecommuniste}
Assume $M \in M_d(\mathbb{R})$ has non-negative coefficients. 

Assume that there is some $\frac{1}{2}>\varepsilon >0$ such that:
\[\sum_{i} \left| \sum_{j} m_{i,j}-1 \right| + \sum_{j} \left| \sum_{i} m_{i,j}-1\right| < d\varepsilon. \]
Then there is some constant $C$ that does not depend on $n$ such that:
\[ d_{\|\cdot\|_1}(M,\mathrm{Birk}_d) \leq C\varepsilon d .\]
\end{prop} 
By Proposition~\ref{prop:algorithmecommuniste}, the distance of $Q_{H,n,x}$ to $\text{Birk}_{|V_n|}$ is thus $o(|V_n|)$.
We need another geometric lemma that controls the distance to symmetric bistochastic matrices from the distance to bistochastic matrices and the distance to symmetric matrices:
\begin{lemma}\label{lem:symbirk}
Let $M \in M_d(\mathbb{R})$:
\begin{equation}\label{eq:symbirk}
d(M, \mathrm{Sym}_d \cap \mathrm{Birk}_d ) \leq   d(M, \mathrm{Sym}_d) +d(M,\mathrm{Birk}_d ). 
\end{equation}
\end{lemma}
Combining Equation~\eqref{eq:symbirk} with Lemmas~\ref{lem:asymptoticperf} and \ref{lem:intersectionconv}, the 
 expected distance of $Q_{H,n,x}$ to symmetric bistochastic matrices is $o(|V_n|)$ as $n,H \rightarrow \infty$. 

\subsubsection{Proof of the second part of Proposition~\ref{prop:infinifini}}
Fix $\varepsilon,\eta >0$. 

Applying Lemma $\ref{lem:asymptoticperf}$ ,  there exists $H_{\eta}$, $x_{\eta}$ and $N_{\eta}$ big enough such that for any $n \geq N_{\eta}$,
\[ \mathbb{E}\left[\frac{1}{|V_n|} \sum_{(i,j) \in V_n^2}  q_{i,j,H_{\eta},n,x_{\eta }}w_{i,j}\right] > \perfV(\mathbb{T},\mathbb{M}_h)- \eta.\]
Now, applying what we just showed at the end of the previous subsection, by taking $H_{\eta,\varepsilon}$ and $N_{\eta,\varepsilon}$ big enough, such that for any $n \geq N_{\eta,\varepsilon}$, the previous inequality holds along with
\[ \mathbb{E}\left[\frac{1}{|V_n|} d(Q_{H_{\eta,\varepsilon},n,x_{\eta}},\mathrm{Birk}_{|V_n|} \cap \mathrm{Stoc}_{|V_n|}) \right] \leq \frac{\varepsilon^2}{x_{\eta}} .  \]
Applying Markov inequality,

\[ \mathbb{P}\left[ \frac{1}{|V_n|} d(Q_{H_{\eta,\varepsilon},n,x_{\eta}},\mathrm{Birk}_{|V_n|} \cap \mathrm{Stoc}_{|V_n|}) \geq \frac{\varepsilon}{x_\eta}  \right] \leq \varepsilon. \]

We drop the dependency of $H$ and $x$ in $\eta$ and $\varepsilon$ from now on to ease notation. 
With probability at least $(1-\varepsilon)$, there exists some bistochastic symmetric matrix $\hat{Q}_{H,n,x}$ such that $\|\hat{Q}_{H,n,x}-Q_{H,n,x}\|_1 \leq \frac{\varepsilon}{x_\eta} |V_n|$.
Now apply Birkhoff-Von Neumann theorem on $\hat{Q}_{H,n,x}$ and $\hat{Q}_{H,n,x}^{T}$ to show that $\hat{Q}_{H,n,x}$ is the convex combination of matrices of the form $\frac{P+P^T}{2}$ where $P$ is a permutation matrix. 
This shows that $\hat{Q}_{H,n,x}$ is a convex combination of involution matrices. 

Involution matrices are in bijection with the set of matchings on the complete graph $K_{|V_n|}$. 
Therefore $\hat{Q}_{H,n,x}$ encodes a random matching $M(\hat{Q}_{H,n,x})$ on $G_n$ (we just leave edges unmatched if they are not in $G_n$). On the event where $\hat{Q}_{H,n,x}$ exists, we have
\begin{align*} \perfV(G_n,M(\hat{Q}_{H,n,x}))  &= \mathbb{E} \left[ \frac{1}{|V_n|}\sum_{(i,j)\in V_n^2} \hat{q}_{i,j,H,n,x} w_{i,j} \right] \\
&> \mathbb{E}\left[ \frac{1}{|V_n|}\sum_{(i,j)^2\in V_n^2}q_{i,j,H,n,x}w_{i,j} \right] - \mathbb{E}\left[ \frac{1}{|V_n|}\sum_{(i,j)^2\in V_n^2}|q_{i,j,H,n,x}-\hat{q}_{i,j,H,n,x}|w_{i,j} \right] \\
&> \perfV(\mathbb{T},\mathbb{M}_h)- \eta - \mathbb{E}\left[\frac{1}{|V_n|} \sum_{(i,j)^2\in V_n^2}|q_{i,j,H,n,x}-\hat{q}_{i,j,H,n,x}|w_{i,j} \right] \\
&> \perfV(\mathbb{T},\mathbb{M}_h) - \eta - \mathbb{E}  \left[ \frac{1}{ |V_n|} \sum_{(i,j) \in V_n^2}|q_{i,j,H,n,x}-\hat{q}_{i,j,H,n,x}| w_{i,j}(\mathbbm{1}_{w_{i,j}\leq x})\right]  \\
&> \perfV(\mathbb{T},\mathbb{M}_h) - \eta-x \frac{d(Q_{H,n,x},\hat{Q}_{H,n,x})}{|V_n|}  \\
&> \perfV(\mathbb{T},\mathbb{M}_h) - \eta - x \frac{\varepsilon}{x}\\
&> \perfV(\mathbb{T},\mathbb{M}_h) - \eta - \varepsilon.\\
\end{align*}

We have shown that with probability at least $(1-\varepsilon)$ there exists a random matching $M(\hat{Q}_{H,n,x})$ on $G_n$ whose performance is at least $\perfV(\mathbb{T},\mathbb{M}_h) - \eta - \varepsilon $. In particular, $\mathbb{M}_{\mathrm{opt}}(G_n)$ has to do better. So
\[ \perfV(G_n,\mathbb{M}_\mathrm{opt}(G_n)) \geq (1-\varepsilon)\perfV(G_n,M(\hat{Q}_{H,n,x})) \geq (1-\varepsilon)\left[ 
\perfV(\mathbb{T},\mathbb{M}_h) - \eta - \varepsilon  \right].\]
Taking the limit inferior we get
\[ \liminf_{n \rightarrow \infty} \perfV(G_n,\mathbb{M}_{\mathrm{opt}}(G_n)) \geq \perfV(\mathbb{T},\mathbb{M}_h).   \]
In conclusion, we showed that
\[ \lim_{n \rightarrow \infty}\perfV(G_n,\mathbb{M}_{\mathrm{opt}}(G_n)) = \perfV(\mathbb{T},\mathbb{M}_h),  \]
and as such, all locally convergent subsequences of $(G_n,\mathbb{M}_{\mathrm{opt}}(G_n))$ converge to an optimal matching couple $(\mathbb{T},\mathbb{M})$ as desired. 
\subsubsection{Proof of the technical lemmas}\label{sec:lemmas}

We start by proving a stronger version of Lemma~\ref{lem:asymptoticperf}, where we show the asymptotic correspondence not just for $w_{i,j}$ but for a well-behaved local function.
\begin{definition}
    Let $f: \hat{\mathcal{G}}^{\star} \mapsto \mathbb{R}$. We say that $f$ is a local function if there exists some $H>0$ such that for any $T,T' \in \hat{\mathcal{G}}^{\star}$,
    \[ B_H(T)=B_H(T') \Rightarrow f(T)=f(T'). \]
\end{definition}

\begin{prop}\label{prop:asymptoticlocalfunction}
Let $f$ be a non-negative local function. For every ${i,j} \in E_n$, set $f_{i,j}=f(G_n,(i,j))$, which is f evaluated on $G_n$ rooted at $(i,j)$ and $f_{i,j}=0$ if $\{i,j\} \notin E_n$. Assume $f$ is either bounded or Lipschitz with respect to the weight of the root $w(o)$, then:
\begin{equation}
  \exists H_f \in \mathbb{N},  \forall H>H_f , \lim_{x \rightarrow \infty} \lim_{n \rightarrow \infty} \mathbb{E} \left[ \frac{1}{|V_n|} \sum_{(i,j) \in V_n^2} q_{i,j,H,n,x} f_{i,j}  \right] = \frac{\mathbb{E}_{N \sim \pi}[N]}{2} \mathbb{E}\left[ f(\mathbb{T},o)\mathbbm{1}_{o \in \mathbb{M}_h} \right].
\end{equation}
\end{prop}
Lemma \ref{lem:asymptoticperf} follows by taking $f(G,o)=w(o)$ that is local with $H=0$. 
\begin{proof}[Proof of Proposition \ref{prop:asymptoticlocalfunction}]

Let $H_f$ be the integer arising from the locality of $f$.
First we will condition on $(i,j)$ being in $E_n$, else $f_{i,j}=0$. 
\[ \mathbb{E} \left[ \frac{1}{|V_n|} \sum_{(i,j) \in V_n^2} q_{i,j,H,n,x} f_{i,j}\mathbbm{1}_{(i,j) \in E_n}  \right] =  \mathbb{E} \left[\frac{2|E_n|}{|V_n|} \frac{1}{2|E_n|} \sum_{(i,j) \in V_n^2} q_{i,j,H,n,x} f_{i,j}\mathbbm{1}_{(i,j) \in E_n}  \right] \]
Notice that $2|E_n|$=$|\overset{\rightarrow}{E}_n|$. We will now integrate with respect to $\hat{\lambda}_{n,H}$ for $H > H_f$: 

\[ \mathbb{E} \left[ \frac{1}{|\overset{\rightarrow}{E}_n|} \sum_{(i,j) \in V_n^2} q_{i,j,H,n,x} f_{i,j}  \right] = \int_{\hat{\mathcal{G}}^{\star}} g_{H}(T)\mathbbm{1}_{w(o)\leq x}f(T)\mathrm{d}\hat{\lambda}_{n,H}(T) \]

Let $M$ be the Lipschitz constant of $f$ on the weights of $B_H(T)$.
As $f\mathbbm{1}_{w(o)\leq x}$ is bounded by $Mx$ and $g_H$ by $1$, we have that
\[ \int_{\hat{\mathcal{G}}^{\star}}g_H(T)\mathbbm{1}_{w(o)\leq x}f(T)\mathrm{d}\hat{\lambda}_{n,H}(T) \underset{n \rightarrow \infty}{\longrightarrow }  \int_{\hat{\mathcal{G}}^{\star}}g_H(T)\mathbbm{1}_{w(o)\leq x}f(T)\mathrm{d}\hat{\lambda}_{H}(T) . \]
On the other hand:
\[\mathbb{E} \left[ \frac{2|E_n|}{|V_n|} \right] \underset{n \rightarrow \infty}{\overset{\mathcal{P}}{\longrightarrow}} \mathbb{E}_{N \sim \pi}[N] \]
So by Slutsky's lemma:
\[ \mathbb{E} \left[\frac{2|E_n|}{|V_n|} \frac{1}{|\overset{\rightarrow}{E}_n|} \sum_{(i,j) \in V_n^2} q_{i,j,H,n,x} f_{i,j}\mathbbm{1}_{(i,j) \in E_n}  \right] \underset{n \rightarrow \infty}{\longrightarrow} \mathbb{E}_{N \sim \pi}[N]\int_{\hat{\mathcal{G}}^{\star}}g_H(T)\mathbbm{1}_{w(o)\leq x}f(T)\mathrm{d}\hat{\lambda}_{H}(T) \]
Now taking $x \rightarrow \infty$ by monotone convergence theorem:
\[ \int_{\hat{\mathcal{G}}^{\star}}g_H(T)\mathbbm{1}_{w(o)\leq x}f(T)\mathrm{d}\hat{\lambda}_{H}(T) \underset{n \rightarrow \infty}{\longrightarrow } \int_{\hat{\mathcal{G}}^{\star}}g_H(T)f(T)\mathrm{d}\hat{\lambda}_{H}(T)  . \]
Since we only look up to height $H$, we can substitute $f$ by $f_H$ so by law of total expectation:
\[\int_{\hat{\mathcal{G}}^{\star}}g_H(T)f(T)\mathrm{d}\hat{\lambda}_{H}(T) = \int_{\hat{\mathcal{G}}^{\star}}g(T)f(T)\mathrm{d}\hat{\lambda}(T) = \mathbb{E}[f(\mathbb{T},o)\mathbbm{1}_{o \in \mathbb{M}}]. \]
So the limit is
\[ \mathbb{E}_{N \sim \pi}[N]\mathbb{E}[f(\mathbb{T},o)\mathbbm{1}_{o \in \mathbb{M}}].\]
When $f(T,o)=w(o)$, then we recover
\[ \mathbb{E}_{N \sim \pi}[N]\mathbb{E}[w(o)\mathbbm{1}_{o \in \mathbb{M}}] = \mathbb{E}_{N \sim \pi}[N]\perfE[\mathbb{T},\mathbb{M}_h] = \perfV[\mathbb{T},\mathbb{M}_h]  \]
by Proposition $\ref{prop:chgmtpdv}$. 
\end{proof}
\begin{proof}[Proof of Lemma \ref{lem:intersectionconv}]

For this equation, we adopt the vertex-rooted point of view, let $\lambda_n$ the law of uniformly vertex rooted $G_n$ and $\lambda$ the law of vertex rooted UBGW, then:
\[\mathbb{E}\left[\frac{1}{|V_n|} \sum_{j \in V_n} \left| q_{j,j,H,n,x}  \right|_{-}  \right] = \int_{\mathcal{G}^{\star}} \left|1-\sum_{v \sim o}g_H(T,(o,v))\mathbbm{1}_{w(o)\leq x}\right|_-\mathrm{d}\lambda_{n,H}(T)    .\]
Neglecting $\mathbbm{1}_{w(o)\leq x}$ can only increase the negative part, so: 

\[\int_{\mathcal{G}^{\star}} \left|1-\sum_{v \sim o}g_H(T,(o,v))\mathbbm{1}_{w(o)\leq x}\right|_-\mathrm{d}\lambda_{n,H}(T) \leq \int_{\mathcal{G}^{\star}} \left|1-\sum_{v \sim o}g_H(T,(o,v))\right|_-\mathrm{d}\lambda_{n,H}(T) . \] 
Now with a proof almost identical to the one of Lemma~\ref{lem:BFSconverge}, we can prove that for any $H$, $\lambda_{n,H}$ converges weakly to $\lambda_H$. 
Since $\mathbb{E}_{N \sim \pi}[N]< \infty$, we can first integrate on $\deg(o)$ to bound the integrand by $\mathbb{E}_{N \sim \pi}[N]$. 
So as $n \rightarrow \infty$, we get:
\[ \int_{\mathcal{G}^{\star}} \left|1-\sum_{v \sim o}g_H(T,(o,v))\right|_-\mathrm{d}\lambda_H(T) \]
Now:
\[ \sum_{v \sim o}g(T,(o,v))= \mathbb{P}\left(o \text{ is matched by }\mathbb{M} | \mathbb{T}=T \right)\leq 1 \, \, \, \, \lambda-\text{almost surely on }\mathcal{G}^{\star}\]
So applying $L^1$ martingale convergence theorem as $H \rightarrow \infty$, we recover:
\[ \int_{\mathcal{G}^{\star}}\left|1- \sum_{v \sim o}g_H(T,(o,v)) \right|_-\mathrm{d}\lambda(T)=0.  \]
\end{proof}

\begin{proof}[Proof of Proposition \ref{prop:algorithmecommuniste}]
We will modify $M$ into a bistochastic matrix with a load-balancing algorithm. The idea is to redistribute the load of the heaviest rows to the lightest rows while preserving the columns, then redistribute the load of the heaviest columns to the lightest columns while preserving the rows. At the end of such an algorithm, every row and every column sums to the average of the columns or rows at the start.

Set $M^{(0)}=M$, $L_i^{(0)}=\sum_{j} m_{i,j}$, $C_j^{(0)}=\sum_{i} m_{i,j}$. For $t \leq d$, we will do the following: 
\begin{itemize}
\item Rank $L_i^{(t)}$ such that $L_{k_1}^{(t)} = \ldots = L_{k_p}^{(t)} > L_{k_{p+1}}^{(t)}
\geq \ldots \geq
L_{k_n-l}^{(t)} = \ldots = L_{k_n}^{(t)}$. 

\item Set $(L_i^{(t+1)},C_j^{(t+1)},m_{i,j}^{(t+1)}) \leftarrow (L_i^{(t)},C_j^{(t)},m_{i,j}^{(t)})$ for now.

\end{itemize}
Now, while $L_{k_p}^{(t+1)}> L_{k_{p+1}}^{(t+1)}$. 
As each row has a positive sum, there are some coefficients $m_{k_1,j_1'}^{(t+1)},m_{k_2,j_2'}^{(t+1)},...,m_{k_p,j_p'}^{(t+1)}>0$. 

We will redistribute a fraction of these coefficients along their respective columns to the rows with the smallest sums. Note that by doing so, the sum of each column stays unchanged. However, we want the distribution to:
\begin{enumerate}
    \item never distribute more than one of the $m_{k_i,j_i'}^{(t+1)}$, 
    \item not exceed $L_{k_{p+1}}^{(t+1)} - L_{k_{p+1}}^{(t)}$ otherwise $p$ has to be increased,
    \item not exceed $L_{k_{n-l-1}}^{(t+1)}-L_{k_{n-l}}^{(t+1)}$ otherwise $l$ has to be increased.
\end{enumerate}
The maximal weight that can be distributed is then
\[
\alpha^{(t+1)}:=\min\left(\min_{1 \leq i \leq p}\left( m_{k_i,j_i'}^{(t+1)}\right), L_{k_p}^{(t+1)}-L_{k_{p+1}}^{(t+1)},\frac{p}{l+1} \left[L_{k_{n-l-1}}^{(t+1)}- L_{k_{n-l}}^{(t+1)}   \right]\right).\]

Thus, we set $m_{k_i,j_i'}^{(t+1)} \leftarrow m_{k_i,j_i'}^{(t+1)}-\alpha^{(t+1)}$ and $L_{k_i}^{(t+1)}\leftarrow L_{k_i}^{(t+1)}-\alpha^{(t+1)}$ for $1 \leq i \leq p$  , $m_{k_i,j_i'}^{(+1)} \leftarrow m_{k_i,j_i'}^{(+1)} + \frac{p}{l+1} \alpha^{(t+1)}$ and $L_{k_i}^{(t+1)} \leftarrow L_{k_i}^{(t+1)}+\frac{p}{l+1} \alpha^{(t+1)}$ for $d-l \leq i \leq d$. 

As long the second condition is not satisfied, we recompute $l$, redraw coefficients $m$ and repeat. This is possible as long as the sum of each line is positive, which is the case since we only increased the minimum of the $L_i^{(t+1)}$ away from zero. 
At every step, the algorithm increases $p$ by at least 1, and at step $n$ every line sums to the same amount.
Furthermore, it has not changed the value of any $C_j$.

We then apply the same algorithm to the columns in order to obtain a matrix $\Tilde{M}$ whose rows and columns sum to the same amount $\hat L$, which is $\hat L:= \frac{1}{d} \sum_{i,j} m_{i,j}$.  Finally we output $S:=\frac{1}{\hat{L}} \Tilde{M}$ that is a bistochastic matrix. 

\bigskip

It remains to bound $ \|S-M\|_1$.
Tracking the values of $L_i$ over the course of the algorithm, we see that the total modification to its coefficients is $L_i-\hat{L}$ during the first step and then $C_j-\hat{L}$ during the second step. 
So the total amount of modification we did before dividing by $\hat{L}$ is:
\begin{equation} \label{eq:firstsum}
\sum_{i}|L_i^{(0)}-\hat{L}| + \sum_{j}|C_j^{(0)}-\hat{L}|.
\end{equation}
When we divide by $\hat{L}$, we make a modification of at most 
\begin{equation} \label{eq:secsum}
\sum_{i,j}\Tilde{m}_{i,j} \left|1-\frac{1}{\hat{L}}\right|.
\end{equation}
First, we see that $\hat L$ is close to $1$:
\begin{align*}
    |\hat{L}-1|&=\left| \frac{1}{d}\sum_{i}L_i-1  \right|
    = \frac{1}{d} \left| \sum_{i} L_i-\sum_{i}\frac{1}{d} \right| \\ 
    &\leq \frac{1}{d} \sum_{i} \left|L_i-1\right| \\
    &\leq \varepsilon.
\end{align*}
We can now bound \eqref{eq:firstsum}:
\begin{align*}
 &\sum_{i}|L_i^{(0)}-\hat{L}| + \sum_{j}|C_j^{(0)}-\hat{L}| \\
 &\leq \sum_{i}\left(|L_i^{(0)}-1|+|1-\hat{L}|\right) + \sum_{j} \left(|C_j^{(0)}-1|+|1-\hat{L}| \right). \\
 &< \sum_{i}|L_i-1| + \sum_{j}|C_j-1| +2d|1-\hat{L}| \\
 &< d\varepsilon + 2d \varepsilon = 3d\varepsilon.
 \end{align*}

For \eqref{eq:secsum}, since $\frac{1}{2}>\varepsilon >0$, we have 
$|1-\frac{1}{\hat{L}}|=|\frac{1-\hat{L}}{\hat{L}}|< 4\varepsilon$. Finally, we have that $d\hat{L}=\sum_{i,j} \Tilde{m}_{i,j}  =\sum_{i,j} m_{i,j}^{(0)} < d(1+\varepsilon)$. 

Summing everything up, we established that:
\[ \|S-M\|_1 <  3d\varepsilon + 4d\varepsilon(1+\varepsilon) < 12d\varepsilon\]
and $C=12$ works. 
\end{proof}

\begin{proof}[Proof of Lemma~\ref{lem:symbirk}]
This is true because the Birkhoff polytope is symmetric through the subspace of symmetric matrices.
\begin{align*}
\inf_{S \text{ bistochastic symmetric}} \|M-S\|_1 &\leq \left\|M-\frac{M+M^{T}}{2} \right\|_1 + \inf_{S \text{ bistochastic symmetric}} \left\| \frac{M+M^T}{2}-S \right\|_1 \\
&\leq d(M,\text{Sym}_d) + \inf_{S \text{ bistochastic symmetric}} \left\| \frac{M+M^T}{2}-\frac{S+S^{T}}{2} \right\|_1.
\end{align*}
Now take $S^{\star}$ such that $d(M,\text{Birk}_n)=\|M-S^{\star}\|_1$, then $\frac{S^{\star}+S^{\star T}}{2}$ is bistochastic symmetric so:
\begin{align*}
\inf_{S \text{ bistochastic symmetric}} \|M-S\|_1 &\leq d(M,\text{Sym}_d) + \inf_{S \text{ bistochastic symmetric}} \left\| \frac{M+M^T}{2}-\frac{S+S^{T}}{2} \right\|_1 \\
&\leq d(M,\text{Sym}_d) + \left\| \frac{M+M^{T}}{2} - \frac{S^{\star}+S^{\star T}}{2}  \right\|_1 \\
&\leq d(M,\text{Sym}_d) + \left\|\frac{M-S^{\star}}{2}\right\|_1 + \left\|\frac{M^T-S^{\star T}}{2}\right\|_1 \\
&\leq d(M,\text{Sym}_d) + d(M,\text{Birk}_d).
\end{align*}
\end{proof}

\subsection{Combining everything: convergence of optimally matched graphs}


Proposition~\ref{prop:infinifini}
 shows that every subsequence of $(G_n,o_n,M_{\mathrm{opt}}(G_n))$ has a subsequence that converges locally to some  $(\mathbb{T},\mathbb{M})$ where $\mathbb{M}$ is an optimal matching on $\mathbb{T}$. Meanwhile, Theorem~\ref{unicityUBGW} states that the only optimal matching in law on the UBGW tree $\mathbb{T}$ is $\mathbb{M}_h$. So every subsequence of $(G_n,o_n,M_{\mathrm{opt}}(G_n))$ has a subsequence that converges to $(\mathbb{T},\mathbb{M}_h)$. We deduce that $(G_n,o_n,M_{\mathrm{opt}}(G_n))$ converges locally to $(\mathbb{T},\mathbb{M}_h)$. 
 This concludes the proof of Theorem~\ref{maintheorem}.

\section{Uniqueness of the message passing distribution : Theorem \ref{th:uniqueness}}\label{sec:thuniqueness}

We now turn to the proof of Theorem~\ref{th:uniqueness}. 
 Take any solution $h$ to equation $\eqref{eq:hequation}$, the idea is to recover the distribution $\zeta_h$ from families of statistics on the optimal matching $\mathbb{M}_h$ and use the uniqueness of the distribution of optimal matchings established in Proposition~\ref{prop:optimal}.

As a warm-up, we first inspect the simple case where $\mathrm{supp}(\omega)=\mathbb{R}_+$.
 We will look at the probability that the root is matched conditioned on the weight being $x \in \mathbb{R}_+$:
    \[\mathbb{P}\left( o \in \mathbb{M} | w(o)=x  \right) \]
    as a Radon-Nikodym derivative. 
    
From the edge-decision rule of Proposition \ref{prop:constructionMZ_h},   
    \[ \mathbb{P}\left(Z_{h}+Z_{h}' < x\right)=\mathbb{P}\left( o \in \mathbb{M} | w(o)=x  \right)\]
    where $Z_{h},Z_{h}'$ are independent variables of distribution $\zeta_h$. Hence, this quantity does not depend on $h$. 
    This shows that for any $t \in \mathbb{R}$:
    \[ \mathbb{E}\left[e^{itZ_{h}} \right]^2 \]
    does not depend on $h$.
    By continuity the characteristic function is uniquely determined in a neighbourhood of $0$, and outside of pathological cases, must determine the law of $Z_h$.
    
\bigskip     
    This reasoning breaks down when the $\mathrm{supp}(\omega)$ has "holes" as we have no access to the information of $Z$ at $x$ if $x \notin \mathrm{supp}(\omega)$. To deal with generic supports for $\omega$, we will consider conditioning on the weights of edges along a path originating from the root rather than just the weight of the root edge.

\bigskip

\begin{proof}\renewcommand{\qedsymbol}{}
    Let us assume that $\hat{\pi}_1>0$, we will condition on the event that the weight of the root edge is $w_0$ and that the $+$ side of the root is a simple path $v_1,v_2,....,v_H$ of length $H$ of weights $w_1,....,w_H$.
    We will compute the probability that $v_H$ is unmatched in the matching $\mathbb{M}_h$. The total event of conditioning on the sequence of weights and $v_H$ being unmatched is measurable in the matching and graph, Figure \ref{fig:huniqueness} gives a depiction of the situation.
    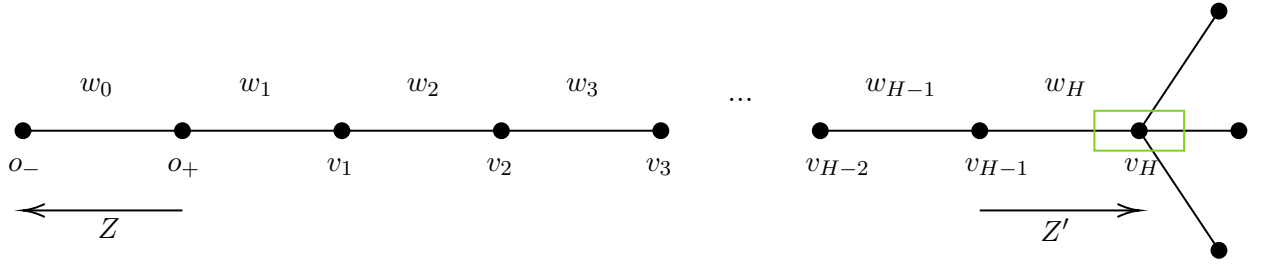
\begin{figure}[t!]
        \centering

\tikzset{every picture/.style={line width=0.75pt}} 

\begin{tikzpicture}[x=0.75pt,y=0.75pt,yscale=-1,xscale=1]

\draw [color={rgb, 255:red, 0; green, 0; blue, 0 }  ,draw opacity=1 ]   (110,120) -- (30,120) ;
\draw  [fill={rgb, 255:red, 0; green, 0; blue, 0 }  ,fill opacity=1 ] (33.73,120) .. controls (33.73,117.94) and (32.06,116.27) .. (30,116.27) .. controls (27.94,116.27) and (26.27,117.94) .. (26.27,120) .. controls (26.27,122.06) and (27.94,123.73) .. (30,123.73) .. controls (32.06,123.73) and (33.73,122.06) .. (33.73,120) -- cycle ;
\draw  [fill={rgb, 255:red, 0; green, 0; blue, 0 }  ,fill opacity=1 ] (113.73,120) .. controls (113.73,117.94) and (112.06,116.27) .. (110,116.27) .. controls (107.94,116.27) and (106.27,117.94) .. (106.27,120) .. controls (106.27,122.06) and (107.94,123.73) .. (110,123.73) .. controls (112.06,123.73) and (113.73,122.06) .. (113.73,120) -- cycle ;
\draw    (110,120) -- (190,120) ;
\draw    (190,120) -- (270,120) ;
\draw    (270,120) -- (350,120) ;
\draw    (430,120) -- (510,120) ;
\draw    (510,120) -- (590,120) ;
\draw    (590,120) -- (630,60) ;
\draw    (590,120) -- (640,120) ;
\draw    (590,120) -- (630,180) ;
\draw  [fill={rgb, 255:red, 0; green, 0; blue, 0 }  ,fill opacity=1 ] (193.73,120) .. controls (193.73,117.94) and (192.06,116.27) .. (190,116.27) .. controls (187.94,116.27) and (186.27,117.94) .. (186.27,120) .. controls (186.27,122.06) and (187.94,123.73) .. (190,123.73) .. controls (192.06,123.73) and (193.73,122.06) .. (193.73,120) -- cycle ;
\draw  [fill={rgb, 255:red, 0; green, 0; blue, 0 }  ,fill opacity=1 ] (273.73,120) .. controls (273.73,117.94) and (272.06,116.27) .. (270,116.27) .. controls (267.94,116.27) and (266.27,117.94) .. (266.27,120) .. controls (266.27,122.06) and (267.94,123.73) .. (270,123.73) .. controls (272.06,123.73) and (273.73,122.06) .. (273.73,120) -- cycle ;
\draw  [fill={rgb, 255:red, 0; green, 0; blue, 0 }  ,fill opacity=1 ] (353.73,120) .. controls (353.73,117.94) and (352.06,116.27) .. (350,116.27) .. controls (347.94,116.27) and (346.27,117.94) .. (346.27,120) .. controls (346.27,122.06) and (347.94,123.73) .. (350,123.73) .. controls (352.06,123.73) and (353.73,122.06) .. (353.73,120) -- cycle ;
\draw  [fill={rgb, 255:red, 0; green, 0; blue, 0 }  ,fill opacity=1 ] (433.73,120) .. controls (433.73,117.94) and (432.06,116.27) .. (430,116.27) .. controls (427.94,116.27) and (426.27,117.94) .. (426.27,120) .. controls (426.27,122.06) and (427.94,123.73) .. (430,123.73) .. controls (432.06,123.73) and (433.73,122.06) .. (433.73,120) -- cycle ;
\draw  [fill={rgb, 255:red, 0; green, 0; blue, 0 }  ,fill opacity=1 ] (513.73,120) .. controls (513.73,117.94) and (512.06,116.27) .. (510,116.27) .. controls (507.94,116.27) and (506.27,117.94) .. (506.27,120) .. controls (506.27,122.06) and (507.94,123.73) .. (510,123.73) .. controls (512.06,123.73) and (513.73,122.06) .. (513.73,120) -- cycle ;
\draw  [fill={rgb, 255:red, 0; green, 0; blue, 0 }  ,fill opacity=1 ] (593.73,120) .. controls (593.73,117.94) and (592.06,116.27) .. (590,116.27) .. controls (587.94,116.27) and (586.27,117.94) .. (586.27,120) .. controls (586.27,122.06) and (587.94,123.73) .. (590,123.73) .. controls (592.06,123.73) and (593.73,122.06) .. (593.73,120) -- cycle ;
\draw  [fill={rgb, 255:red, 0; green, 0; blue, 0 }  ,fill opacity=1 ] (633.73,60) .. controls (633.73,57.94) and (632.06,56.27) .. (630,56.27) .. controls (627.94,56.27) and (626.27,57.94) .. (626.27,60) .. controls (626.27,62.06) and (627.94,63.73) .. (630,63.73) .. controls (632.06,63.73) and (633.73,62.06) .. (633.73,60) -- cycle ;
\draw  [fill={rgb, 255:red, 0; green, 0; blue, 0 }  ,fill opacity=1 ] (643.73,120) .. controls (643.73,117.94) and (642.06,116.27) .. (640,116.27) .. controls (637.94,116.27) and (636.27,117.94) .. (636.27,120) .. controls (636.27,122.06) and (637.94,123.73) .. (640,123.73) .. controls (642.06,123.73) and (643.73,122.06) .. (643.73,120) -- cycle ;
\draw  [fill={rgb, 255:red, 0; green, 0; blue, 0 }  ,fill opacity=1 ] (633.73,180) .. controls (633.73,177.94) and (632.06,176.27) .. (630,176.27) .. controls (627.94,176.27) and (626.27,177.94) .. (626.27,180) .. controls (626.27,182.06) and (627.94,183.73) .. (630,183.73) .. controls (632.06,183.73) and (633.73,182.06) .. (633.73,180) -- cycle ;
\draw  [color={rgb, 255:red, 126; green, 211; blue, 33 }  ,draw opacity=1 ] (567.5,110) -- (612.5,110) -- (612.5,130) -- (567.5,130) -- cycle ;
\draw    (110,160) -- (32,160) ;
\draw [shift={(30,160)}, rotate = 360] [color={rgb, 255:red, 0; green, 0; blue, 0 }  ][line width=0.75]    (10.93,-3.29) .. controls (6.95,-1.4) and (3.31,-0.3) .. (0,0) .. controls (3.31,0.3) and (6.95,1.4) .. (10.93,3.29)   ;
\draw    (510,160) -- (588,160) ;
\draw [shift={(590,160)}, rotate = 180] [color={rgb, 255:red, 0; green, 0; blue, 0 }  ][line width=0.75]    (10.93,-3.29) .. controls (6.95,-1.4) and (3.31,-0.3) .. (0,0) .. controls (3.31,0.3) and (6.95,1.4) .. (10.93,3.29)   ;

\draw (382,102) node [anchor=north west][inner sep=0.75pt]   [align=left] {{\large ...}};
\draw (21,132.4) node [anchor=north west][inner sep=0.75pt]    {$o_{-}$};
\draw (101,132.4) node [anchor=north west][inner sep=0.75pt]    {$o_{+}$};
\draw (181,132.4) node [anchor=north west][inner sep=0.75pt]    {$v_{1}$};
\draw (261,132.4) node [anchor=north west][inner sep=0.75pt]    {$v_{2}$};
\draw (341,132.4) node [anchor=north west][inner sep=0.75pt]    {$v_{3}$};
\draw (421,132.4) node [anchor=north west][inner sep=0.75pt]    {$v_{H-2}$};
\draw (57,92.4) node [anchor=north west][inner sep=0.75pt]    {$w_{0}$};
\draw (137,92.4) node [anchor=north west][inner sep=0.75pt]    {$w_{1}$};
\draw (221,92.4) node [anchor=north west][inner sep=0.75pt]    {$w_{2}$};
\draw (301,92.4) node [anchor=north west][inner sep=0.75pt]    {$w_{3}$};
\draw (451,92.4) node [anchor=north west][inner sep=0.75pt]    {$w_{H-1}$};
\draw (541,92.4) node [anchor=north west][inner sep=0.75pt]    {$w_{H}$};
\draw (501,132.4) node [anchor=north west][inner sep=0.75pt]    {$v_{H-1}$};
\draw (581,132.4) node [anchor=north west][inner sep=0.75pt]    {$v_{H}$};
\draw (66,162.4) node [anchor=north west][inner sep=0.75pt]    {$Z$};
\draw (539,162.4) node [anchor=north west][inner sep=0.75pt]    {$Z'$};

\end{tikzpicture}

        \caption{We look at the event where $v_H$ (in the green box) is unmatched conditionally on the sequence of weights $w_0,....,w_H$.}
        \label{fig:huniqueness}
    \end{figure}

    Now writing $Z=Z_h(o_+,o_-)$ and $Z'=Z_h(v_{H-1},v_H)$, $v_H$ being unmatched is equivalent to:
    \begin{align*}
        &\left\{ \max_{u \sim v_H}(w(v_H,u)-Z_h(v_H,u)) \leq0         \right\} \\
        =&\left\{ Z'=\max\left(0,\max_{u \sim v_H, u \neq v_{H-1}}(w(v_H,u)-Z_h(v_H,u))\right)=0 , w(v_H,v_{H-1})-Z_h(v_H,v_{H-1}) \leq 0  \right\} \\
        =&\left\{ Z'=0, w_H \leq Z_h(v_H,v_{H-1})  \right\}.
    \end{align*}
    By independence we can write that the probability of this event is the product
    \[ \mathbb{P}(Z'=0)\mathbb{P}(w_h<Z_h(v_H,v_{H-1}))  \]
    which does not depend on $h$.
    We will show in Corollary~\ref{coro:perf} that $\mathbb{P}(Z'=0)=h(0)$ only depends on $\mathbb{M}_h$ , which then implies that:
    \begin{align*}
        &\mathbb{P}(w_h \leq Z_h(v_H,v_{H-1})) \\
        =&\mathbb{P}\left( w_H \leq \max\left(0,w_{H-1}-\max\left(0,w_{H-2}-.... -\max\left(0,w_1-\max\left(0,w_0-Z\right) \right)... \right)  \right)    \right).
    \end{align*}
    does not depend on $h$.
    The probability in the previous display can be written as the following integral:
    \[G(w_{0:H}):=G(w_0,...,w_H)= \int_{\mathbb{R}_+} \mathbbm{1}_{\left( w_H \leq \max\left(0,w_{H-1}-\max\left(0,w_{H-2}-.... -\max\left(0,w_1-\max\left(0,w_0-x\right) \right)... \right)  \right)    \right)}  \mathrm{d}\mathbb{P}_h(x)   .\]
    To simplify notations we will denote the function integrated by
    \[g(x,w_{0:H}):=g(x_,w_0,...,w_H)=\mathbbm{1}_{\left( w_H \leq \max\left(0,w_{H-1}-\max\left(0,w_{H-2}-.... -\max\left(0,w_1-\max\left(0,w_0-x\right) \right)... \right)  \right)    \right)}, \]
    so that
    \[ G(w_{0:H})=\int_{\mathbb{R}_+} g(x,w_{0:H})\mathrm{d}\mathbb{P}_h(x). \]
    To prove Theorem~\ref{th:uniqueness}, we will show that for any $x_0 \in [0, \sup(\mathrm{supp}(\omega))]$, for $\varepsilon>0$ small enough, we can recover
     $\mathbb{P}(x_0 \leq Z < x_0+\varepsilon)$
    as a difference $G(w_{0:H})-G(w_{0:H}')$ for $(w_{0:H})$ and $(w_{0:H}')$ with $(w_{0:H})$ chosen according to the following Lemma and $(w_{0:H}')=(w_0,...,w_{H-1},w_H+\varepsilon)  $.
\end{proof}
    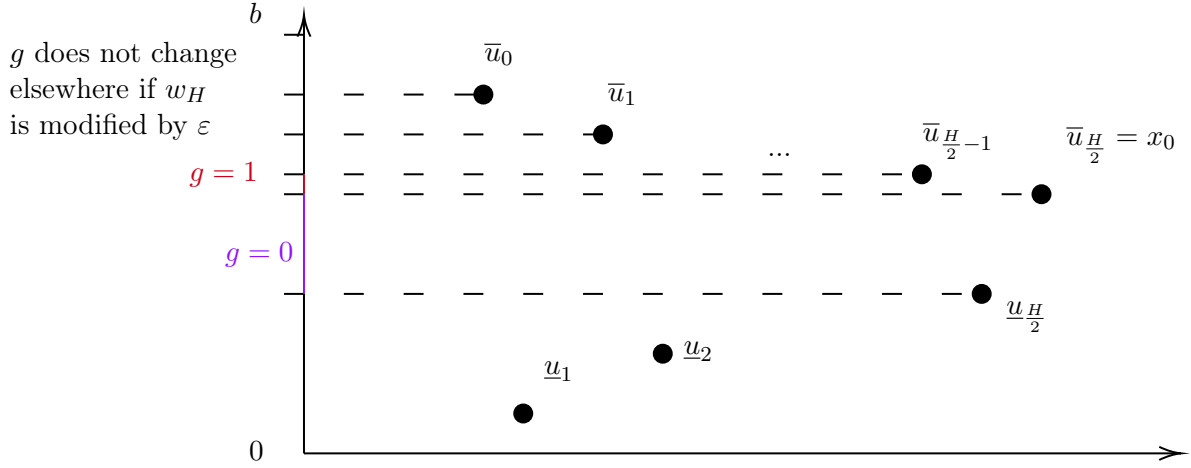
\begin{figure}[t!]
        \centering

\tikzset{every picture/.style={line width=0.75pt}} 

\begin{tikzpicture}[x=0.75pt,y=0.75pt,yscale=-1,xscale=1]

\draw    (150,100) -- (150,22) ;
\draw [shift={(150,20)}, rotate = 90] [color={rgb, 255:red, 0; green, 0; blue, 0 }  ][line width=0.75]    (10.93,-3.29) .. controls (6.95,-1.4) and (3.31,-0.3) .. (0,0) .. controls (3.31,0.3) and (6.95,1.4) .. (10.93,3.29)   ;
\draw    (150,240) -- (588,240) ;
\draw [shift={(590,240)}, rotate = 180] [color={rgb, 255:red, 0; green, 0; blue, 0 }  ][line width=0.75]    (10.93,-3.29) .. controls (6.95,-1.4) and (3.31,-0.3) .. (0,0) .. controls (3.31,0.3) and (6.95,1.4) .. (10.93,3.29)   ;
\draw  [fill={rgb, 255:red, 0; green, 0; blue, 0 }  ,fill opacity=1 ] (244.5,60) .. controls (244.5,57.51) and (242.49,55.5) .. (240,55.5) .. controls (237.51,55.5) and (235.5,57.51) .. (235.5,60) .. controls (235.5,62.49) and (237.51,64.5) .. (240,64.5) .. controls (242.49,64.5) and (244.5,62.49) .. (244.5,60) -- cycle ;
\draw  [fill={rgb, 255:red, 0; green, 0; blue, 0 }  ,fill opacity=1 ] (264.5,220) .. controls (264.5,217.51) and (262.49,215.5) .. (260,215.5) .. controls (257.51,215.5) and (255.5,217.51) .. (255.5,220) .. controls (255.5,222.49) and (257.51,224.5) .. (260,224.5) .. controls (262.49,224.5) and (264.5,222.49) .. (264.5,220) -- cycle ;
\draw  [fill={rgb, 255:red, 0; green, 0; blue, 0 }  ,fill opacity=1 ] (304.5,80) .. controls (304.5,77.51) and (302.49,75.5) .. (300,75.5) .. controls (297.51,75.5) and (295.5,77.51) .. (295.5,80) .. controls (295.5,82.49) and (297.51,84.5) .. (300,84.5) .. controls (302.49,84.5) and (304.5,82.49) .. (304.5,80) -- cycle ;
\draw  [fill={rgb, 255:red, 0; green, 0; blue, 0 }  ,fill opacity=1 ] (334.5,190) .. controls (334.5,187.51) and (332.49,185.5) .. (330,185.5) .. controls (327.51,185.5) and (325.5,187.51) .. (325.5,190) .. controls (325.5,192.49) and (327.51,194.5) .. (330,194.5) .. controls (332.49,194.5) and (334.5,192.49) .. (334.5,190) -- cycle ;
\draw  [fill={rgb, 255:red, 0; green, 0; blue, 0 }  ,fill opacity=1 ] (524.5,110) .. controls (524.5,107.51) and (522.49,105.5) .. (520,105.5) .. controls (517.51,105.5) and (515.5,107.51) .. (515.5,110) .. controls (515.5,112.49) and (517.51,114.5) .. (520,114.5) .. controls (522.49,114.5) and (524.5,112.49) .. (524.5,110) -- cycle ;
\draw  [fill={rgb, 255:red, 0; green, 0; blue, 0 }  ,fill opacity=1 ] (464.5,100) .. controls (464.5,97.51) and (462.49,95.5) .. (460,95.5) .. controls (457.51,95.5) and (455.5,97.51) .. (455.5,100) .. controls (455.5,102.49) and (457.51,104.5) .. (460,104.5) .. controls (462.49,104.5) and (464.5,102.49) .. (464.5,100) -- cycle ;
\draw  [fill={rgb, 255:red, 0; green, 0; blue, 0 }  ,fill opacity=1 ] (494.5,160) .. controls (494.5,157.51) and (492.49,155.5) .. (490,155.5) .. controls (487.51,155.5) and (485.5,157.51) .. (485.5,160) .. controls (485.5,162.49) and (487.51,164.5) .. (490,164.5) .. controls (492.49,164.5) and (494.5,162.49) .. (494.5,160) -- cycle ;
\draw    (235.5,60) -- (225.5,60) ;
\draw    (210,60) -- (200,60) ;
\draw    (180,60) -- (170,60) ;
\draw    (260,80) -- (270,80) ;
\draw    (240,80) -- (230,80) ;
\draw    (210,80) -- (200,80) ;
\draw    (180,80) -- (170,80) ;
\draw    (150,80) -- (140,80) ;
\draw    (150,60) -- (140,60) ;
\draw    (150,30) -- (140,30) ;
\draw    (500,110) -- (510,110) ;
\draw    (480,110) -- (470,110) ;
\draw    (450,110) -- (440,110) ;
\draw    (420,110) -- (410,110) ;
\draw    (380,110) -- (390,110) ;
\draw    (360,110) -- (350,110) ;
\draw    (330,110) -- (320,110) ;
\draw    (300,110) -- (290,110) ;
\draw    (260,110) -- (270,110) ;
\draw    (240,110) -- (230,110) ;
\draw    (210,110) -- (200,110) ;
\draw    (180,110) -- (170,110) ;
\draw    (440,100) -- (450,100) ;
\draw    (420,100) -- (410,100) ;
\draw    (390,100) -- (380,100) ;
\draw    (360,100) -- (350,100) ;
\draw    (320,100) -- (330,100) ;
\draw    (300,100) -- (290,100) ;
\draw    (270,100) -- (260,100) ;
\draw    (240,100) -- (230,100) ;
\draw    (200,100) -- (210,100) ;
\draw    (180,100) -- (170,100) ;
\draw    (150,100) -- (140,100) ;
\draw    (150,110) -- (140,110) ;
\draw    (290,80) -- (300,80) ;
\draw    (390,160) -- (380,160) ;
\draw    (350,160) -- (360,160) ;
\draw    (330,160) -- (320,160) ;
\draw    (300,160) -- (290,160) ;
\draw    (270,160) -- (260,160) ;
\draw    (230,160) -- (240,160) ;
\draw    (210,160) -- (200,160) ;
\draw    (180,160) -- (170,160) ;
\draw    (150,160) -- (140,160) ;
\draw    (480,160) -- (470,160) ;
\draw    (450,160) -- (440,160) ;
\draw    (410,160) -- (420,160) ;
\draw [color={rgb, 255:red, 208; green, 2; blue, 27 }  ,draw opacity=1 ]   (150,100) -- (150,110) ;
\draw [color={rgb, 255:red, 144; green, 19; blue, 254 }  ,draw opacity=1 ]   (150,110) -- (150,160) ;
\draw    (150,160) -- (150,240) ;

\draw (381,88) node [anchor=north west][inner sep=0.75pt]   [align=left] {...};
\draw (121,12.4) node [anchor=north west][inner sep=0.75pt]    {$b$};
\draw (121,232.4) node [anchor=north west][inner sep=0.75pt]    {$0$};
\draw (239,32.4) node [anchor=north west][inner sep=0.75pt]    {$\overline{u}_{0}$};
\draw (269,192.4) node [anchor=north west][inner sep=0.75pt]    {$\underline{u}_{1}$};
\draw (301,52.4) node [anchor=north west][inner sep=0.75pt]    {$\overline{u}_{1}$};
\draw (339,182.4) node [anchor=north west][inner sep=0.75pt]    {$\underline{u}_{2}$};
\draw (459,72.4) node [anchor=north west][inner sep=0.75pt]    {$\overline{u}_{\frac{H}{2} -1}$};
\draw (531,74.4) node [anchor=north west][inner sep=0.75pt]    {$\overline{u}_{\frac{H}{2}} =x_{0}$};
\draw (501,162.4) node [anchor=north west][inner sep=0.75pt]    {$\underline{u}_{\frac{H}{2}}$};
\draw (1,32) node [anchor=north west][inner sep=0.75pt]   [align=left] {$g$ does not change \\ elsewhere if $w_H$ \\ is modified by $\varepsilon$};
\draw (109,132.4) node [anchor=north west][inner sep=0.75pt]  [color={rgb, 255:red, 144; green, 19; blue, 254 }  ,opacity=1 ]  {$g=0$};
\draw (91,92.4) node [anchor=north west][inner sep=0.75pt]  [color={rgb, 255:red, 208; green, 2; blue, 27 }  ,opacity=1 ]  {$g=1$};

\end{tikzpicture}

        \caption{Illustration of Lemma \ref{lem:huniqueness} and its consequences.}
        \label{fig:huniquenesslemma}
    \end{figure}
    
    \begin{lemma}\label{lem:huniqueness}
        Fix $a,b \in \mathbb{R}_+$ with $a<b$, $x_0 \in ]0,b[$.
        There exists $H>0$, even if $x_0>\frac{b}{2}$ and odd if not, and a sequence $(w_0,w_1,....,w_{H})\in [a,b]^{H+1}$ satisfying the following:
        \begin{enumerate}
            \item \[\sum_{i=0}^{H}(-1)^i w_{H-i}=x_0  .\]
  \item The sequence $\overline{u}_k:=\sum_{i=0}^{2k}(-1)^i w_i$ is positive decreasing.
        \item The sequence $\underline{u}_k:=\sum_{i=0}^{2k+1}(-1)^i w_i$ is increasing.
            \item For any $1 < k \leq H$, $\sum_{i=k}^{H} (-1)^iw_{H-i} \neq 0$.
            \item Finally, $\overline{u}_{\lfloor\frac{H}{2}\rfloor} > \underline{u}_{\lfloor \frac{H}{2} \rfloor}$.
        \end{enumerate}
    \end{lemma}

    \begin{proof}[Proof of Lemma \ref{lem:huniqueness}]
        Let's first show the lemma while disregarding condition 4.
        Let us assume $x_0 > \frac{b}{2}$, the other case being symmetrical.
        Take $N=\lfloor\frac{(b-x_0)}{b-a}\rfloor$, for $k<N$, pick $w_{2k}=a $ and $w_{2k+1}=b$.
        Then $\overline{u}_k=b-k(b-a)$ is in decreasing and $\underline{u}_k=k(b-a)$ is increasing.
        Finally, pick $w_{2N}= x_0+\underline{u}_{N} = x_0+\lfloor\frac{(b-x_0)}{b-a}\rfloor (b-a) $.
        This choice is possible as $(b-x_0)-(b-a) < \underline{u}_{N}\leq (b-x_0)$.
        Then $\overline{u}_N=w_{2N}-\underline{u}_{N}=x_0$ by construction, and $u_N < u_{N-1}$.
        The assumption $x_0 > \frac{b-a}{2}$ also implies that $\overline{u}_N >\underline{u}_{N}$.

        Now if we want to obtain condition 4., we just need to pick $w_{2k}=a+\epsilon_k$ and $w_{2k+1}=b-\epsilon_k'$ with $(\epsilon_k, \epsilon_k')$ very small but linearly independent over $\mathbb{Z}$. The distance loss can just be recovered by doing two more steps on $N$. 
    \end{proof}

\bigskip

We return to the proof of Theorem~\ref{th:uniqueness}.

\begin{proof}[Continuation of proof of Theorem~\ref{th:uniqueness}]
       By assumption on $\omega$ take $b=\sup(\mathrm{supp}(\omega))$ and $a$ such that $[a,b] \subset \mathrm{supp}(\omega)$.
    Now for any $h$, take any $x_0 \in [0,b[$ and take $H>0$ and a sequence $(w_0,....,w_H)$ given by Lemma $\ref{lem:huniqueness}$.
    If $b=+\infty$ then do the same reasoning with an arbitrary large $b>0$.
    We will reuse the notations $\overline{u}_k$ and $\underline{u}_k$ from Lemma $\ref{lem:huniqueness}$.
    We will only treat the case when $x_0 > \frac{b}{2}$ and $H$ even, the other case being similar.
    
    Recall $g(x,w_0,...,w_H)=g(x_,w_{0:H})=\mathbbm{1}_{\left( w_H \leq \max\left(0,w_{H-1}-\max\left(0,w_{H-2}-.... -\max\left(0,w_1-\max\left(0,w_0-x\right) \right)... \right)  \right)    \right)}$.
    Decompose the invariant integral we found earlier:
    \begin{align*}
    \int_{\mathbb{R}_+} g(x,w_{0:H})  \mathrm{d}\mathbb{P}_h(x)
    &= \int_{0}^{\underline{u}_0} g(x,w_{0:H}) \mathrm{d}\mathbb{P}_h(x) +\sum_{k=0}^{\frac{H}{2}-1} \int_{\underline{u}_k}^{\underline{u}_{k+1}} g(x,w_{0:H}) \mathrm{d}\mathbb{P}_h(x)  \\
     \qquad &+ \int_{\underline{u}_{\frac{H}{2}-1}}^{\overline{u}_{\frac{H}{2}}} g(x,w_{0:H}) \mathrm{d}\mathbb{P}_h(x) \\
    \qquad &+ \sum_{k=0}^{\frac{H}{2}-1} \int_{\overline{u}_{\frac{H}{2}-k}}^{\overline{u}_{\frac{H}{2}-1-k}} g(x,w_{0:H}) \mathrm{d}\mathbb{P}_h(x) + \int_{\overline{u}_0}^{b} g(x,w_{0:H}) \mathrm{d}\mathbb{P}_h(x).
    \end{align*}
    Our goal is to show that:
    \begin{enumerate}
        \item $g$ evaluates to $0$ on the middle interval: \[\forall x \in \left(\underline{u}_{\frac{H}{2}-1},\overline{u}_\frac{H}{2} \right),  \, g(x,w_{0:H})=0.\]
        \item $g$ evaluates to $1$ on the interval just above: \[ \forall x \in \left[\overline{u}_{\frac{H}{2}},\overline{u}_{\frac{H}{2}-1}\right),  g(x,w_{0:H})=1.\]
        \item If we modified $w_H$ to $w_H+\varepsilon$ for \[\varepsilon < \min_{k \leq H} \left| \sum_{i=k}^{H} (-1)^iw_{H-i} \right| ,\]
        then none of the $g(x,w_{0:H})$'s value would change on the other intervals.
    \end{enumerate}
    First, in the middle, we have the system of inequalities:
    \begin{align}
        \forall k \leq \frac{H}{2}, \overline{u}_k&=\sum_{i=0}^{2k}(-1)^i w_{i}>x ,\\
        \forall k \leq \frac{H}{2}, \underline{u}_k&=\sum_{i=0}^{2k-1}(-1)^{i+1}w_{i}<x .
    \end{align}
    Looking at the inequality for $\overline{u}_0=w_0>x$ it means that the last nested $\max(0,w_0-x)$ intervening in $g(x,w_{0:H})$ evaluates at $w_0-x$.
    
    Then looking at the next inequality $\underline{u}_1=w_0-w_1<x$, it is equivalent to $w_1-w_0+x>0$ so the next nested $\max(0,w_1-w_0+x)$ evaluates at $w_1-w_0+x$.
    
    Unravelling the maximums, we end up with the inequality inside $g$ being:
    \[w_H \leq \sum_{i=0}^{H-1}(-1)^{i} w_{H-1-i}+x  . \]
    But this is precisely the opposite of $\overline{u}_\frac{H}{2} > x$, hence $g(x,w_{0:H})=0$.

    Second, if $x$ is between $\overline{u}_{\frac{H}{2}}$ and $\overline{u}_{\frac{H}{2}-1}$, then the unravelling we did in the previous case stays true, the only inequality that becomes reversed is the final one, so we have:
    \[w_H \leq \sum_{i=0}^{H-1}(-1)^{i} w_{H-1-i}+x   \]
    which is the event on which $g=1$, hence $g=1$ on this interval.

    Third, if $x$ is not in the previous two intervals, then let $i_0=i_0(x,w_{0:H})$ be the last time the nested maximums evaluate to zero.
    Then $g$ being one is equivalent to:
    \[w_H \leq  \sum_{i=0}^{H-i_0-1} (-1)^{i+1} w_{H-1-i} \]
    and $i_0$'s value does not depend on the last values of $w_i$, so it does not depend on $w_H$ on those intervals (one can see $i_0$ as a stopping time on $x$ and $w_i$).
    By condition 4. of the lemma, none of these inequalities change when we modify $w_H$ by the prescribed $\varepsilon$, so $g(x,w_{0:H})$'s value does not change either.

    To conclude, the only change in the integral occurs around $x_0$, hence for all $\varepsilon$ small enough,
    \[ \mathbb{P}_{Z_{h} \sim \zeta_{h}}(x_0 \leq Z_{h} < x_0+\varepsilon)     .\]
    does not depend on $h$.
    This being true for all $x_0 \in \mathrm{supp}(Z_{h})$, we have thus shown that the law of $h$ is unique and that the density of $h$ at $x_0$ can be recovered by as $\frac{\partial G(w_0,...,w_H)}{\partial w_H}$ for correctly chosen $w_0,....,w_{H}$.

    If $\hat{\pi}_1=0$, take the smallest $p \geq 1$ such that $\hat{\pi}_p>0$, then condition on weights such that every path from $o_+$ on the $+$ side sees the same weights as in the case $\hat{\pi}_1>0$. A little bit of analysis shows that a similar proof holds.
\end{proof}

\begin{remark}
    It should not be too hard to lessen the hypothesis on this proof, the open question would be whether it remains true had we picked an entirely singular distribution for $\omega$.
\end{remark}

\section{Applications}\label{sec:applications}

In this section, we show a few applications of Theorem $\ref{maintheorem}$. Namely, we compute asymptotics of local statistics using the local convergence.

\subsection{Optimal matching performance and density}

In this entire subsection, we reuse the previous notation: $(G_n)$ is a sequence of random graphs that converges locally to a UBGW of reproduction law $\pi$ and weight law $\omega$ on the edges with $\pi$ and $\omega$ integrable and $\omega$ atomless. 
Let $\phi$ be the generating function of $\pi$ and $\hat{\phi}$ be $\frac{\phi'}{\phi'(1)}$.
 In particular, $\mathbb{E}_{N \sim \pi}[N]=\phi'(1)$.
The function $h$ is the solution to Equation~\eqref{eq:hequation}, and the family $Z_h(u,v)$ is the associated process on $\mathbb{T}$. Finally, $\zeta_h$ the law for which $h$ is the cumulative distribution function. 

Let $\mathbb{M}_{\mathrm{opt}}(G_n)$ be any optimal matching on $G_n$ then applying the theorem to the local function $f=w(o)\mathbbm{1}_{o \in \mathbb{M}}$ and $f=\mathbbm{1}_{o \in \mathbb{M}}$ gives Corollary~\ref{coro:perf}, which we restate with slightly more detail:

\begin{coro}\label{coro:perf2}
Under the same assumptions as Theorem~\ref{maintheorem}, the asymptotics of the average cost per edge and edge density can be computed on the limiting tree. Let $W$ of law $\omega$, $Z$ and $Z'$ of law $\zeta_h$ such that $(W,Z_h,Z_h')$ are mutually independent. Then:
    \begin{align}
        \lim_{n \rightarrow \infty} \mathbb{E}\left[ \frac{1}{|E_n|} \sum_{e \in \mathbb{M}_{\mathrm{opt}}(G_n)}w(e)  \right] &= \mathbb{E}[w(o)\mathbbm{1}_{o \in \mathbb{M}}] = \mathbb{E}\left[W\mathbbm{1}_{Z_h+Z_h'<W}  \right].\\
        \lim_{n \rightarrow \infty} \mathbb{E}\left[ \frac{1}{|E_n|} \sum_{e \in \mathbb{M}_{\mathrm{opt}}(G_n)}1  \right] &= \mathbb{E}[\mathbbm{1}_{o \in \mathbb{M}}] = \mathbb{E}\left[\mathbbm{1}_{Z_h+Z_h'<W}\right] = \frac{1-\phi(\hat{\phi}^{-1}(h(0)))}{\phi'(1)}  .
    \end{align}
\end{coro} 

\begin{proof}
The only thing to prove is the last equality $\mathbb{E}\left[\mathbbm{1}_{Z_h+Z_h'<W}\right] = \frac{1-\phi(\hat{\phi}^{-1}(h(0)))}{\phi'(1)}$. Indeed, the other identities are direct consequences of Theorem~\ref{maintheorem}. Let us start with a simple calculation:
\begin{align*}
\mathbb{P}(Z_h+Z_h'<W)&= \mathbb{E}_W[P(Z_h+Z_h'<W)] \\ 
&= \mathbb{E}_W\left[\int_\mathbb{R} \mathbb{P}(Z_h+s<W)\mathrm{d}\mathbb{P}_{Z_h'}(s) \right] \\
&=\mathbb{E}_W\left[\int_\mathbb{R} \mathbb{P}(Z_h<W-s)\mathrm{d}h(s)\right] \\
&= \mathbb{E}_W\left[ \int_{\mathbb{R}} h(W-s)\mathrm{d}h(s)    \right] \\
&= \int_{\mathbb{R}} \mathbb{E}_W[h(W-s)]\mathrm{d}h(s).
\end{align*}
Using the fact that $h$ is solution to Equation~\eqref{eq:hequation} we get
\begin{align*}
\mathbb{P}(Z_h+Z_h'<W)
&= \int_{\mathbb{R}_+} (1-\hat{\phi}^{-1}(h(s)))\mathrm{d}h(s).
\end{align*}

Let $g(u)=\hat{\phi}(1-u)$, $g^{-1}(u)=1-\hat{\phi}^{-1}(u)$, $G(u)= -\frac{\phi(1-u)}{\phi'(1)}$ be an antiderivative of $g$.
From there, elementary calculus gives
\begin{align*}
\mathbb{P}(Z_h+Z_h'<W)
&= \int_{\mathbb{R}_+}g^{-1}(h(s)) \mathrm{d}h(s) \\
&= h(0)g^{-1}(h(0)) + \int_{h(0)}^{1} g^{-1}(u)\mathrm{d}u \\
&= h(0)g^{-1}(h(0)) + \left[ ug^{-1}(u) - G(g^{-1}(u)) \right]_{h(0)}^{1} \\
&= -G(0)+ G(g^{-1}(h(0))) \\
&= \frac{\phi(1)}{\phi'(1)} - \frac{\phi(1-(1-\hat{\phi}^{-1}(h(0))))}{\phi'(1)} \\
&= \frac{1-\phi(\hat{\phi}^{-1}(h(0)))}{\phi'(1)}.
\end{align*}
\end{proof}

We can also state the following, which generalises Corollary~\ref{coro:perf2}.

\begin{coro}\label{coro:onesidedegree}
    Under the same assumptions as Theorem~\ref{maintheorem}, the asymptotics of the law of a vertex being matched in $\mathbb{M}_{\mathrm{opt}}(G_n)$ when its degree is conditioned to be $(k+1)$, for a vertex $v$ chosen uniformly in $G_n$, is given by :
\begin{equation}
    \lim_{n \rightarrow \infty} \mathbb{P}\left(v \text{ is matched in } \mathbb{M}_{\mathrm{opt}}(G_n) | \deg(v)=k+1\right) =1- \hat{\phi}^{-1}(h(0))^{k+1}.
\end{equation}    
Equivalently, for the edge-rooted version, let $(u,v)$ be an uniform directed edge of $G_n$:
\begin{equation}\label{eq:onesidedegree}
    \lim_{n \rightarrow \infty} \mathbb{P}\left((u,v) \in \mathbb{M}_{\mathrm{opt}}(G_n) | \deg(v)=k+1\right) = \frac{1- \hat{\phi}^{-1}(h(0))^{k+1}   }{k+1}.
\end{equation}
\end{coro}
\begin{proof}
    
    We will prove the edge-rooted Equation~\eqref{eq:onesidedegree}.
    Convergence is obtained by applying Theorem~\ref{maintheorem} to the local function $f=\mathbbm{1}_{o=(o_-,o_+) \in \mathbb{M}, \deg(o_+)=k }$. The function $f$ is $1-$local so we condition on the $1-$neighbourhood. 
    
    Let us call $Z_{-,i}$ the outwards $Z$ on the minus side and $Z_{+,i}$ the same on the plus side. They are an i.i.d family of law $\zeta_h$. Let us also write $w_{+,i}$ and $w_{-,i}$ the weights of the corresponding edges that are also i.i.d of law $\omega$. Finally, let us write $N_-$ the number of children of $o_-$ and $N_+$ the number of children of $o_+$, the family $(Z_{+,i},w_{+,i},N_{+},Z_{-,i},w_{-,i},N_{-})$ is independent. Refer to Figure~\ref{fig:onesidedegree} for an illustration.
    \begin{figure}
        \centering
        \includegraphics[scale=0.65]{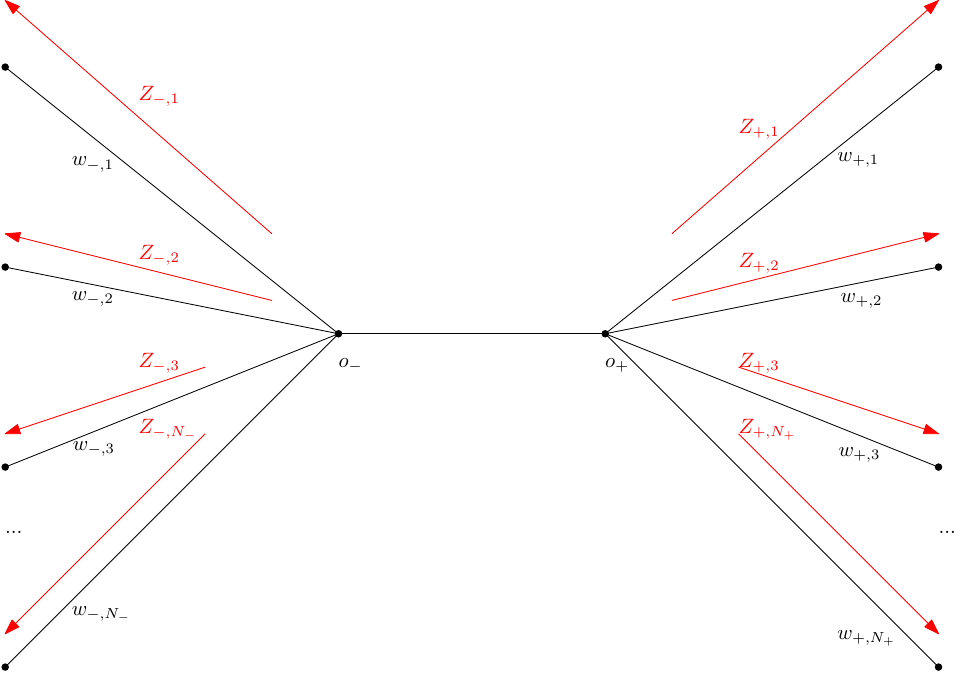}
        \caption{Illustration of the notations in the proof of Corollary~\ref{coro:onesidedegree}.}
        \label{fig:onesidedegree}
    \end{figure}
    The limiting quantity is thus
    \[ \mathbb{E}\left[ \mathbbm{1}_{Z(o_-,o_+)+Z(o_+,o_-)<w(o)} | N_+=k \right]  .\]
    $N_-$ is unrestricted so $Z(o_+,o_-)$ is simply of law $\zeta_h$, whereas $Z(o_-,o_+)=\max_{1 \leq i \leq N_+}(0,w_{+,i}-Z_{+,i})$.
    This leads to computing
    \begin{align*}
    \mathbb{P}\left(Z(o_-,o_+\right)\leq t| N_+=k) &= \mathbb{P} \left( \max_{1 \leq i \leq N_+}(0,w_{+,i}-Z_{+,i}) \leq t  | N_+=k \right) \\
    &= \mathbb{P}\left( t \geq 0\right) \prod_{i=1}^{k} \mathbb{P}\left( w_{+,i}-Z_{+,i} \leq t \right) \\
    &= \mathbbm{1}_{t \geq 0}\left( 1-\mathbb{E}_W[h(W-t)] \right)^k \\
    &=\mathbbm{1}_{t \geq 0 } \hat{\phi}^{-1}(h(t))^k
    \end{align*}
    by noticing that $ 1-\mathbb{E}_W[h(W-t)] = \hat{\phi}^{-1}(h(t))$ according to Equation~\eqref{eq:hequation}. Define:
    \[h_k(t)=\mathbbm{1}_{t \geq 0 } \hat{\phi}^{-1}(h(t))^k . \]
    We have just seen that $h_k$ is the cdf of $Z_{o_-,o_+}$ conditionally on $N_+=k$
    Now we integrate with respect to the value of $Z_{o_-,o_+}$:
    \begin{align*}
        &\mathbb{E}\left[ \mathbbm{1}_{Z(o_-,o_+)+Z(o_+,o_-)<w(o)} | N_+=k \right]  \\
        &= \hat{\phi}^{-1}(h(0))^k\mathbb{P}(Z_{o_+,o_-}< w(o)) + \int_{0}^{1} \mathbb{P}\left( Z_{o_-,o_+}+ s < w(o) \right) h_k'(s)\mathrm{d}s. \\
        &=\hat{\phi}^{-1}(h(0))^{k}\mathbb{E}_W[h(W)]+ \int_{0}^{1}\mathbb{E}_W[h(W-s)]k \frac{h'(s)}{\hat{\phi}'(h(s))} \hat{\phi}^{-1}(h(s))^{k-1}\mathrm{d}s. \\
        &=\hat{\phi}^{-1}(h(0))^{k}\left(1-\hat{\phi}^{-1}(h(0))\right)+ \int_{0}^{1}\mathbb{E}_W\left( 1-\hat{\phi}^{-1}(h(s))\right)k \frac{h'(s)}{\hat{\phi}'(h(s))} \hat{\phi}^{-1}(h(s))^{k-1}\mathrm{d}s. \\
        &=\hat{\phi}^{-1}(h(0))^{k}\left(1-\hat{\phi}^{-1}(h(0))\right) + \int_{\hat{\phi}^{-1}(h(0))}^{1}k(1-u)u^{k-1}\mathrm{d}u \\
        &=\hat{\phi}^{-1}(h(0))^{k}\left(1-\hat{\phi}^{-1}(h(0))\right)  + \left[u^{k}-\frac{k}{k+1}u^{k+1}\right]_{\hat{\phi}^{-1}(h(0))}^{1} \\
        &=\hat{\phi}^{-1}(h(0))^{k}\left(1-\hat{\phi}^{-1}(h(0))\right) + \left(1-\frac{k}{k+1}\right)-\left( \left(\hat{\phi}^{(-1)}(h(0))\right)^{k} -\frac{k}{k+1}\left(\hat{\phi}^{(-1)}(h(0))\right)^{k+1} \right) \\
        &=\frac{1}{k+1}-\frac{1}{k+1}\left(\hat{\phi}^{-1}(h(0))\right)^{k+1} \\
        &= \frac{1-\left(\hat{\phi}^{-1}(h(0))\right)^{k+1}}{k+1}.
    \end{align*}
\end{proof}

We did not manage to find further nice formulas that only depend on $h(0)$. Indeed, if we push the conditioning further, the expressions depend on the entire shape of the function $h$. This is the case for instance for the probability of presence of a gap in the matching:
\begin{coro}\label{coro:probagap}
    Under the same assumptions of Theorem~\ref{maintheorem}, we can compute the asymptotics of the law of a "gap" on one side of an edge of $\mathbb{M}_{\mathrm{opt}}(G_n)$ when the side of the gap is conditioned to be of degree k. Let $h_w$ be the cdf of $\omega$, then:
    \begin{align*}
    &\lim_{n \rightarrow \infty}\mathbb{P}\left( \text{no children of }v \text{ in } T_{(u,v)} \text{ is matched by } \mathbb{M}_{\mathrm{opt}}(G_n) | \deg(v)=k+1, (u,v) \in \mathbb{M}_{\mathrm{opt}}(G_n)  \right) \\
    &=\frac{h(0)^k}{1-\hat{\phi}^{-1}(h(0))^{k}} k\mathbb{E}_{(W,Z)\sim \omega \otimes \zeta_h}[\mathbbm{1}_{W-Z \geq 0}(h_w(W-Z))^k]. 
    \end{align*}
\end{coro}
\begin{proof}[Proof of Corollary \ref{coro:probagap}]
    
    \begin{figure}
        \centering
        \includegraphics[scale=0.65]{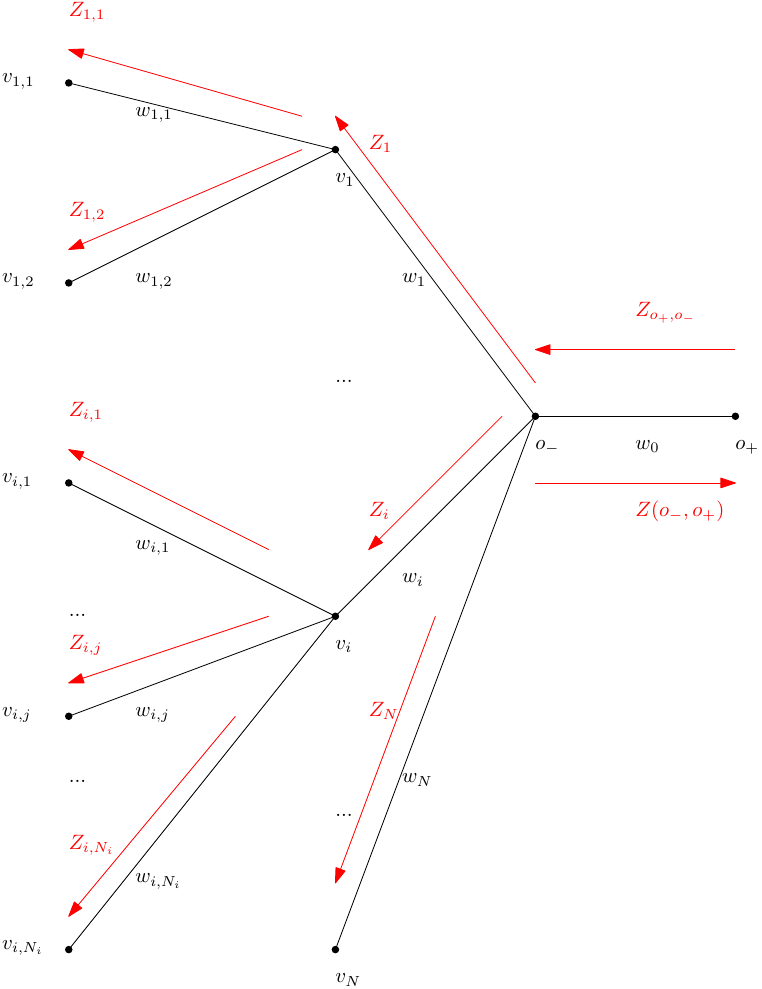}
        \caption{Illustration of the notations used in the proof of Corollary~\ref{coro:probagap}.}
        \label{fig:probagap}
    \end{figure}
    As in the previous corollaries, we apply Theorem~\ref{maintheorem} to the corresponding $2-$local indicator function of the event to obtain the convergence to the corresponding probability on the tree.
    Write $v_i$ the children of $o_-$ in $\mathbb{T}_{(o_+,o_-)}$, $w_i$ the weight of $(o_-,v_i)$,  $v_{i,j}$ the children of $v_i$ and $w_{i,j}$ the weight of $(v_i,v_{i,j})$, $Z_{i,j}=Z(v_i,v_{i,j})$, $Z_i=Z(o_-,v_i)$, $Z_{o_+,o_-}=Z(o_+,o_-)$ and $w_o=w(o_+,o_-)$. Finally, write $N$ the number of children of $o_-$ and $N_i$ the number of children of $v_i$. See Figure~\ref{fig:probagap} for an illustration.
    Fix $k \in \mathbb{N}$, the condition rewrites as:
    \[ B=\left\{Z_{o_+,o_-}+Z(o_-,o_+)< w_o, N=k\right\}  \]
    Now, apply the recursive equation on $Z_{o_+,o_-}$ twice to obtain:
    \begin{align*}
        B&=\left\{ Z(o_+,o_-)+Z(o_-,o_+)< w_o, N=k\right\} \\
        &=\left\{ \max\left(0,\max_{i=1}^{k}\left( w_i-Z_i \right)\right) + Z(o_-,o_+)<w_o,N=k     \right\} \\
        &=\left\{ \max\left(0,\max_{i=1}^{k}\left( w_i-\max\left(0,\max_{i=1}^{N_i}\left(w_{i,j}-Z_{i,j}\right)\right) \right)\right) + Z(o_-,o_+)<w_o,N=k     \right\}.
    \end{align*}
    The event we are interested in is:
    \[ A=\bigcap_{i \leq N} \left\{ \max_{u \sim v_i}(w(v_i,u)-Z(v_i,u))<0    \right\} . \]
    On $B$, we know that $N=k$, but we also know that either $w(v_i,o_-)-Z(v_i,o)<\max_{u \sim v_i}(w(v_i,u)-Z(v_i,u))$, or both are negative, because else, $v_i$ would be matched with $o_-$ who is already matched to $o_+$. So $A \cap B$ can be rewritten as:
    \[\bigcap_{i \leq k} \left\{ \max_{i=1}^{N_i}\left( w_{i,j}-Z_{i,j}\right) <0   \right\} \cap B . \]
    This means that the maximum inside the maximums in the expression in $B$ is simply 0, so $A \cap B$ can be rewritten as :
    \[\bigcap_{i \leq k} \left\{ \max_{i=1}^{N_i}\left( w_{i,j}-Z_{i,j}\right) <0   \right\} \cap \left\{ \max\left(0,\max_{i=1}^{k}w_i \right) +Z(o_-,o_+)< w_o , N=k \right\} .  \]
    Finally, we can notice by the recursive equation on $Z$ that the left event can be rewritten for $i \leq N$ into
    \[ \left\{ Z_i=0 \right\} . \]
    We also notice that the event on the right can be written as
    \[ \left\{ 0<w_o-Z(o_-,o_+), N=k \right\} \cap \bigcap_{i=1}^{k} \left\{ w_i<w_o-Z(o_-,o_+)  \right\} . \]
    In the end, we obtain that:
    \[A \cap B =  \left\{ 0<w_o-Z(o_-,o_+) , N=k   \right\}  \cap \bigcap_{i=1}^{k}\left\{  w_i<w_o-Z(o_-,o_+), Z_i=0\right\}.  \]
    Every variable appearing in this expression are now independent, namely, $w_i,w_o$ of law $\omega$, $N$ of law $\hat{\pi}$, $ Z_{i},Z(o_-,o_+)$ of law $\zeta_h$. So we can compute its probability. 
    \begin{align*} \mathbb{P}(A \cap B)&= \mathbb{P}(N=k) \mathbb{E}_{(W,Z)\sim \omega \otimes \zeta_h}[\mathbbm{1}_{W-Z \geq 0}(h_w(W-Z))^k]      \prod_{i=1}^{k}\mathbb{P}\left( Z_i=0\right)  \\
    &= \hat{\pi}_k h(0)^k \mathbb{E}_{(W,Z)\sim \omega \otimes \zeta_h}[\mathbbm{1}_{W-Z \geq 0}(h_w(W-Z))^k].
    \end{align*}
    By Corollary \ref{coro:onesidedegree}, 
    \[\mathbb{P}(B)= \frac{1-\hat{\phi}^{-1}(h(0))^{k}}{k} \mathbb{P}(N=k)= \frac{ \hat{\pi}_k\left(1-\hat{\phi}^{-1}(h(0))^{k}\right)}{k}.  \]
    In the end,
    \[ \mathbb{P}(A|B)= \frac{h(0)^k}{1-\hat{\phi}^{-1}(h(0))^{k}} k\mathbb{E}_{(W,Z)\sim \omega \otimes \zeta_h}[\mathbbm{1}_{W-Z \geq 0}(h_w(W-Z))^k].\]
\end{proof}

\subsection{Exponential edge weights}
The case where $\omega$ is an exponential law is one of the only cases where $h$ is explicitly computable, and it has been solved by Gamarnik et al. in \cite{gamarnik2003maximum} in the case of Erdös-Renyi graphs and $d-$regular graphs. Our method generalizes their result to a general reproduction distribution as follows.  Let us look at Equation~\eqref{eq:hequation} when $W$ is of law $\text{Exp}(1)$.
\[ h(t)=\mathbbm{1}_{t \geq 0 } \hat{\phi}(1-\mathbb{E}[h(W-t)]) . \]
Write 
\[ \mathbb{E}[h(W-t)]= \int_{\mathbb{R}_+} h(w-t)e^{-w}\mathrm{d}w = \int_{-t}^{+\infty}h(u)e^{-u-t}\mathrm{d}u=e^{-t}\int_{\mathbb{R}_+}h(u)e^{-u}\mathrm{d}u.  
\]
Set \[K=\int_{\mathbb{R}_+}h(u)e^{-u}\mathrm{d}u,\] then
\[h(t)=\mathbbm{1}_{t \geq 0 } \hat{\phi}(1-e^{-t}K), \]
and $K$ has to solve
\[K=\int_{\mathbb{R}_+} \hat{\phi}(1-e^{-u}K)e^{-u}\mathrm{d}u. \]
Now, consider the map $f:x \in [0,1] \mapsto  \int_{\mathbb{R}_+} \hat{\phi}(1-e^{-u}x)e^{-u}\mathrm{d}u $. 
We have
$f(0)=1, f(1)=\int_{\mathbb{R}_+} \hat{\phi}(1-e^{-u})e^{-u}\mathrm{d}u<1$ as long $\hat{\phi}\neq 1$ (in which case the tree is empty),
and $f$ is continuous and strictly decreasing. Thus, there exists a unique $K$ satisfying the equation and the solution $h$ is unique and explicit, as expected.

\section{Extensions} \label{sec:extensions}
To conclude, we now give a few possible generalisations. We only give the main ideas and stay purposely light on details.

\subsection{Multi-type UBGW tree and Stochastic Block Model}

Our results extends to multi-type UBGW trees which also appear in \cite{bordenave2012matchings} but in the unweighted case. These trees appear as local limits for Stochastic Block Models defined as follows.
Let $V_n=\{1,...,n\}$, and let $k \in \mathbb{N}^{\star}$ be the number of types.
Let $(\alpha_1,....,\alpha_k)$ be strictly positive numbers such that $\sum_{i=1}^k \alpha_k=1$.
For $(i,j)\in \{1,...,k\}^2$ fix $c_{i,j}>0$.
Partition $V_n$ into $V_{n,1}...,V_{n,k}$ such that $\left||V_{n,i}|-\alpha_i n\right|<1$, we may have to add one dummy vertex at the end to deal with the fact that $\alpha_i n$ may not be an integer, but it doesn't change the asymptotic behaviour. 
The (sparse) stochastic block model is the random graph generated on $V_n$ so that independently for $x \in V_{n,i}$ and $y \in V_{n,j}$, $\mathbb{P}((x,y)\in E_n) = p_{i,j}=\frac{c_{i,j}}{n}$. 






Then if we follow the same proof, the equation on the $Z(u,v)$ can be rewritten as a system with $k^2$ equations once we condition on the types of $u$ and $v$. It should be possible to use Schauder's fixed point theorem on this system to get existence and then continue with the techniques developed in this work.

\subsection{Vertex weights}
Our results extend to the case where the weights are no longer on the edges but on the vertices instead. It is equivalent to setting $w(u,v)=w(u)+w(v)$ where the $w(v)$ are i.i.d, the weights on edges are no longer independent as they are correlated as soon as they share a vertex. 

Going back to the heuristic:
\[  Z(u,v)=\max(0,\max_{\substack{u' \sim v \\ u' \neq u}} w(v,u')-Z(v,u')) \]
becomes
\[  Z(u,v)=\max(0,\max_{\substack{u' \sim v \\ u' \neq u}} w(v)+w(u')-Z(v,u')) .  \]
To recover independence inside the maximum, we introduce the alternative variables
\[ \hat{Z}(u,v)=Z(u,v)-w(v), \] then the equation becomes:
\begin{equation}\label{eq:Zrecursionvertex}
\hat{Z}(u,v)= \max(-w(v), \max_{\substack{u' \sim v \\ u' \neq u}} -\hat{Z}(v,u')).   
\end{equation}
Now the list of variables inside the maximum is independent, but we can see that $\hat{Z}(u,v)$ is correlated with $w(v)$ with the correlation structure appearing in the equation. 
This leads to the following RDE over a law $\zeta$: For $W$ of law $\omega$, $N$ of law $\hat{pi}$ and $Z,Z_i$ of law $\zeta$, all independent:
\begin{equation}\label{eq:Zvertexweightlaw}
Z \overset{\mathrm{(law)}}{=} \max(-W,\max_{1 \leq i \leq N}(-Z_i))
\end{equation}
Similarly as before we can define the map:
\begin{alignat*}{2}
\mathrm{F}\colon X&\rightarrow X\\
           f&\mapsto \mathrm{F}(f)\colon {}&\mathbb{R} &\rightarrow [0,1] \\
               &&t & \mapsto  {} \mathbb{P}_{W\sim \omega}(t \geq W) \hat{\phi}(1- f(-t))
\end{alignat*}
with $f$ the lower continuous version of $f$ and use Schauder's theorem to recover a solution to the RDE~\eqref{eq:Zvertexweightlaw}. 

We then proceed in a similar fashion as with edge weights by applying Kolmogorov's extension Theorem, with the difference that we need to correlate $\hat{Z}$ on a boundary with the weights of the vertices at the boundary, then use the recursion to define it on an entire neighbourhood. 

We recover $Z$ by setting $Z(u,v)=\hat{Z}(u,v)+w(v)$, the remaining geometric considerations are then identical as $Z$ still satisfies equation $\eqref{eq:Z_hrecursion}$. 

\subsection{Maximum subgraph satisfying random capacity constraints}
As noted by previous authors \cite{aldous2000zeta2,salez2011cavity}, it is possible to study a slightly more general type of problem with the approach developed in this work. 

Let $\mathcal{C}$ be a law on $\mathbb{N}$ with finite expectation. Consider a random graph on $G_n$ that converges locally to a UBGW tree $\mathbb{T}$. We decorate $G_n$ by adding random independent decorations on vertices $c(v)$ of law $\mathcal{C}$ that we call the capacity of $v$. The local limit is then the previous tree with additional independent capacities drawn on its vertices. 

The maximum subgraph under capacity constraint is the subgraph $M$ of $G$ maximising
\[ \sum_{e \in M} w(e) \]
subject to
\[ \deg_M(v) \leq c(v),  \,\,\,\forall  v \in V  .\]
In the deterministic case, when $c(v)\equiv 1$, we recover a matching, when $c(v) \equiv k$, a maximal subgraph with degrees less than or equal to $k$.

The variables $Z$ for this problem can be defined as:
\[Z(u,v)=OPT(T^v) - OPT(T^v \text{ where $c(v)$ is decreased by $1$}).   \]
Let $\overset{[k]}{\max}$ be the operator that returns the $k^{th}$ largest value of a set. The recursion becomes: 
\[ Z(u,v)=\max\left(0,\underset{\substack{u' \sim v \\ u' \neq u}}{\overset{[c(v)]}{\max}} \, (w(v,u')-Z(v,u'))\right),  \]
where we set by convention that the $\overset{[0]}{\max}$ of a list is $+ \infty$.
Let $c_l=\mathbb{P}_{C \sim \mathcal{C}}[C=l] $.
This translates into an equation on its cdf $h$ of the form:
\[h(t)=\mathbbm{1}_{t \geq 0} \,  \sum_{l=0}^{\infty}c_l \sum_{p=0}^{l-1} \frac{\left(\mathbb{E}_{W\sim \omega}\left[h(W-t)\right]\right)^p}{p!} \hat{\phi}^{(p)}(1-\mathbb{E}_{W\sim \omega}\left[h(W-t)\right])  \]
Which is still of the form:
\[ h(t)=\mathbbm{1}_{t \geq 0} \Phi(\mathbb{E}_{W \sim \omega}\left[h(W-t)\right])  \]
where $\Phi$ is continuous and decreasing.
The decision rule remains: 
\[(u,v) \in \mathbb{M}_{\mathrm{opt}} \Leftrightarrow Z(u,v)+Z(v,u) < w(u,v).  \]

\bigskip

\printbibliography

\bigskip

{Nathanaël Enriquez: \href{mailto:nathanael:enriquez@universite-paris-saclay.fr}{nathanael.enriquez@universite-paris-saclay.fr}}\\
{Laboratoire de Mathématiques d’Orsay, CNRS, Université Paris-Saclay, 91405, Orsay, France
and DMA, \'Ecole Normale Supérieure – PSL, 45 rue d’Ulm, F-75230 Cedex 5 Paris, France}

\medskip

{Mike Liu: \href{mailto:mike.liu@ensae.fr}{mike.liu@ensae.fr}}\\
{GENES, Fairplay joint team, CREST France}

\medskip

{Laurent Ménard: \href{mailto:laurent.menard@normalesup.org}{laurent.menard@normalesup.org}}\\
{Modal'X, UMR CNRS 9023, UPL, Univ. Paris-Nanterre, F92000 Nanterre, France and ENSAE, Criteo AI Lab \& Fairplay joint team, France.}  

\medskip

{Vianney Perchet: \href{mailto:vianney.perchet@ensae.fr}{vianney.perchet@ensae.fr}}\\
{ENSAE \& Criteo AI Lab (Fairplay team), France}

\end{document}